%% file: main.tex
\documentclass[10pt]{article}

\usepackage{authblk}

\usepackage{amsmath,amssymb}
\usepackage{isotope}
\usepackage{dirtytalk}
\usepackage{amsbsy}
\usepackage{newclude}
\usepackage{nicefrac,xfrac}
\usepackage{bbm}
\usepackage{makecell}
\usepackage[utf8]{inputenc}

\usepackage[utf8]{inputenc}
\usepackage{hyperref}
\usepackage{tikz}
\usepackage{authblk}
\usepackage[maxbibnames=99,
backend=biber,
sorting=nyt
]{biblatex}
\newcommand\norm[1]{\lVert#1\rVert}
\usepackage{accents}
\newlength{\dhatheight}

\usepackage[T1]{fontenc}
\usepackage{amsfonts, amsthm, caption, framed, enumerate, mathtools, fullpage, subcaption, siunitx, tabularx, verbatim, relsize, stackrel, cancel, multicol, media9, xcolor}

\usepackage[normalem]{ulem}

\providecommand{\keywords}[1]
{	
  \textbf{\textit{Keywords:}} #1
}

\providecommand{\MSC}[1]
{	
  \textbf{\textit{MSC2020 subject classification:}} #1
}

\usepackage{booktabs}

\newtheorem{theorem}{Theorem}[section]
\newtheorem{lemma}[theorem]{Lemma}
\newtheorem{example}[theorem]{Example}
\newtheorem{corollary}[theorem]{Corollary}

\newtheorem{proposition}[theorem]{Proposition}
\theoremstyle{definition}
\newtheorem{definition}[theorem]{Definition}
\theoremstyle{remark}
\newtheorem{remark}[theorem]{Remark}

\newtheorem{innercustomass}{Assumption}
\newenvironment{customass}[1]
  {\renewcommand\theinnercustomass{#1}\innercustomass}
  {\endinnercustomass}

\title{Branching Interval Partition Diffusions}
\author{Matthew Buckland}
\affil{Department of Statistics, University of Oxford. Email: matthew.buckland@bnc.ox.ac.uk}
\begin{document}

\maketitle

\begin{abstract}
    We introduce and study branching interval partition diffusions in their natural generality. We let interval widths evolve independently according to a general real-valued diffusion subject only to conditions that ensure finite lifetimes of intervals and allow the continuous generation of new intervals. The latter is governed by the Pitman-Yor excursion measure of the real-valued diffusion and an associated spectrally positive L\'evy process to order both these excursions and their start times. This generalises previous work by Forman, Pal, Rizzolo and Winkel on the self-similar case and gives rise to a new class of general Markovian homogeneous Crump-Mode-Jagers-type branching processes with characteristics varying diffusively during their lifetimes.
\end{abstract}
\small
\keywords{interval partition; excursion theory; branching processes; Ray-Knight theorem; diffusion; L\'evy process}
\newline
\MSC{Primary 60J25; 60J60; 60J80, Secondary 60G51; 60G55}
\normalsize

\include*{introduction}
\include*{preliminaries}
\include*{initialresults}
\include*{continuity}
\include*{appendix}

\vspace{0.5cm}
\renewcommand{\abstractname}{Acknowledgements}
\begin{abstract}
    This research has been supported by the EPSRC Centre for Doctoral Training in Mathematics of Random Systems: Analysis, Modelling and Simulation (EP/S023925/1). The author would also like to thank his supervisor Matthias Winkel for many re-readings of drafts as well as very helpful discussions with the technicalities in this project.
\end{abstract}

\printbibliography

\end{document}

%% file: introduction.tex
\section{Introduction}\label{section introduction}
Continuous-time branching processes with finite birth rates and independent identically distributed (iid) lifetimes of individuals are known as Crump-Mode-Jagers (CMJ) processes, as studied in Jagers \cite[Chapter 6]{jagers1975branching}. For these processes individuals are born from single parents, and during their lives they can in turn give birth to more individuals until they die after some iid positive random times after their births; those with constant birth rates are known as homogeneous and those for which all births produce exactly one offspring are known as binary. Specifically, CMJ processes are the population size processes and so take values in $\mathbb{N}$. In the homogeneous binary cases individuals produce single offspring after iid exponential random times and this repeats until their deaths (or never occurs if they die before their first reproduction event). We call the difference between birth and death the lifetime of the individual and the interval of time between the birth and death the life of the individual. To produce Markovian homogeneous binary CMJ processes we require either the lifetimes to be infinite or for the lifetimes to be exponentially distributed.
\par
We study homogeneous binary branching processes with infinite birth rates where individuals are marked with non-negative characteristics that continuously vary during their lives, and we construct continuous Markovian interval-partition-valued processes by the method shown in Forman, Pal, Rizzolo, and Winkel \cite{construction}. Interval partitions are ordered collections of intervals with summable lengths; we define these below following Aldous \cite[Section 17]{aldous1985} and Pitman \cite[Chapter 4]{pitman2006combinatorial}.
\begin{definition}[Interval partitions]\label{definition interval partition}
For $M>0$, an \emph{interval partition of length $M$} is a set $\beta$ of intervals of $[0,M]$ that are open and disjoint, and which cover $[0,M]$ up to a Lebesgue-null set. We write $\norm{\beta}$ for the length $M$. An \emph{interval partition} is an interval partition of length $M$ for some $M$. We define $\mathcal{I}_H$ to be the space of interval partitions.
\end{definition}
In our model individuals do not live forever and the characteristics are $0$ at the births and deaths and are positive during the rest of the lives. Indeed, the characteristics, as functions of the ages of the individuals they mark, are continuous non-negative excursions away from $0$ of length given by the lifetime of the individual they mark. For this reason we refer to a characteristic at a point in time as the level of health, or just as the health, for the individual it marks at that point in time, and we refer to the continuous excursion of marks over an individual's life as its health excursion.
\par
Lambert \cite{lambert2010contour} and Lambert and Uribe Bravo \cite{lambert2018totally} study spectrally positive L\'evy processes (SPLPs) that are the contour functions of classes of splitting trees and are known as jumping chronological contour processes (JCCPs). For any $x \geq 0$ the jumps across level $x$ represent the individuals alive at that point in time. We construct interval partitions from the health excursions and the SPLP by the method from \cite{construction}. Definition \ref{definition aggregate mass} formally describes this construction and Figure \ref{figure scaffolding and spindles} uses their terminology which we explain in Section \ref{subsection construction}. The result is that at every point in time an interval partition is produced with a one-to-one correspondence of intervals and individuals alive at that point in time, and each interval has a length given by its corresponding individual's health at that point in time. These interval partitions therefore, as functions of the time for the branching process, form an interval-partition-valued process. This paper continues the investigation into these processes, which we refer to as branching interval partition evolutions (branching IP-evolutions), where new intervals are generated between existing intervals. We generalise the class that are known to be diffusions beyond the self-similar cases that have already been established in \cite{construction}. We use a complete and separable metric on interval partitions $d'_H$, from Forman et al. \cite{forman2018interval}, which we define later in Definition \ref{definition d'_H metric}.
\par
When we consider the health excursions rather than just the lifetimes of the individuals, the birth rate is a measure on the space of continuous non-negative excursions away from $0$, and it is this rate that we are referring to when we say the birth rate is constant. We set the birth rate to be the Pitman-Yor excursion measure of a non-negative diffusion, from Pitman and Yor \cite{pitman1982decomposition}, and we explore the range of possible non-negative diffusions for which this is applicable. The Pitman-Yor excursion measure is more general than the standard It\^{o} excursion measure as it includes cases where $0$ is an exit boundary of the diffusion whereas the It\^o excursion measure only includes cases where $0$ is a regular boundary point of the diffusion (we give the relevant details in Section \ref{section main results}). Under certain conditions for summability of the health levels, that we specify in Theorem \ref{theorem aggregate mass summability}, this leads to branching IP-evolutions where the interval lengths all evolve independently with the same continuous Markovian dynamics given by the specified non-negative diffusion. It is these dynamics on the interval lengths which make possible the continuity and Markov property for the branching IP-evolutions. This is also why we choose excursions away from $0$, as this preserves continuity at the points of birth and death of individuals; when an individual is born a new interval grows continuously out of $0$ in the IP-evolution, with the parent's interval on the left, and it dies when this health level returns to $0$.
\par
Excursion measures of diffusions are not finite but $\sigma$-finite, and this is why the processes that we study have an infinite birth rate. However, most contribute small health levels and have small lifetimes in the sense that, for any positive threshold $\varepsilon>0$ there is a finite rate of births with lifetimes or maximum health levels greater than $\varepsilon$. The continuity and the Markov property of the branching IP-evolutions do not follow immediately from the continuity and Markov property of the one-dimensional diffusions, indeed, it is not even the case that we can produce an IP-evolution from the Pitman-Yor excursion measure of any diffusion. This is because we may not be able to construct the necessary SPLP \say{scaffolding} or, even if we can, the health levels may not be summable. In Section \ref{section main results}, Theorem \ref{theorem l\'evy measure condition} gives conditions for the construction of the SPLP, Theorem \ref{theorem aggregate mass summability} conditions for the summability of health levels, and Theorems \ref{theorem diffusion} and \ref{theorem interval partitions regular} give conditions for the existence of continuous Markovian branching IP-evolutions.
\subsection{Construction of the SPLP "scaffolding" and the "skewer process"}\label{subsection construction}
Let $\Lambda$ be a L\'evy measure on $(0,\infty)$ for a recurrent L\'evy process, i.e. a measure such that $\int_{x=0}^\infty (x^2 \wedge x) \Lambda(dx)$ is finite, and let $\mathbf{J}$ be a Poisson random measure (PRM) with intensity $\text{Leb} \otimes \Lambda$ (where \say{Leb} is the Lebesgue measure on $[0,\infty)$). If $\int_{x=0}^\infty x\Lambda(dx)$ is finite then we can construct a recurrent SPLP of bounded variation
\begin{equation*}
    \mathbf{X}_t \coloneqq -t \int_{x=0}^\infty x\Lambda(dx) + \sum_{\substack{(s,j) \text{    atom of   } \mathbf{J}:\\ 0<s<t}} j.
\end{equation*}
Now suppose that $\Lambda$ instead satisfies the condition $\int_{x=0}^\infty (x^2 \wedge x)\Lambda(dx)<\infty=\int_{x=0}^\infty x\Lambda(dx)$. In this case, for any $\varepsilon>0$, we can construct a recurrent SPLP of bounded variation using the measure $\Lambda_\varepsilon(dx)\coloneqq \mathbbm{1}_{x>\varepsilon}\Lambda(dx)$ as we did above and take the limit $\varepsilon \downarrow 0$ is to obtain a recurrent SPLP of unbounded variation, see \eqref{equation scaffolding process}. For our applications the L\'evy measure $\Lambda$ we will use is of the form $\nu(\zeta \in \cdot)$ where $\nu$ is the Pitman-Yor excursion measure for our specified diffusion (a measure on $\mathcal{E}$, the space of continuous positive excursions) and $\zeta$ is a function, $\zeta:\mathcal{E} \rightarrow [0,\infty)$, which gives the length of excursions which we call the lifetimes. 
\par
More formally, continuous positive excursions are functions $f:\mathbb{R} \rightarrow \mathbb{R}$ for which $f(y)=0$ for all $y \leq 0$, and there exists some $x>0$ such that $f(y)=0$ for all $y > x$, and finally for $y \in (0,x)$ that $f(y)>0$ and $f$ is continuous on $[0,x]$; the lifetime $\zeta$ of $f$ is then $\zeta(f)=x$. We use the lifetimes as the jump sizes for the SPLP. Suppose $\int_{x=0}^\infty(x^2 \wedge x) \nu(\zeta \in dx)<\infty$, and let $\mathbf{N}'$ be a PRM with intensity $\text{Leb} \otimes \nu$. We define the recurrent SPLP $\mathbf{X}= (\mathbf{X}_t, t \geq 0)$ by the almost sure limits
\begin{equation}\label{equation scaffolding process}
    \mathbf{X}_t \coloneqq  \lim_{\varepsilon \downarrow 0}\Bigg(-t\int_{x=\varepsilon}^\infty x \nu(\zeta \in dx) + \sum_{\substack{(s,f) \text{    atom of   } \mathbf{N}'\\ 0<s<t}}\zeta(f)\mathbbm{1}_{\zeta(f)>\varepsilon}\Bigg).
\end{equation}
The PRM $\mathbf{N}'$ takes values in a certain class of $\sigma$-finite point measures on $[0,\infty) \times \mathcal{E}$ which we denote by $\mathcal{N}$, and the SPLP $\mathbf{X}$ is a c\`adl\`ag process and we denote the set of c\`adl\`ag paths by $\mathcal{D}$. The class $\mathcal{N}$ consists of all the point measures for which the limit in \eqref{equation scaffolding process} exists for all $t>0$. The L\'evy measure for $\mathbf{X}$ is $\nu(\zeta \in \cdot)$ which we refer to as the lifetime under the Pitman-Yor excursion measure. More generally, for a suitable point measure $N \in \mathcal{N}$, let $\xi(N)$ be the c\`adl\`ag function, $(X_t, t \geq 0)$, obtained when replacing $\mathbf{N}'$ by $N$ in \eqref{equation scaffolding process}.
\begin{figure}
    \centering
    \includegraphics[width=\textwidth]{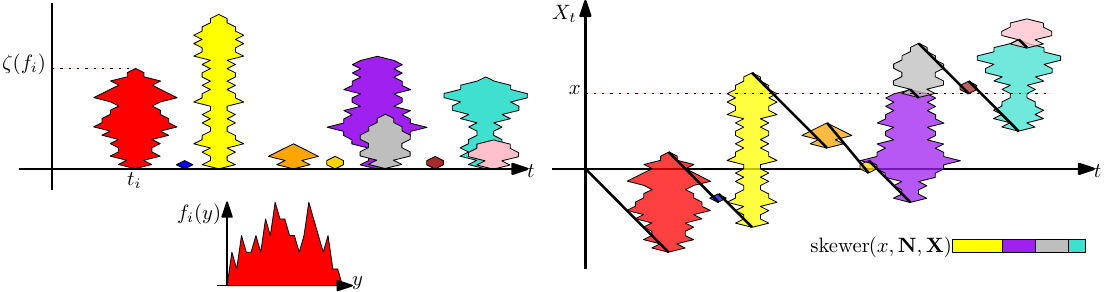}
    \caption{Left: A PRM of health excursions represented with time vertically moving up page and reflected to produce spindle-like shapes, with the excursion of a single spindle underneath. Right: the scaffolding process (black) constructed from the PRM, along with the interval partition generated by the skewer process at a fixed level $x$.}
    \label{figure scaffolding and spindles}
\end{figure}
\par
We call $\mathbf{X}$ the scaffolding and $\mathbf{N}'$ the point measure of spindles. The reason for these names can be seen in the simplified version $(N,X)$ in Figure \ref{figure scaffolding and spindles}. In the SPLP, there is a one-to-one correspondence between individuals and jumps. Moreover, the jump heights are equal to the individuals' lifetimes, and levels over which each jump crosses are the times for which the corresponding individual is alive in our model. Therefore, in the figure, the health excursion can be placed vertically in each jump. The excursions themselves look like spindles as we have width equal to health and they are symmetric for better visual representation on the diagram. We next define the skewer function.
\begin{definition}[Skewer function]\label{definition aggregate mass}
Let $I \subseteq J$ be index sets and $N\coloneqq  \sum_{i \in I} \delta_{(r_i, f_i)}$, $N'\coloneqq  \sum_{i \in J} \delta_{(r_i, f_i)} \in \mathcal{N}$, and let $X=\xi(N')$. We define the aggregate health process $((M^y_{N,X})(t), t \in [0,\infty))$ by
\begin{equation*}
    M^y_{N,X}(t) \coloneqq  \sum_{\substack{(r,f)\text{   atom of   } N: \\ 0 \leq r < t}} f(y-X_{r-}).
\end{equation*}
We also define the total health process at $t \in [0,\infty)$ as $(M^y_{N,X}(t), y \in \mathbb{R})$. When $M^y_{N,X}(\infty)<\infty$ we define the skewer for $(N,X)$ by $\text{skewer}(y,N,X)\coloneqq  \{ (M^y_{N,X}(t-), M^y_{N,X}(t)): t \in [0,\infty), M^y_{N,X}(t-) < M^y_{N,X}(t) \}$. Note that in this case $\text{skewer}(y,N,X)$ is an interval partition of length $M^y_{N,X}(\infty)$. We define $\text{skewer}(y,N,X)$ to be $\emptyset$, the empty interval partition, if $M^y_{N,X}(\infty)$ is infinite. For PRMs $\mathbf{N} \leq \mathbf{N}'$ and scaffolding $\mathbf{X} = \xi(\mathbf{N}')$ we define the scaffolding-and-spindles construction of $(\mathbf{N},\mathbf{X})$ to be the IP-evolution given by $\overline{\text{skewer}}(\mathbf{N},\mathbf{X})\coloneqq (\text{skewer}(y,\mathbf{N},\mathbf{X}), y \geq 0)$. Unless stated otherwise, we take $I=\{j \in J:r_j \leq T\}$ where $T=T(X)$ is such that $T(\mathbf{X})$ is an a.s finite stopping time and in this case we refer to $\mathbf{N}$ as a stopped PRM.
\end{definition} 
\subsection{Main results}\label{section main results}
We define one-dimensional diffusions with the framework given by Revuz and Yor \cite[Chapter VII.3]{revuz2013continuous}. In our setting we consider a Markov process $\mathbf{Y}$ whose state space $E$ is of the form $E=[0,c]$ where $c$ is finite or of the form $E=[0,c)$ where $c$ is either finite or infinity. We assume that the paths of $\mathbf{Y}$ are continuous on $E$ and that $\mathbf{Y}$ has the strong Markov property. We further assume that our diffusion is regular on $(0,c)$, which means that for all points $x$, $y \in (0,c)$ the probability of a diffusion reaching point $y$ having started from point $x$ is positive, i.e. $\mathbb{P}_x(T_y<\infty)>0$ where $T_y$ is the hitting time of the point $y$ by the diffusion. With these assumptions, there exists a strictly increasing function $s:E \rightarrow [0,\infty)$ such that for any $a$, $w$, $x$ in $E$ with $0 \leq a<x<w\leq c$ we have 
\begin{equation}\label{equation scale function property}
    \mathbb{P}_x(T_w<T_a) = \frac{s(x)-s(a)}{s(w)-s(a)},
\end{equation}
where the function $s$ is called the scale function of $\mathbf{Y}$, and we also define for $v \in (a,w)$
\begin{equation}\label{equation greens function original}
    G_{(a,w)}(x,v) =\frac{s(x)-s(a)}{s(w)-s(a)} \left(s(w)-s(v)\right) \mathbbm{1}_{x<v<w} + \frac{s(w)-s(x)}{s(w)-s(a)} \left(s(v)-s(a)\right) \mathbbm{1}_{a<v<x}
\end{equation}
which we call the Green's function for the diffusion on $(a,w)$. Then there exists a unique measure $M$ on $(0,c)$, which is known as the speed measure, such that for any $0 \leq a<x<w\leq c$ we have
\begin{equation}\label{equation expected hitting time}
    \mathbb{E}_x(T_a \wedge T_w) = \int_{v=a}^w G_{(a,w)}(x,v) M(dv).
\end{equation}
For certain diffusions there exist limits which are known as the infinitesimal drift $(\mu(x), x \in (0,c))$ and the infinitesimal variance $(\sigma^2(x), x \in (0,c))$ which are given by
\begin{equation}\label{equation drift and variance}
    \mu(x) = \lim_{h \downarrow 0} \frac{1}{h} \mathbb{E}_x(Y_h - x);  \hspace{0.5cm} \sigma^2(x) = \lim_{h \downarrow 0} \frac{1}{h} \mathbb{E}_x((Y_h - x)^2).
\end{equation}
We now state assumptions on the diffusions that we consider in this paper.
\begin{customass}{A}\label{assumption a}
A diffusion on state space $E$ satisfies \emph{Assumption \ref{assumption a}} if:
\begin{enumerate}
    \item We have that $\mathbb{P}_x(T_0<\infty)>0$ for some $x \in (0,c)$ which is equivalent to $\int_{y=0}^x s(y) M(dy) < \infty$ (equivalently for all $x \in (0,c)$ in both cases).
    \item For all $0<y<x<c$ we have that $\mathbb{E}_x(T_y)$ is finite. This is satisfied if and only if $M(E \backslash [0,x])$ is finite for all $x \in (0,c)$.
    \item For all $x \in (0,c)$, the limits $\mu(x)$ and $\sigma^2(x)$ in \eqref{equation drift and variance}  exist and $\mu:(0,c) \rightarrow \mathbb{R}$ and $\sigma^2:(0,c) \rightarrow (0,\infty)$ are locally bounded Borel-measurable functions.
    \item The functions $\mu$ and $\sigma^2$ are continuous on $(0,b]$ for some $b \in E \backslash \{0\}$.
\end{enumerate}
We label these conditions $A1$, $A2$, $A3$, $A4$.
\end{customass}
The first condition is equivalent to $0$ being an attainable boundary point for the diffusion, i.e. a regular or exit boundary point according to the Feller boundary point classification, see \cite[Table 15.6.1]{karlin1981second}. The boundaries of the diffusion are the points $0$ and $c$. In words, a regular boundary point is one where you can restart the diffusion from the boundary once the diffusion first hits the boundary, i.e. $\mathbf{P}_0(T_x<\infty)>0$; a typical example is the boundary point $0$ for a Brownian motion in $[0,\infty)$ that is reflected at $0$. An exit boundary point is where this behaviour does not occur, and instead, you have that the diffusion starting in the interior of the state space can reach the point $0$ with positive probability, but the diffusion cannot reenter the interior after hitting $0$. We note that the boundary point $c$ is regular exactly when $s(c)\coloneqq \lim_{x \uparrow c}s(x)$ is finite equivalently when the state space $E$ is equal to $[0,c]$ and not $[0,c)$. In later sections this case sometimes needs special consideration, and we may identify it by any of the equivalent conditions just stated depending on which is most relevant to the section. Furthermore, as some of the measures that we use in this paper are not necessarily atom-free we take care when we define the limits of integration. We write $\int_{w=a}^b$ for an integral over $(a,b)$ and we write $\int_{w \in (a,b]}$ to integrate over $(a,b]$ (and then similarly to integrate over $[a,b)$ or $[a,b]$).
\begin{theorem}[L\'evy measure]\label{theorem l\'evy measure condition}
Consider a diffusion that satisfies $A1$, $A2$ of Assumption \ref{assumption a}. Then it has a Pitman-Yor excursion measure $\nu$ and we have that
\begin{equation}\label{equation levy measure bounded variation}
    \int_{x=0}^\infty x \nu(\zeta \in dx) < \infty
\end{equation}
if and only if
\begin{equation}\label{equation 0 is regular}
    M(E) < \infty,
\end{equation}
which is equivalent to $0$ being a regular boundary point for the diffusion. In this setting the lifetime under $\nu$ is a L\'evy measure for a recurrent SPLP of bounded variation.
\par
On the other hand, if we have that: the diffusion satisfies all of Assumption \ref{assumption a}; the condition in \eqref{equation 0 is regular} does not hold (i.e. $0$ is an exit boundary of the diffusion); and we have
\begin{equation}\label{equation x^2 bound}
    \int_{v=0}^b \int_{y=0}^v s(y) M(dy) M(dv) < \infty  \text{   where   } b \in (0,c) \text{   is as in } A4
\end{equation}
then we have
\begin{equation}\label{equation levy measure unbounded variation}
    \int_{x=0}^\infty (x^2 \wedge x) \nu(\zeta \in dx) < \infty = \int_{x=0}^\infty x \nu(\zeta \in dx).
\end{equation}
In this setting the lifetime under $\nu$ is a L\'evy measure for a recurrent SPLP of unbounded variation.
\end{theorem}
We believe that, under $A1$, $A2$, the condition in \eqref{equation x^2 bound} is necessary as well as sufficient for the lifetime under $\nu$ to be a L\'evy measure for a recurrent process (i.e. we also believe that the regularity conditions $A3$ and $A4$ are not necessary), but proving this is beyond the scope of this paper. The L\'evy measure condition of Theorem \ref{theorem l\'evy measure condition} means that we can construct a scaffolding-and-spindles pair $(\mathbf{N},\mathbf{X})$. Later, in Proposition \ref{proposition starting interval partition}, we define and construct $(\mathbf{N}_\beta,\mathbf{X}_\beta)$, for $\beta \in \mathcal{I}_H$. This is a scaffolding and spindles pair for which we prove that $\overline{\text{skewer}}(\mathbf{N}_\beta, \mathbf{X}_\beta)$ is the modification of an IP-evolution such that $\text{skewer}(\mathbf{N}_\beta, \mathbf{X}_\beta,0)=\beta$ a.s.
\par
Throughout the rest of the paper we focus on diffusions that satisfy $A1$ and $A2$ and transformations $g:E\rightarrow g(E)$ that are strictly increasing and continuous with $g(0)=0$. We define $g(c)=\lim_{x \uparrow c}g(x) \leq \infty$ and therefore $g(E)=[0,g(c))$ or $g(E)=[0,g(c)]$. Let $\mathbf{Y}=(\mathbf{Y}_t, t \geq 0)$ be a diffusion in state space $E$ that satisfies $A1$ and $A2$ and let $\mathbf{Z} =(\mathbf{Z}_t, t \geq 0)$ be the process defined by $\mathbf{Z}_t=g(\mathbf{Y}_t)$ for all $t\geq 0$. Then firstly, $\mathbf{Z}$ is a diffusion on $g(E)$ that satisfies $A1$ and $A2$ and which has Pitman-Yor excursion measure $\nu_{\mathbf{Z}}$ such that $\nu_{\mathbf{Z}}(\zeta \in \cdot)=\nu_{\mathbf{Y}}(\zeta \in \cdot)$ where $\nu_{\mathbf{Y}}$ is the Pitman-Yor excursion measure for $\mathbf{Y}$. We refer to $\mathbf{Y}$ as the untransformed diffusion throughout the paper to avoid confusion between $\mathbf{Y}$ and $\mathbf{Z}$.
\begin{theorem}[Aggregate health summability]\label{theorem aggregate mass summability}
In the setting described above, consider a diffusion $\mathbf{Z}$ constructed from an untransformed diffusion $\mathbf{Y}$ along with a spatial transformation $g$. Suppose that $\nu=\nu_\mathbf{Z}$ (equivalently $\nu=\nu_\mathbf{Y}$) satisfies \eqref{equation levy measure bounded variation} or \eqref{equation levy measure unbounded variation} and let $(\mathbf{N},\mathbf{X})$ be scaffolding-and-spindles pair where $\mathbf{N}$ has intensity $\text{Leb} \otimes \nu_{\mathbf{Z}}$ stopped at an $(0,\infty)$-valued random variable $T$.
\par
If $\nu=\nu_\mathbf{Z}$ satisfies \eqref{equation levy measure bounded variation}, then the aggregate health process $M^0_{\mathbf{N},\mathbf{X}}(T)$ is finite a.s. If $\nu=\nu_{\mathbf{Z}}$ satisfies \eqref{equation levy measure unbounded variation}, then $M^0_{\mathbf{N},\mathbf{X}}(T)$ is finite a.s. if and only if
\begin{equation}\label{equation IP existence condition}
     \int_{y=0}^b g(y) M_\mathbf{Y}(dy) < \infty \text{   for some (equivalently all)   } b \in (0,c),
\end{equation}
or equivalently that $\int_{z=0}^b z M_\mathbf{Z}(dz)$ is finite for some (or all) $b \in (0,g(c))$. Furthermore in this case we have that $M^0_{\mathbf{N},\mathbf{X}}(T)$ is infinite a.s. when \eqref{equation IP existence condition} does not hold.
\end{theorem}
We now give a stronger set of assumptions in Assumption \ref{assumption b}.
\begin{customass}{B}\label{assumption b}
A diffusion $\mathbf{Z}$ satisfies \emph{Assumption \ref{assumption b}} if it is of the form $\mathbf{Z}=g(\mathbf{Y})$, where $\mathbf{Y}$ is a diffusion that satisfies $A1$, $A2$, $A3$ of Assumption \ref{assumption a} with state space $E \subseteq [0,\infty)$ of the form $E=[0,c]$ or $E=[0,c)$ for some $c \leq \infty$ and furthermore is such that:
\begin{enumerate}
    \item the Pitman-Yor excursion measure $\nu_\mathbf{Y}$ satisfies \eqref{equation levy measure unbounded variation} where $\nu=\nu_\mathbf{Y}$ (equivalently $\nu=\nu_\mathbf{Z}$);
    \item $\sigma_{\mathbf{Y}}^2(x)=4x$;
\end{enumerate}
there exists $\varepsilon_0 \in (0,c)$ such that
\begin{enumerate}
    \setcounter{enumi}{2}
    \item for some $0<\alpha^- \leq \alpha^+<1$ we have
        \begin{equation*}
            -2\alpha^+ \leq \mu_{\mathbf{Y}}(x) \leq -2\alpha^- \text{   for all   } x \in (0,\varepsilon_0]
        \end{equation*}
\end{enumerate}
(where the subscript $\mathbf{Y}$ in $\mu$ and $\sigma^2$ denote that they are the infinitesimal drift and variance for $\mathbf{Y}$), and $g:E \rightarrow g(E)$ is a strictly increasing continuous function with $g(0)=0$, $g(E) \subseteq [0,\infty)$ and
\begin{enumerate}
    \setcounter{enumi}{3}
        \item there exists $q_0>\alpha^+$ such that the transformation $g$ is $q$-H\"older continuous on $[0,\varepsilon_0)$ for all $q<q_0$, i.e. for all $q<q_0$ there exists $C_q \geq 0$ such that for all $x$, $y \in [0,\varepsilon_0)$ we have $|g(x)-g(y)|\leq C_q|x-y|^q$. 
\end{enumerate}
Finally, we say $\mathbf{Z}=g(\mathbf{Y})$ satisfies Strong Assumption \ref{assumption b} if as well as the above we have that
\begin{enumerate}
    \setcounter{enumi}{4} 
    \item the transformation $g$ satisfies $0<\liminf_{x \rightarrow 0} \frac{g(x)}{x} \leq \infty$ or $0 \leq \limsup_{x \rightarrow 0} \frac{g(x)}{x} < \infty$.
\end{enumerate}
\end{customass}
We believe that $B1$ follows from $B2$ and $B3$ but we could not show this. We do show however in Proposition \ref{corollary levy measure assumption b} that if the untransformed diffusion additionally satisfies $A4$ then we have that $B2$ and $B3$ implies $B1$. We give intuition behind condition $B5$ with an example of a transformation $g$ that does not satisfy the condition in Example \ref{example g counterexample}. 
\begin{theorem}[Interval Partition Diffusions when $0$ is an exit point]\label{theorem diffusion}
Consider a diffusion $\mathbf{Z}$ that satisfies Assumption \ref{assumption b} and let $\nu$ be its Pitman-Yor excursion measure. Let $(\mathbf{N},\mathbf{X})$ be a scaffolding-and-spindles pair where $\mathbf{N}$ has intensity $\text{Leb} \otimes \nu$ stopped at an $(0,\infty)$-valued random variable $T$. Then there exists a modification of $\overline{\text{skewer}}(\mathbf{N},\mathbf{X})$ which is a continuous branching IP-evolution. Furthermore the process $\mathbf{X}$ has a local time that is continuous in time and level. Let $\tau^0$ be the right-continuous inverse of the local time of $\mathbf{X}$ at level $0$. Then for $T=\tau^0(s)$ for some $s >0$ the continuous branching IP-evolution defined above is a simple Markov process in an appropriate filtered probability space and if furthermore $\mathbf{Z}$ satisfies Strong Assumption \ref{assumption b} then this IP-evolution is strongly Markovian.
\end{theorem}
We prove Theorem \ref{theorem diffusion} in Section \ref{section continuity} and we establish continuity of branching IP-evolutions from a sufficiently large class of possible initial states in Section \ref{section continuity initial partition}. We now state Assumption \ref{assumption c} which considers a disjoint set of diffusions to Assumption \ref{assumption b}. 
\begin{customass}{C}\label{assumption c}
    A diffusion $\mathbf{Z}$ satisfies \emph{Assumption \ref{assumption c}} if it is of the form $\mathbf{Z}=g(\mathbf{Y})$, where $\mathbf{Y}$ where $\mathbf{Y}$ is a diffusion that satisfies $A1$, $A2$, $A3$ of Assumption \ref{assumption a} with state space $E \subseteq [0,\infty)$ of the form $E=[0,c]$ or $E=[0,c)$ for some $c \leq \infty$ and furthermore is such that:
    \begin{enumerate}
        \item $\sigma_{\mathbf{Y}}^2(x)=4x$;
    \end{enumerate}
    there exists $\varepsilon_1 \in (0,c)$ such that
    \begin{enumerate}
        \setcounter{enumi}{1}
        \item for some $0<\beta^- \leq \beta^+<1$ we have
        \begin{equation*}
            2\beta^- \leq \mu_{\mathbf{Y}}(x) \leq 2\beta^+ \text{   for all   } x \in (0,\varepsilon_1];
        \end{equation*}
        \item and $g:E \rightarrow g(E)$ is a strictly increasing continuous function with $g(0)=0$, $g(E) \subseteq [0,\infty)$ and
        \begin{equation*}
            \int_{y=0}^{\varepsilon_1} \frac{g(y^{1/(1-\beta^-)})}{y^2} dy < \infty.
        \end{equation*}
    \end{enumerate}
\end{customass}
Theorem \ref{theorem interval partitions regular} states our second main result.
\begin{theorem}[Interval Partition Diffusions when $0$ is a regular point]\label{theorem interval partitions regular}
Consider a diffusion $\mathbf{Z}$ under Assumption \ref{assumption c} with excursion measure $\nu$ and let $(\mathbf{N},\mathbf{X})$ be scaffolding-and-spindles pair where $\mathbf{N}$ has intensity $\text{Leb} \otimes \nu$ stopped at a $(0,\infty)$-valued random variable $T$. Then $\overline{\text{skewer}}(\mathbf{N},\mathbf{X})$ is a continuous branching IP-evolution. Furthermore we have that $\mathbf{X}$ is a L\'evy process of bounded variation and when $T$ denotes the $n^{\text{th}}$ hitting time for $\mathbf{X}$ of level $0$ for some $n \in \mathbb{N}$ then $\overline{\text{skewer}}(\mathbf{N},\mathbf{X})$ is a simple Markov process in an appropriate filtered probability space. 
\end{theorem}
We believe in this setting that $\overline{\text{skewer}}(\mathbf{N},\mathbf{X})$ is a continuous strongly Markovian process but showing this is beyond the scope of the paper. We do not believe that Theorems \ref{theorem diffusion} and \ref{theorem interval partitions regular} cover all the one-dimensional diffusions from which it is possible to construct continuous Markovian branching IP-evolutions by the construction of scaffolding and spindles described in Section \ref{subsection construction}, though we believe that all such diffusions must satisfy either \eqref{equation 0 is regular} or both \eqref{equation x^2 bound} and \eqref{equation IP existence condition}. 
\subsection{Examples}
In this section we list a few examples of diffusions that satisfy the theorems in the previous section. The diffusions in our selected examples satisfy $A3$ and we begin with some consequences of the existence of the infinitesimal drift and infinitesimal variance. We call an IP-evolution with intervals that evolve as independent diffusions with infinitesimal drift $\mu$ and infinitesimal variance $\sigma^2$ a $(\mu, \sigma^2)$-IP-evolution. 
\par
Consider a diffusion that has well-defined infinitesimal drift and variance, see \eqref{equation drift and variance}, and has state space $E$ of the form $[0,c)$ or $[0,c]$ for some $c \leq \infty$. In this case, for any $b \in (0,c)$, we can define a scale function $s$ with derivative $s'$ and its speed measure $M$ with density $m$ by, for $x \in (0,c)$,
\begin{equation}\label{equation scale and speed measure}
    s'(x) = \exp\left(-\int_{z=b}^x \frac{2 \mu (z)}{\sigma^2(z)}dz\right); \hspace{0.25cm}  \hspace{0.25cm} m(x) = \frac{2}{\sigma^2(x)s'(x)} = \frac{2}{\sigma^2(x)} \exp \left(\int_{z=b}^x \frac{2 \mu(z)}{\sigma^2(z)}dz\right);
\end{equation}
$s(x)=\int_{y=0}^x s'(y) dy$; $M((l,r))=\int_{y=l}^r m(y)dy$. We were free to choose the lower limit of $b$ here and in Section \ref{section model assumptions} we take $b$ to be $\varepsilon_0$ when under Assumption \ref{assumption b} and to be $\varepsilon_1$ when under Assumption \ref{assumption c}. Different values in the interior only change $s'$, $s$, $m$ by the same constant scale factor and such changes do not affect any of the conditions in our results.
\par
If $g$ is a twice differentiable function with continuous second derivative, the parameters $\mu_{\mathbf{Z}}$ and $\sigma^2_{\mathbf{Z}}$ are defined and we have, for $x \in (0,c)$, that
\begin{equation}\label{equation infinitesimal transformations}
    \mu_{\mathbf{Z}}(g(x))= \frac{1}{2}\sigma^2_{\mathbf{Y}}(x)g''(x) + \mu_{\mathbf{Y}}(x)g'(x); \hspace{0.5cm} \sigma^2_{\mathbf{Z}}(g(x))=\sigma^2_{\mathbf{Y}}(x)(g'(x))^2,
\end{equation}
We further have that
\begin{equation*}
    s_{\mathbf{Z}}(g(x)) = s_{\mathbf{Y}}(x); \hspace{0.5cm} s'_{\mathbf{Z}}(g(x)) = \frac{1}{g'(x)} s'_{\mathbf{Y}}(x); \hspace{0.5cm} m_{\mathbf{Z}}(g(x)) = \frac{1}{g'(x)} m_{\mathbf{Y}}(x).
\end{equation*}
This means that the untransformed diffusion $\mathbf{Y}$ satisfies $A1$, $A2$, $A3$ if and only if the diffusion $\mathbf{Z}$ satisfies $A1$, $A2$, $A3$. Similarly, the results in Theorem \ref{theorem l\'evy measure condition} hold for the diffusion $\mathbf{Z}$ if and only if they hold for the diffusion $\mathbf{Y}$.
\par
This first example concerns self-similar diffusions. Our main reference in this work, \cite{construction}, determined all the self-similar diffusions for which it is possible to use the scaffolding and spindles construction for branching interval partition diffusions.
\begin{example}[Self-similar diffusions]
    This case of branching interval partition diffusions has already been established in \cite[Theorems 1.3, 1.4]{construction}. Let $k>0$, $0<\alpha<1,q$ and consider an untransformed diffusion $\mathbf{Y}$ and a transformation $g$ with
    \begin{equation*}
        \mu_{\mathbf{Y}}(y)=-2\alpha; \hspace{0.5cm} \sigma^2_{\mathbf{Y}}(y)=4y; \hspace{0.5cm} g(y)=ky^q.
    \end{equation*}
    Here $\mathbf{Y}$ is a squared Bessel process with dimension $-2\alpha$. Let $\mathbf{Z}$ be the diffusion $\mathbf{Z}=g(\mathbf{Y})$. Then $\mathbf{Z}$ satisfies Assumption \ref{assumption b} and there exist continuous strong-Markov $(\mu_{\mathbf{Z}}, \sigma^2_\mathbf{Z})$-IP-evolutions. These are all the self-similar diffusions for which the construction of branching IP-evolutions that we use in this paper is possible, see \cite[Theorem 5.1]{lamperti1972semi} for a parametrisation of self-similar Markov processes. 
\end{example}
\par
This second example concerns Wright-Fisher diffusion, which are well-known across the literature, and we quote the relevant information about them from \cite[Example 15.6.8]{karlin1981second}. The Wright-Fisher model is a gene frequency model with mutation parameters that describes a population with two different types, say type $1$ and type $2$. The diffusion measures the proportion of the population which is type $1$, and therefore takes values in $[0,1]$. There are two parameters, $\gamma_1$, $\gamma_2 \geq 0$ which determine the rate of mutation from one type to another (i.e. the parameter $\gamma_1$ determines the rates of mutation from type $1$ to type $2$).  They have infinitesimal drift and infinitesimal variance given by
    \begin{equation}
        \mu_{\mathbf{Z}}(z) = \gamma_2 (1-z) - \gamma_1 z, \hspace{0.5cm} \sigma_{\mathbf{Z}}^2(z) =  z(1-z).
    \end{equation}
\begin{example}[Wright-Fisher gene frequency model involving mutation parameters]\label{example wright fisher}
    The Wright-Fisher diffusion with mutation parameters $\gamma_1 \geq 0$ and $0<\gamma_2<1/2$ for interval partition diffusions can be used to generate continuous simple Markov branching $(\mu_{\mathbf{Z}}, \sigma^2_\mathbf{Z})$-IP-evolutions.
    \par
    We prove the existence of IP-evolutions in this case in Proposition \ref{proposition wright fisher} below.
\end{example}
\begin{proposition}\label{proposition wright fisher}
    The Wright-Fisher diffusions as defined in Example \ref{example wright fisher} satisfy Assumption \ref{assumption c} and therefore can be used to generate continuous simple Markov branching IP-evolutions.
    \begin{proof}
    Let $\mathbf{Y}$ be defined on $[0,\pi^2]$ with 
    \begin{equation*}
        \mu_{\mathbf{Y}}(y) = 1 +\sqrt{y} \Big(4\gamma_2\sec\Big(\frac{\sqrt{y}}{2}\Big) - \cot(\sqrt{y}) - 2(\gamma_1+\gamma_2)\tan\Big(\frac{\sqrt{y}}{2}\Big)\Big); \hspace{0.5cm} \sigma_\mathbf{Y}^2(y) = 4y.
    \end{equation*}
    With the transformation $g(y)=\sin^2\big(\frac{\sqrt{y}}{2}\big)$ we have the Wright-Fisher diffusion as described above; this can be checked with the formulae in \eqref{equation infinitesimal transformations}. We now verify the diffusion satisfies Assumption \ref{assumption c} and therefore, by Theorem \ref{theorem interval partitions regular}, there exist continuous simple Markov $(\mu_{\mathbf{Z}},\sigma_{\mathbf{Z}}^2)$-IP evolutions. It is equivalent to check that $\mathbf{Z}$ satisfies $A1$, $A2$, $A3$ rather than check that diffusion $\mathbf{Y}$ does. The diffusion $\mathbf{Z}$ has scale function and speed measure with the following asymptotic properties:
    \begin{equation*}
        s'(z)= O(z^{-2\gamma_2}), \hspace{0.25cm} m(z)= O(z^{2\gamma_2-1}) \hspace{0.2cm} \text{   as    } z \rightarrow 0; \hspace{0.25cm} s'(z) = O((1-z)^{-2\gamma_1}), \hspace{0.2cm} m(z) = O((1-z)^{1-2\gamma_1}) \text{   as   } z \rightarrow 1.
    \end{equation*}
    Therefore
    \begin{equation*}
        \int_{z=0}^b s(z)m(z)dz = O\Big(\int_{z=0}^b z^{1-2\gamma_2} z^{2\gamma_2-1} dz\Big) = O(b) \hspace{0.25cm} \text{   as    } b \rightarrow 0.
    \end{equation*}
    Similarly, the condition $\int^{1}_{z=b}m_{\mathbf{Z}}(z)dz<\infty$ can be checked.
    \par
    The drift satisfies $\mu_{\mathbf{Y}}(y) = 4\gamma_2 + o(1)$ near $0$ and therefore there exists $\varepsilon>0$ such that $2\gamma_2 \leq \mu_{\mathbf{Y}}(y) \leq 2\gamma_2 + 1$ for $y \in (0,\varepsilon)$. Furthermore, $g(y)=O(y)$ near $0$ and so
    \begin{equation*}
        \int_{y=0}^b \frac{g(y^{1/(1-\gamma_2)})}{y^2} dy = O\Big(\int_{y=0}^b y^{-2+1/(1-\gamma_2)} dy\Big) = O(b^{\gamma_2/(1-\gamma_2)}) \hspace{0.25cm} \text{   as    } b \rightarrow 0
    \end{equation*}
    and this is again finite.
    \end{proof}
\end{proposition}
We believe that Wright-Fisher diffusions with $\gamma_2=0$ should also lead to continuous Markov branching IP-evolutions, but our methods were not powerful enough to prove this. On the other hand, the requirement that $\gamma_1>0$ is necessary here to satisfy $A2$.
\par
The next example concerns Cox-Ingersoll-Ross processes, Hutzenthaler \cite{cox-ingersoll-ross}, with financial and biological applications. This family of diffusions is an alternative gene frequency model with two mutation parameters in a scaling where the mutation likelihood is larger in magnitude in the scaling limit, and in this context it is known as Laguerre diffusions after the Laguerre orthogonal polynomials involved, see \cite[Example 15.2.F(d)]{karlin1981second} for details. It also appears as a population growth model in the same reference, see \cite[Example 15.13.C]{karlin1981second}. Finally, they are referred to as Cox-Ingersoll-Ross processes in \cite[Equation (17)]{cox-ingersoll-ross} and have financial applications.
\begin{example}[Cox-Ingersoll-Ross processes/Laguerre Diffusions]
    For $a>0$ and $b$, $c \in \mathbb{R}$ the infinitesimal drift and infinitesimal variance is given by
    \begin{equation*}
        \mu_{\mathbf{Z}}(z) = bx + c; \hspace{0.5cm} \sigma^2_\mathbf{Z}(z) = az.
    \end{equation*}
    The diffusion can be constructed with the transformation $g(y)=(a/4)y$ and the untransformed diffusion with $\mu_{\mathbf{Y}}(y) = by + 4c/a$ and $\sigma^2_{\mathbf{Y}}(y)=4y$. This satisfies Assumption \ref{assumption c} when $0<2c<a$ and $b<0$ and therefore in these cases, by Theorem \ref{theorem interval partitions regular}, there exist continuous simple Markov $(\mu_{\mathbf{Z}},\sigma_{\mathbf{Z}}^2)$-IP-evolutions. 
\end{example}
This final example is something of a middle case between Assumptions \ref{assumption b} and \ref{assumption c} for which we have some results but no proof for the existence of continuous Markovian IP-evolutions.
\begin{example}\label{example squared bessel dim 0}
    Consider diffusions $\mathbf{Z}$ of the form $\mathbf{Z}=g(\mathbf{Y})$, where $\mathbf{Y}$ is a diffusion that satisfies $A1$, $A2$, $A3$ with state space $E \subseteq [0,\infty)$ of the form $E=[0,c]$ or $E=[0,c)$ for some $c \leq \infty$ and furthermore is such that: 
    \begin{enumerate}
        \item $\sigma^2_{\mathbf{Y}}(x)=4x$;
    \end{enumerate}
    there exists $\varepsilon_2 \in (0,c)$ such that
    \begin{enumerate}
        \setcounter{enumi}{1}
        \item $\mu_{\mathbf{Y}}(x)=0$ for all $x \in (0,\varepsilon_2]$. 
    \end{enumerate}
    (where the subscript $\mathbf{Y}$ in $\mu$ and $\sigma^2$ denote that they are the infinitesimal drift and variance for $\mathbf{Y}$), and $g:E \rightarrow g(E)$ is a strictly increasing continuous function with $g(0)=0$, $g(E)\subseteq [0,\infty)$ and
    \begin{enumerate}
        \setcounter{enumi}{2}
        \item there exists $q'_0>0$ such that the transformation $g$ is $q$-H\"older continuous on $[0,\varepsilon_0)$ for all $q<q'_0$, i.e. for all $q<q'_0$ there exists $C'_q \geq 0$ such that for all $x$, $y \in [0,\varepsilon_2)$ we have $|g(x)-g(y)|\leq C'_q|x-y|^q$
    \end{enumerate}
    Note that $\mathbf{Y}$ satisfies $A4$ with $b=\varepsilon_2$. We define a scale function $s_{\mathbf{Y}}$ a speed measure $m_{\mathbf{Y}}$ by setting $b=\varepsilon_2$ in \eqref{equation scale and speed measure} so that
    \begin{equation}\label{equation speed scale squared bessel dim 0}
        s'_{\mathbf{Y}}(x) = 1, \hspace{0.5cm} s_{\mathbf{Y}}(x) = x, \hspace{0.5cm} m_{\mathbf{Y}}(x) = \frac{1}{2x} \hspace{0.5cm} \text{ for } x \in (0,\varepsilon_2).
    \end{equation}
    Condition \eqref{equation 0 is regular} is not satisfied: $M_{\mathbf{Y}}((0,\varepsilon_2))=\infty$; furthermore we have that \eqref{equation x^2 bound} is satisfied as 
    \begin{equation*}
        \int_{v=0}^{\varepsilon_2} \int_{y=0}^v s_{\mathbf{Y}}(y) m_{\mathbf{Y}}(y) dy m_{\mathbf{Y}}(v)dv = \int_{v=0}^{\varepsilon_2} \int_{y=0}^v \frac{1}{2} dy \frac{1}{2v} dv = \frac{\varepsilon_2}{4}.
    \end{equation*}
    Therefore the Pitman-Yor excursion measure satisfies \eqref{equation levy measure unbounded variation} by Theorem \ref{theorem l\'evy measure condition}. Furthermore, for $q<q'_0$
    \begin{equation*}
        \int_{v=0}^{\varepsilon_2} g(v)m_{\mathbf{Y}}(v)dv \leq \int_{v=0}^{\varepsilon_2} \frac{C'_q}{2}v^{q-1} dv = \frac{C'_q}{2q} < \infty
    \end{equation*}
    and so $M^0_{\mathbf{N}, \mathbf{X}}(T)$ is finite for a stopped scaffolding-and-spindles pair $(\mathbf{N},\mathbf{X})$ by Theorem \ref{theorem aggregate mass summability}.
\end{example}
Further examples that satisfy Assumption \ref{assumption a} but not Assumptions \ref{assumption b} or \ref{assumption c} that could in theory generate branching IP-evolutions include reflected Brownian motion with or without a constant drift on $E=[0,c]$ for some finite $c$; reflected Brownian motion on $E=[0,\infty)$ with a constant negative drift; reflected Ornstein-Uhlenbeck process on $E=[0,\infty)$. All these processes satisfy \eqref{equation levy measure bounded variation} (as $0$ is a regular boundary point for these processes) and so we have aggregate health summability by Theorem \ref{theorem aggregate mass summability}. 
\subsection{Related work}\label{subsection related work}
Interval partitions can be found in several different areas of mathematics. As noted in \cite[Introduction]{construction}, interval partitions are the ranges of subordinators \cite{pitman92subordinator}, the complements of the zero-set of a Brownian bridges \cite[Example 3]{gnedinpitman2005}, the products of stick-breaking schemes \cite[Example 2]{gnedinpitman2005}, as well as the limits of compositions of natural numbers \cite{gnedin1997}. Furthermore, in our setting the class of interval partitions, $\mathcal{I}_H$, along with the metric $d'_H$ that we define in Definition \ref{definition d'_H metric}, serves as a useful complete and separable metric space to explore branching processes with diffusively evolving characteristics and to produce continuous Markovian processes. That said, we feel that these diffusions on the space of interval partitions should be taken as objects of interest in their own right.
\par
Using only the lifetimes of the individuals, and not the health excursions, the branching processes we study give rise to splitting trees, which can be found in the literature in e.g. Geiger and Kersting \cite{Geiger1997DepthFirstSO}. A splitting tree is a rooted planar tree that represents the type of branching processes that we study. It is a rooted tree where branch points represent births, leaves represent deaths, and the edges of varying length represent the times between birth or death events. Starting with the ancestor at the root, one can trace the path along the edge until a branch or boundary point. In the case of a boundary point, the ancestor had no offspring and the tree is a single edge of length given by the ancestors lifetime. In the case of a branch point, the left-most edge denotes the ancestor and the other edge denotes its first offspring. Continuing to follow the left-most edges at each branch point gives the collection of edges that represent the ancestor, and the sum of lengths is the lifetime. All the other edges at the branch points denote the offspring of the ancestor. Repeating this process for all descendants would return all the information about the branching process. The collection of individuals alive in the branching process at time $t$ is therefore in bijection with the collection of points of distance $t$ from the root.
\par
The skewer process draws similarities to the Ray-Knight theorems, see Revuz and Yor \cite[Theorems XI.2.1,XI.2.3]{revuz2013continuous} for a reference, which show that the local time process, indexed by level, of certain stopped Brownian motions are squared Bessel processes. In our setting we have interval partitions instead of local times, although, we find in Lemma \ref{lemma diversity is local time} that the local time at any fixed level is a function of the interval partition at that level a.s. In Forman et al. \cite{forman2018uniform} they show a stronger result for the self-similar cases: this same function applied to the interval partitions at all levels simultaneously gives the local times at those levels a.s. We still believe that this holds in non-self-similar cases but our methods were not powerful enough to prove it. 
\par
A motivation for our primary source, \cite{construction}, and related papers by the same authors, \cite{forman2018interval, forman2018uniform, forman2021diffusions, forman2020diffusions}, is to generate interval-partition-valued diffusions to explore diffusively evolving masses of the sub-trees of a spine of a continuum random tree. This is part of their monograph, \cite{forman2023aldous}, to solve the Aldous conjecture, Aldous \cite{aldous1999wright}, showing the existence of a diffusion on the space of continuum trees with the Brownian Continuum Random tree as the stationary distribution.
\subsection{Structure of this paper}\label{subsection structure of the paper}
Section \ref{section preliminaries} explores the relevant literature we use to prove our results. In Section \ref{subsection ito excursion theory} we look at the general theory of SPLPs with reference Bertoin \cite{bertoin1996levy}, and state key results about excursions of SPLPs from Pardo et al. \cite{spectrally_negative}. In Section \ref{subsection Pitman-Yor theory} we give a summary of the Pitman-Yor excursion measure from \cite{pitman1982decomposition} as well as the Williams disintegration for the Pitman-Yor excursion measure in Theorem \ref{definition Williams decomposition}, and we give an overview of the theory of disintegration. The disintegration theory also allows us to consider conditional measures such as $\nu(\cdot | \zeta=t)$, for $t>0$, which denote the probability distribution of a health excursion for an individual with a lifetime equal to $t$. In Section \ref{subsection Green's functions} we explore Green's functions as these along with the Williams decomposition are useful for calculations for all our main results. 
\par
Then in Section \ref{section initial results} we prove results for diffusions under Assumption \ref{assumption a}. In Section \ref{section l\'evy measure and other conditions} we prove Theorem \ref{theorem l\'evy measure condition} which is necessary to have the SPLP \say{scaffolding}. Then, in Section \ref{subsection general rsplps} we prove Theorem \ref{theorem aggregate mass summability} which determines that the levels of health are summable at a fixed point in time, and this is necessary to produce interval partitions. Section \ref{subsection disintegration proof} proves necessary disintegration results for the excursion measure of the SPLP. In Section \ref{section space of interval partitions} we restrict cases in Assumption \ref{assumption a} to find the possible initial states for an IP-evolution.
\par
In Section \ref{section model assumptions} we consider the parameters for the untransformed diffusion under Assumption \ref{assumption b}, which we use in Section \ref{subsection continuity local time} where we show conditions for the recurrent SPLP of unbounded variation to have local time that is jointly continuous in time and space. This is required in our Kolmogorov-\v{C}entsov continuity criterion argument in Section \ref{subsection continuity in aggregate mass} along with some further calculations in Section \ref{subsection diffusion deviations}. We do not use Kolmogorov-\v{C}entsov argument for diffusions under Assumption \ref{assumption c} as we do not have an equivalent to the continuity of local time result that was given in Section \ref{subsection continuity local time} that only holds for to SPLPs of unbounded variation. We and Theorems \ref{theorem diffusion} (without the strong Markov property) and \ref{theorem interval partitions regular} in Section \ref{section simple markov}. We then show the strong Markov property for Theorem \ref{theorem diffusion} under Strong Assumption \ref{assumption b} in Section \ref{subsection strong markov} by showing continuity in the initial state. 

%% file: preliminaries.tex
\section{Preliminaries}\label{section preliminaries}
In this section we give an overview of the results we use from the existing literature. We explore the excursions of a spectrally one-sided L\'evy process away from level $0$ in Bertoin \cite{bertoin1996levy} and Pardo, P\'erez and Rivero \cite{spectrally_negative}. We also use the excursion measure of a one-dimensional diffusion from Pitman and Yor \cite{pitman1982decomposition}; in particular we use the Williams decomposition \cite[Description (3.3)]{pitman1982decomposition} which describes the Pitman-Yor excursion measure. We require a notion of a conditional expectation for $\sigma$-finite measures, and for this we use the concept of a disintegration from Chang and Pollard \cite{chang1997conditioning}. Finally, we look at Revuz and Yor \cite{revuz2013continuous} for the Green's function for real-valued diffusions and use a result from Karlin and Taylor \cite[Section 15]{karlin1981second} which applies to such diffusions under Assumption $A3$.
\par
We define a complete and separable metric space of interval partitions $(\mathcal{I}_H,d'_H)$ from Forman et al. \cite{metric2020diversity}.
\begin{definition}[Interval partition metric space $(\mathcal{I}_H, d'_H)$]\label{definition d'_H metric}
Recall $\mathcal{I}_H$, the set of interval partitions, from Definition \ref{definition interval partition}. To define $d'_H$, we define correspondences between two interval partitions and distortions of the correspondences. We use the standard notation of $[n]=\{1, \dots, n\}$ for $n \in \mathbb{N}$.
\begin{enumerate}
    \item A \emph{correspondence} between two interval partitions $\beta$ and $\gamma$ is two collections of intervals $U_1, \dots, U_N$ from $\beta$ and $V_1,\dots, V_N$ from $\gamma$, for some $N \in \mathbb{N}$, such that each are ordered with respect to the left-to-right ordering of their interval partitions.
    \item The \emph{distortion} of a correspondence $(\beta, \gamma, (U_i,V_i)_{i \in [N]})$ is defined to be
    \begin{equation*}
        \text{dis}(\beta,\gamma,(U_i,V_i)_{i \in [N]}) \coloneqq  \min \left\{
            \begin{matrix}
                \sum_{i=1}^N | \text{Leb}(U_i) - \text{Leb}(V_i)| + \norm{\beta} - \sum_{i=1}^N \text{Leb}(U_i), \\
                \sum_{i=1}^N | \text{Leb}(U_i) - \text{Leb}(V_i)| + \norm{\gamma} - \sum_{i=1}^N \text{Leb}(V_j)  \\
            \end{matrix}
        \right\}  
    \end{equation*}
    \item The distance \emph{$d'_H(\beta,\gamma)$} between interval partitions $\beta$ and $\gamma$ is defined to be the infimum of distortions over the set of all possible correspondences:
    \begin{equation*}
        d'_H(\beta,\gamma)\coloneqq  \inf\limits_{N \geq 0, (U_i,V_i)_{i \in [N]}} \text{dis}(\beta,\gamma,(U_i,V_i)_{i \in [N]}). 
    \end{equation*}
\end{enumerate}
\end{definition}
\par
We define path spaces and various functions on these spaces.
\begin{definition}\label{definition excursion spaces}
We define the following metric spaces of functions:
\begin{enumerate}
    \item Let $\mathcal{D}$ be the Skorokhod space of one-dimensional c\`adl\`ag functions $g:\mathbb{R} \rightarrow \mathbb{R}$. Let $\mathcal{C}$ be the metric space of continuous functions $f:\mathbb{R} \rightarrow \mathbb{R}$ under the uniform metric. 
    \item Let $X: \mathcal{D} \rightarrow \mathcal{D}$ be the canonical process on $\mathcal{D}$, and let $Y$ be its restriction to $\mathcal{C}$, given by
    \begin{equation*}
        X_t(\omega) = \omega(t),\hspace{0.25cm} \omega \in \mathcal{D}; \hspace{0.5cm} Y_t(\omega) = \omega(t),\hspace{0.25cm} \omega \in \mathcal{C}.
    \end{equation*}
    We have defined two processes here, $X$ and $Y$ for clarity later on.
    \item Let $\mathcal{D}_{\text{exc}}\subset \mathcal{D}$ be the subspace of $\mathcal{D}$ consisting of c\`adl\`ag excursions away from $0$, i.e.
    \begin{equation*}
        \mathcal{D}_{\text{exc}}\coloneqq  \Big\{g \in \mathcal{D}\Big|\exists z \in (0,\infty) : g|_{(-\infty,0]\cup[z,\infty)}=0, g(x)>0 \text{   for all   } x \in (0,z)\Big\}.
    \end{equation*}
    Let $\mathcal{E}\subset \mathcal{C}$ of be the subspace of $\mathcal{D}_{exc}$ of continuous positive excursions from $0$, given by
    \begin{equation*}
        \mathcal{E}=\Big\{f \in \mathcal{D}\Big|\exists z \in (0,\infty): f|_{(-\infty,0]\cup[z,\infty)}=0; f \text{   positive, continuous on   } [0,z]\Big\}.
    \end{equation*}
\end{enumerate}
\end{definition}
\begin{definition}\label{definition stopping times}
We define functions on the spaces given in Definition \ref{definition excursion spaces}:
\begin{enumerate}
    \item For $g \in \mathcal{D}_{\text{exc}}$, we define $\zeta(g)\coloneqq \sup \{s\geq0: g(s)>0\}$, and we call $\zeta(g)$ the lifetime of $g$.
    \item For $g \in \mathcal{D}_{\text{exc}}$, we define $A(g)\coloneqq \sup \{g(s), 0 \leq s \leq \zeta(g) \}$, and we call $A(g)$ the amplitude of $g$. 
    \item For $g \in \mathcal{D}_{\text{exc}}$, we define $R(g)\coloneqq  \inf \{ s>0 : g(s)=A(g) \text{   or   } g(s-)=A(g) \}$.
    \item For $h \in \mathcal{D}$, we define $T^+_0(h) = \inf \{ t>0: h(t) \geq 0 \}$. We call $T^+_0(h)$ the crossing time of $0$ for $h$, and we call $h(T^+_0(h)-)$ and $h(T^+_0(h))$ the undershoot of $h$ and the overshoot of $h$ respectively.
    \item For $h \in \mathcal{D}$, and $E \subset \mathbb{R}$ measurable, we define $T_E(h)\coloneqq  \inf \{t > 0, h(t) \in E\} \in [0,\infty]$. For $y \in \mathbb{R}$, we define $T_y(h)\coloneqq  T_{\{y\}}(h)$ and we call $T_y(h)$ the hitting time of $y$ for $h$.
\end{enumerate}
\end{definition}
\subsection{It\^o excursion theory for recurrent SPLPs of bounded and of unbounded variation}\label{subsection ito excursion theory}
As discussed in Section \ref{subsection construction}, the scaffolding $\mathbf{X}$ is a recurrent SPLP either of bounded variation or of unbounded variation but with no Gaussian component. We discuss the excursions of the SPLPs in both cases.    
\par
Firstly, let $\mathbf{X}$ be a recurrent SPLP of bounded variation. For $x \in \mathbb{R}$, let $\mathbb{P}_x$ be the distribution on $\mathcal{D}$ of $\mathbf{X}$ started at $\mathbf{X}_0=x$, and let $\mathbb{E}_x$ denote the expectation operator. For all $y \in \mathbb{R}$, the level $y$ is irregular ($\mathbb{P}_x(T_x(X)=0)=1$), Bertoin \cite[Corollary VII.5]{bertoin1996levy}. For any fixed level $y \in \mathbb{R}$ the number of returns to $y$ after some time $t>0$ is finite a.s., \cite[Lemma 4.2]{lambert2010contour}, but the set of levels for which the number of returns to a level is infinite is dense in the range of the SPLP attains a.s., \cite[Theorem 5.3]{lambert2010contour}.
\par
Now instead consider a recurrent SPLP $\mathbf{X}$ of unbounded variation and with no Gaussian component. For all $y \in \mathbb{R}$, the level $y$ is regular ($\mathbb{P}_y(T_y(X)=0)=1$) and instantaneous ($\mathbb{P}_y(T_{\mathbb{R}\backslash y}(X)=0)=1$) for $\mathbf{X}$, \cite[Chapter IV.1-4]{bertoin1996levy}. Therefore we can define the local time $\ell^y$ for $\mathbf{X}$ at level $y$. The local time at all levels $y \in \mathbb{R}$ is an occupation density for the L\'evy process such that for any bounded and measurable function $h$ and any $t>0$ we have
\begin{equation*}
    \int_{y=-\infty}^\infty h(y) \ell^y(t)dy = \int_{s=0}^t h(\mathbf{X}_s) ds \hspace{0.5cm} \text{a.s.}
\end{equation*}
The local time is a continuous increasing function $\ell^y:[0,\infty) \rightarrow [0,\infty)$, and so we can define a right inverse $\tau^y$. Define $\tau^y:[0,\infty) \rightarrow [0,\infty)$ as, for $s>0$, $\tau^y(s) = \inf \{ t>0 : \ell^y(t) > s \}$. Note that $\tau^0$ is a subordinator, \cite[Theorem IV.8]{bertoin1996levy}. Then, following \cite[Equation IV.6]{bertoin1996levy}, we can define the excursion process $e\coloneqq (e_s, s \geq 0 : \tau^0(s)-\tau^0(s-)>0)$ for $\mathbf{X}$,
\begin{equation}\label{equation excursion process}
    e_s = \left(\mathbf{X}_{\tau^0(s-)+r}, 0 \leq r < \tau^0(s)-\tau^0(s-)\right),
\end{equation}
where each $e_s \in \mathcal{D}_{\text{exc}}$. There are only countably many points where $\tau(s)-\tau(s-)>0$ a.s. From \cite[Theorem IV.10]{bertoin1996levy} we note that $e$ is a Poisson point process.
\par
Define $\mathbf{S}=(\mathbf{S}_t)_{t \in [0,\infty)}$, where, for $t>0$, $\mathbf{S}_t \coloneqq  \sup_{0\leq s\leq t} \mathbf{X}_s$, i.e. $\mathbf{S}$ is the running supremum of $\mathbf{X}$. The process $\mathbf{S}-\mathbf{X}$ is non-negative for all $t>0$ and only hits the level $0$ at the points where $\mathbf{X}$ attains new suprema. Note that the level $0$ is also regular and instantaneous for $\mathbf{S}-\mathbf{X}$ and we can define $\overline{e}$ to be the excursion process for $\mathbf{S}-\mathbf{X}$ in the same way that we defined $e$ for $\mathbf{X}$. Let $n$ be the characteristic measure for $e$; i.e. the excursion measure for $\mathbf{X}$, and let $\overline{n}$ the characteristic measure for $\overline{e}$; i.e. the excursion measure for $\mathbf{S}-\mathbf{X}$. The Laplace exponent, $\psi:[0,\infty) \rightarrow \mathbb{R}$, for $\mathbf{X}$ is formally defined to be the function that satisfies
\begin{equation*}
    \mathbb{E}_0(e^{-\lambda X_t}) = e^{t \psi(\lambda)},
\end{equation*}
and in our setting is equal to 
\begin{equation}\label{equation laplace exponent}
    \psi(\lambda) = \int_{x=0}^\infty \left(e^{-\lambda x} - 1 + \lambda x\right) \Lambda(dx),
\end{equation}
see \cite[Chapter VII.1]{bertoin1996levy}.
\par
As the recurrent SPLP $\mathbf{X}$ has no Gaussian component we have that, up to a null set, all excursions under $n$ have the following structure a.s.: initially the excursion takes values in $(-\infty,0)$, before entering $(0,\infty)$ by a jump across $0$, and then taking values in $(0,\infty)$ until hitting $0$, \cite[Theorem 3]{spectrally_negative}. 
\par
We define a relation between the excursion measure $n$ of $\mathbf{X}$ and the excursion measure $\overline{n}$ for $\mathbf{S}-\mathbf{X}$, with the scale function $\hat{W}$ for $\hat{\mathbf{X}}\coloneqq -\mathbf{X}$, defined in \cite[Section VII.2]{bertoin1996levy}. The scale function $\hat{W}:[0,\infty) \rightarrow [0,\infty)$ is a continuous increasing function which relates to the Laplace exponent for $\hat{\mathbf{X}}$ in the following way: for all $\lambda>0$,
\begin{equation*}
    \int_{x=0}^\infty e^{-\lambda x}\hat{W}(x) dx = \frac{1}{\psi(\lambda)}.
\end{equation*}
Let $\hat{\mathbb{P}}_x$ be the probability measure for the law of the dual process $\hat{\mathbf{X}}$ started at $\hat{\mathbf{X}}(0)=x$. The scale function $\hat{W}$ for $\hat{\mathbf{X}}$ solves the two-sided exit problem for $\hat{\mathbf{X}}$, \cite[Theorem VII.8]{bertoin1996levy}: for $x,y>0$, the probability that $\hat{\mathbf{X}}$ makes its first exit from $[-x,y]$ at $y$ is
\begin{equation*}
    \hat{\mathbb{P}}_0\Big( \inf_{0 \leq s \leq T_y}(X_s) \geq -x\Big) = \frac{\hat{W}(x)}{\hat{W}(x+y)}.
\end{equation*}
In the following, we use notation such as
\begin{equation*}
    \overline{n}(g(t) \in A, \zeta>t) \coloneqq  \overline{n}(\{ g \in \mathcal{D}_{\text{exc}}: g(t) \in A, \zeta(g)>t\}),
\end{equation*}
where $A$ is a measurable set of $\mathbb{R}$ in this setting. To relate $n$ to $\overline{n}$, we use an instance of a general result for recurrent SPLPs from \cite[Theorem 3 (iii)]{spectrally_negative}: there exists $k_1>0$ such that for all $y$, $z>0$ 
\begin{equation}\label{equation PPR}
    n\left(-g\left(T^+_0-\right) \in dy, g\left(T^+_0\right) \in dz, T^+_0>t\right) = k_1\int_{x=0}^\infty \Big(\hat{W}(x)-\hat{W}(x-y)\Big) \overline{n}(g(t) \in dx, \zeta>t)\Lambda(dz+y) dy,
\end{equation}
where $\hat{W}(t)\coloneqq 0$ for $t<0$.
\par
For $x>0$, when we consider events up of $\hat{\mathbf{X}}$ up to a finite fixed time $t$ (i.e. any event $A$ such that $A \in \mathcal{F}_t$, the natural filtration for the process $\mathbf{X}$), $\hat{\mathbb{P}}^\uparrow_x$ is the distribution on $\mathcal{D}$ of $\hat{\mathbf{X}}$ when started from $x$ and conditioned to be positive in $(0,t)$, which can be found in \cite[Equation VII.6]{bertoin1996levy}:
\begin{equation*}
    \hat{\mathbb{P}}^\uparrow_x (X_t \in dy) \coloneqq  \frac{\hat{W}(y)}{\hat{W}(x)}\hat{\mathbb{P}}_x(X_t \in dy, T_{(-\infty,0)}(X)>t).
\end{equation*}
Furthermore, \cite[Proposition VII.14]{bertoin1996levy} states that the probability measure for $\mathbf{X}$ conditioned to stay positive when started at level $0$, $\mathbb{P}^\uparrow_0$, can be defined by the limit $\lim_{x \downarrow 0}\mathbb{P}^\uparrow_x$. Then, from \cite[Proposition VII.15]{bertoin1996levy}, we have that for some $k_2>0$ and for any $x$, $t>0$,
\begin{equation}\label{equation overline n}
    \overline{n}(g(t) \in dx, \zeta(g)>t) = \frac{k_2}{\hat{W}(x)} \hat{\mathbb{P}}_0^\uparrow(X_t \in dx).
\end{equation}
\subsection{Pitman-Yor excursion theory for one-dimensional diffusions and disintegration theory}\label{subsection Pitman-Yor theory}
Consider a diffusion that satisfies Assumption $A1$. We construct an excursion measure following Pitman and Yor \cite[Section 3]{pitman1982decomposition}. This is equivalent to the It\^o excursion when $0$ is a regular boundary point, but is a more general excursion measure as it also applies in the case where $0$ is an exit boundary point.
\par
Before we define the excursion measure, we need to define two more diffusions -- the $0$-diffusion and the $\uparrow$-diffusion. For $x>0$, the $0$-diffusion started at $x$ is simply the original diffusion started at $x$ which we stop when it hits level $0$. For $x>0$, the $\uparrow$-diffusion started at the level $x$ has law given by the limit, as $y \uparrow \infty$, of the conditional probability measure for the original diffusion (again started at $x$) conditioned on $\{T_y<T_0\}$; details of this can be found in \cite[Section 3]{pitman1982decomposition}. We can define the $\uparrow$-diffusion starting at $0$ and in this case the process never returns to level $0$. The scale function $s^\uparrow$ and the speed measure $M^\uparrow$ of the $\uparrow$-diffusion can be determined in terms of $s$ and $M$: 
\begin{equation}\label{equation s',s,m uparrow}
    s^\uparrow (x) = \frac{-1}{s(x)}; \hspace{0.5cm} M^\uparrow (dx) = s(x)^2 M(dx).
\end{equation}
The scale function $s^\uparrow$ result follows from the defining characteristic of a scale function in \eqref{equation scale function property} and some simple calculations. The speed measure result is noted in the literature in e.g. Salminen et al. \cite{Salminen2007}. With the $\uparrow$-diffusion and the $0$-diffusion we can define the Pitman-Yor excursion measure $\nu$.
\begin{definition}[Pitman-Yor excursion measure $\nu$]\label{definition Pitman-Yor excursion measure}
We define the $\sigma$-finite excursion measure $\nu$ on the space of positive continuous excursions $\mathcal{E}$ by the following conditions:
\begin{enumerate}
    \item $\nu(f\equiv 0) = \nu(f(0) \neq 0) = 0$;
    \item for all $y \in E \backslash \{0\}$ we have $\nu\left(T_y(f)<\infty\right) = 1/s(y)$;
    \item for all $y \in E \backslash \{0\}$, the conditional distribution $\nu\left(\cdot |T_y(f) < \infty\right)$ is the law of the following $\mathbf{V}=\left(\mathbf{V}_t\right)_{t \in [0,\infty)}$. The process $\left(\mathbf{V}_t\right)_{t \in [0,T_y)}$ has the law of an $\uparrow$-diffusion started at $0$ and run until it hits $y$, and $\left(\mathbf{V}_{T_y +t}\right)_{t \in [0,\infty)}$ has the law of an independent $0$-diffusion started at $y$.
\end{enumerate}
\end{definition}
The interpretation for $\nu$ is that the excursions come into $E \backslash \{0\}$ from the level $0$ according to the $\uparrow$-diffusion and then evolve in $E \backslash \{0\}$ according to the $0$-diffusion.
\par
We give an overview of the results from Chang and Pollard \cite{chang1997conditioning} on the theory of disintegration, which will allow us to condition over infinite dimensional measures such as the Pitman-Yor excursion measure $\nu$.
\begin{figure}
    \centering
    \includegraphics[width=\textwidth]{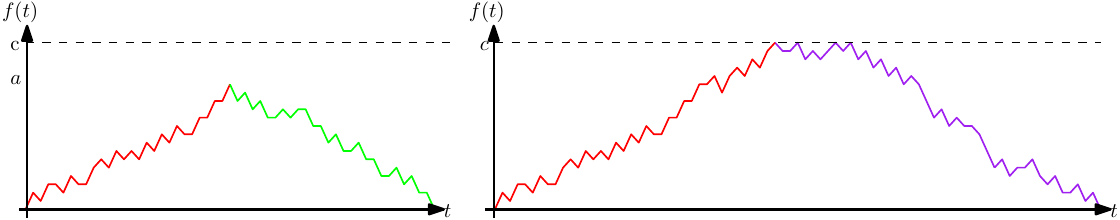}
    \caption{Left: an excursion with amplitude $a \in (0,c)$, constructed from an $\uparrow$-diffusion (red) and a time-reversed $\uparrow$-diffusion (green) both started at $0$ and stopped at $a$. Right: an excursion with amplitude $c$, constructed from an $\uparrow$-diffusion stopped at $c$ (red) and a $0$-diffusion started at $c$ and stopped at $0$ (purple).}
    \label{figure diffusion hit level c}
\end{figure}
\begin{definition}[Disintegration]\label{definition disintegration}
Let $\lambda$ be a $\sigma$-finite measure on a metric space $\mathcal{X}$ and let $S: \mathcal{X} \rightarrow \mathbb{R}$ be a measurable map. Then $\lambda$ has a disintegration $\{\lambda(\cdot| S=t)\}_{t \in [0,\infty)}$ with respect to some measurable $S: \mathcal{X} \rightarrow \mathbb{R}$, known as an $S$-disintegration, if
\begin{enumerate}
    \item $\lambda(\cdot | S=t)$ is a probability measure on $\mathcal{X}$ concentrated on $\{S=t\}$. In other words, $\lambda(\{S \neq t\}|S=t)=0$ for $\lambda(S \in \cdot)$-almost all $t$;
    \item for any non-negative measurable function $F$ on $\mathcal{X}$, $t \mapsto \int_{g \in \mathcal{X}} F(g) \lambda(dg|S=t)$ is measurable and
    \item $\int_{g \in \mathcal{X}} F(g) \lambda(dg) = \int_{t \in \mathbb{R}} \int_{g \in \mathcal{X}} F(g) \lambda(dg|S=t) \lambda(S \in dt)$.
\end{enumerate}
\end{definition}
We abuse notation and write $s(dw)$, for $w \in E \backslash \{0\}$, to denote the Lebesgue-Stieltjes measure on $E \backslash \{0\}$ as defined from the continuous increasing scale function $s$. Formally, $s(dw)$ is the Borel measure $\phi(dw)$ for which we have for any $0\leq x \leq y <c$ that $\phi((x,y))=s(y)-s(x)$. We remark that under Assumption $A3$ we have that $s(dw)=s'(w)dw$ but that $s$ is not in general differentiable. 
\par
The following description of the Williams decomposition of $\nu$ depends on whether $s(c)<\infty$. If $s(c)<\infty$ then the level $c$ is a regular boundary point for the diffusion, and therefore the diffusion can hit $c$ in finite time. In this case the process reflects off the point $c$ and the Williams decomposition of the excursion measure therefore includes a case where the amplitude $A$ is equal to $c$ which can occur with positive probability. The graph on the right in Figure \ref{figure diffusion hit level c} shows an excursion that hits level $c$ in this setting. 
\begin{theorem}[Williams decomposition of $\nu$]\label{definition Williams decomposition}
Consider a diffusion that satisfies $A1$. The its Pitman-Yor excursion measure $\nu$ is the unique measure that satisfies the following:
\begin{enumerate}
    \item $\nu\left(A \in dw \right)= (s(dw)/s(w)^2)$ for $w \in (0,c)$. If $s(c)$ is finite, we further have that $\nu(A=c)=1/s(c)$.
    \item For $w \in (0,c)$, a process $\mathbf{V}$ under $\nu(\cdot| A = w)$ has an a.s. unique time at which it attains its amplitude, given by $R$, so that $\mathbf{V}_R=w$. Furthermore, the processes
    \begin{equation*}
        \left(\mathbf{V}_t, 0 \leq t \leq R\right) \text{   and   } \left(\mathbf{V}_{\zeta - t}, 0 \leq t \leq \zeta - R\right)
    \end{equation*}
    are independent $\uparrow$-diffusions each started at the level $0$ and then stopped at time $T_w$. See the graph on the left in Figure \ref{figure diffusion hit level c}.
    \item A process $\mathbf{V}$ under $\nu(\cdot| A = c)$ evolves as the following. The process $\left(\mathbf{V}_t\right)_{t \in [0,T_c)}$ has the law of an $\uparrow$-diffusion started at $0$ run until it hits $c$, and $\left(\mathbf{V}_{T_c +t}\right)_{t \in [0,\infty)}$ has the law of an independent $0$-diffusion started at $c$. See the graph on the right in Figure \ref{figure diffusion hit level c}.
\end{enumerate}
\end{theorem}
The following theorem states sufficient conditions for a function $S:\mathcal{X} \rightarrow \mathbb{R}$ on a complete and separable space to disintegrate a measure $\lambda$.
\begin{theorem}[Existence Theorem for a disintegration for $\lambda$]\label{theorem disintegration}
Suppose $(\mathcal{X}, \Sigma(\mathcal{X}), \lambda)$ is a complete and separable space and that $S$ is a measurable map. Then $\lambda$ has an $S$-disintegration. Furthermore, the measures $\lambda(\cdot | S=t)$ are uniquely determined up to an almost sure equality: that is to say that if $\lambda_t$ is another $S$-disintegration then $\lambda(S^{-1}(t \in [0,\infty]: \lambda_t \neq \lambda(\cdot | S=t))) = 0$.
\end{theorem}
Definition \ref{definition disintegration} and Theorem \ref{theorem disintegration} are given in greater generality in \cite[Definition 1]{chang1997conditioning} and \cite[Theorem 1]{chang1997conditioning} respectively -- what we refer to as an $S$-disintegration they refer to as a $(S,S\lambda)$-disintegration. The Williams decomposition in Theorem \ref{definition Williams decomposition} satisfies the conditions of Definition \ref{definition disintegration} and so is an $A$-disintegration of $\nu$. Furthermore, by Theorem \ref{theorem disintegration} there exists an $\zeta$-disintegration of $\nu$.
\subsection{Green's functions of diffusions}\label{subsection Green's functions}
In this section we explore the properties of Green's functions of diffusions that we use in this paper. For a probability measure $\mathbb{Q}$ we let $\mathbb{E}_\mathbb{Q}$ denote the expectation operator under $\mathbb{Q}$. For a diffusion on a state space $E$ and for $x \in E$ let $\mathbb{Q}_x$ denote the probability measure of the diffusion started at $x$. This first result is from \cite[Corollary VII.3.8]{revuz2013continuous}.
\begin{theorem}\label{theorem Green's function general}
Consider a diffusion on a state space $E\subseteq [0,\infty)$ which is of the form $E=[0,c)$ or $E=[0,c]$ for some $c \leq \infty$. Let $a$, $w \in E$ with $a<w$, and let $F:\mathbb{R}\rightarrow \mathbb{R}$ be a locally bounded Borel function. We have for $x \in (a,w)$ that
\begin{equation*}
     \mathbb{E}_{\mathbb{Q}_x}\left( \int_{s=0}^{T_a\wedge T_w} F(Y_s) ds \right) = \int_{v=0}^w F(v) G_{(a,w)}(x,v) M(dv)
\end{equation*}
where $G_{(a,w)}$ is the Green's function for the diffusion on $(a,w)$ as in \eqref{equation greens function original}.
\end{theorem}
The next result is from \cite[Proposition VII.3.10]{revuz2013continuous} and includes some immediate consequences that apply to our setting.
\begin{lemma}\label{theorem diffusion started at c}
    Consider a diffusion which satisfies $A1$ and has speed measure $M$ such that $M((0,b))<\infty$ for some $b \in E$ and let $F:\mathbb{R} \rightarrow [0,\infty)$ be a locally bounded Borel function. Then, for $w \in E$, we have
    \begin{equation*}
        \mathbb{E}_{\mathbb{Q}_0}\left( \int_{s=0}^{T_w} F(Y_s) ds \right) = \int_{v=0}^w F(v) (s(w)-s(v)) M(dv).
    \end{equation*}
    In particular, if a diffusion $\mathbf{Y}$ satisfies $A1$, $A2$ and $s(c)$ is finite then $c$ is a regular boundary point for $\mathbf{Y}$ and we can define the diffusion $\hat{\mathbf{Y}} \coloneqq c-\mathbf{Y}$. Also, using \say{hats} for anything associated with $\hat{\mathbf{Y}}$, we have that $\hat{\mathbf{Y}}$ satisfies $A1$, $M_{\hat{\mathbf{Y}}}((0,b))<\infty$ for some $b \in E$ and
    \begin{equation*}
        \begin{split}
            \mathbb{E}_{\mathbb{Q}_c}\left( \int_{s=0}^{T_0} F(Y_s) ds \right) & = \mathbb{E}_{\hat{\mathbb{Q}}_0}\left( \int_{s=0}^{T_c} F(c-Y_s) ds \right) \\
            & =  \int_{v=0}^c F(c-v) (s_{\hat{\mathbf{Y}}}(c)-s_{\hat{\mathbf{Y}}}(v)) M_{\hat{\mathbf{Y}}}(dv) =  \int_{v=0}^c F(v) s_{\mathbf{Y}}(v) M_{\mathbf{Y}}(dv). \\
        \end{split}
    \end{equation*}
\end{lemma}
We now look at results about $\uparrow$-diffusions. For $x \in E$, we define $\mathbb{Q}^\uparrow_x$ to be the distribution of the $\uparrow$-diffusion started from $x$. This final lemma uses results from Karlin and Taylor \cite[Section 15.3]{karlin1981second}.
\begin{lemma}\label{lemma solution green differential}
Consider a diffusion that satisfies $A1$, $A3$ and has infinitesimal drift $\mu$ and infinitesimal variance $\sigma^2$ that are continuous on $E$. Then we can define the $\uparrow$-diffusion in this setting and we have that the infinitesimal drift $\mu^\uparrow$ and infinitesimal variance $(\sigma^2)^\uparrow$ for the $\uparrow$-diffusion are the continuous functions
\begin{equation}\label{equation mu and sigma uparrow}
\mu^\uparrow(x) = \mu(x) + \frac{s'(x)}{s(x)} \sigma^2(x); \hspace{0.5cm} (\sigma^2)^\uparrow (x) = \sigma^2(x).
\end{equation}
Furthermore, for $a>0$, we define functions $U_{(a,w)}$ and $V_{(a,w)}$ from $U$ and $V$ in \cite[Equation 15.3.37]{karlin1981second}: for $x \in (a,w)$ define
\begin{equation}\label{equation U differential solution}
    U_{(a,w)}(x) \coloneqq  \mathbb{E}_{\mathbb{Q}^\uparrow_x} \Big((T_a \wedge T_w)^2\Big), \hspace{0.5cm} V_{(a,w)}(x) \coloneqq  2 \mathbb{E}_{\mathbb{Q}^\uparrow_x} (T_a \wedge T_w).
\end{equation}
Then $U_{(a,w)}(\cdot)$ is the solution to the following differential equation in $x$:
\begin{equation}\label{equation U differential}
    \frac{1}{2} (\sigma^2)^\uparrow(x) U''_{(a,w)}(x) + \mu^\uparrow(x) U'_{(a,w)}(x) + V_{(a,w)}(x) = 0.
\end{equation}
Furthermore, there exists a unique solution $(U_{(a,w)}(x))_{x \in (a,w)}$ to \eqref{equation U differential} with boundary conditions $U_{(a,w)}(a)=0$ and $U_{(a,w)}(w)=0$, which is given by
\begin{equation*}
    U_{(a,w)}(x) = k_{(a,w)} (s^\uparrow (w) - s^\uparrow(x)) + \int_{y=x}^w \int_{z=a}^y (s^\uparrow)'(y) V_{(a,w)}(z)m^\uparrow(z)dydz \text{   for   } x \in (a,w),
\end{equation*}
where
\begin{equation*}
    k_{(a,w)} = \frac{-1}{s^\uparrow(w) - s^\uparrow(a)} \int_{y=a}^w \int_{z=a}^y (s^\uparrow)'(y) V_{(a,w)}(z)m^\uparrow(z)dzdy.
\end{equation*}
\end{lemma}
We prove the solution satisfies the differential equation in Lemma \ref{lemma green differential} in the Appendix. The formulae for $\mu^\uparrow$ and $(\sigma^2)^\uparrow$, stated in \eqref{equation mu and sigma uparrow}, can be found in the literature, e.g. Pinsky \cite[p. 371]{pinsky}. In \eqref{equation expected hitting time} we have $V_{(a,w)}(x)$ in terms of the speed measure and the Green's function for the diffusion on $(a,w)$ which is given in \eqref{equation greens function original}. We can then make substitutions for $s^\uparrow$ and $m^\uparrow$ using \eqref{equation s',s,m uparrow}. We take monotonic limits as $a$ decreases to $0$ and then as $x$ decreases to $0$ to obtain
\begin{equation}\label{equation V final}
    \mathbb{E}_{\mathbb{Q}^\uparrow_0}\big( T_w \big) = \int_{v=0}^w (s^\uparrow(w)-s^\uparrow(v)) M^\uparrow(dv) = \int_{v=0}^w \Big(\frac{1}{s(v)} - \frac{1}{s(w)}\Big) s(v)^2 M(dv) \text{   for   } w \in E,
\end{equation}
\begin{equation}\label{equation U final}
   \mathbb{E}_{\mathbb{Q}^\uparrow_0} \big(T_w^2\big) = 2 \int_{y=0}^w \int_{z=0}^y \int_{v=z}^w \left(\frac{1}{s(v)} - \frac{1}{s(w)}\right) s(v)^2 m(v) dv s(z)^2 m(z) dz \frac{s'(y)}{s(y)^2} dy  \text{   for   } w \in E.
\end{equation}
This ultimately follows as $0$ is not a regular or exit boundary point for the $\uparrow$-diffusion and so we have $\lim_{a \downarrow 0} \mathbb{E}_{\mathbb{Q}^\uparrow_x}(T_a \wedge T_w) = \mathbb{E}_{\mathbb{Q}^\uparrow_x}(T_w)$ and $\lim_{x \downarrow 0} \mathbb{E}_{\mathbb{Q}^\uparrow_x}(T_w) = \mathbb{E}_{\mathbb{Q}^\uparrow_0}(T_w)$ and similarly for $\mathbb{E}_{\mathbb{Q}^\uparrow_0}(T_w^2)$.

%% file: initialresults.tex
\section{Scaffolding and spindles -- Assumption \ref{assumption a}}\label{section initial results}
In this section we give results for the interval partition produced at a fixed level $y \in \mathbb{R}$ by our scaffolding-and-spindles construction.
\par
\subsection{Proof of Theorem \ref{theorem l\'evy measure condition}}\label{section l\'evy measure and other conditions}
We use Green's functions to determine the L\'evy measure condition in Theorem \ref{theorem l\'evy measure condition}. We prove that with the conditions of Theorem \ref{theorem l\'evy measure condition} the lifetime measure under the Pitman-Yor excursion measure is the jump measure of a recurrent SPLP, and the conclusion of the theorem follows. Our diffusion takes values state space $E\subseteq[0,\infty)$ of the form $E=[0,c)$ or $E=[0,c]$ where $c \leq \infty$. For $b \in (0,c)$ we can split the excursions $f \in \mathcal{E}$ of the diffusion based on whether or not they have amplitude $A(f)$ greater than $b$:
\begin{equation}\label{equation x^2 upper split}
    \int_{x=0}^\infty (x^2 \wedge x) \nu(\zeta \in dx) = \int_{f \in \mathcal{E}} \left(\zeta(f)^2 \wedge \zeta(f)\right) \nu(df) \leq \int_{f \in \mathcal{E}} \zeta(f)^2 \mathbbm{1}_{A(f) \leq b} \nu(df) + \int_{f \in \mathcal{E}} \zeta(f) \mathbbm{1}_{A(f) > b} \nu(df).
\end{equation}
\begin{lemma}\label{lemma recurrent upper}
Consider a diffusion that satisfies $A1$ and $A2$ with Pitman-Yor excursion measure $\nu$. Then for $b \in (0,c)$ we have that
\begin{equation*}
    \int_{f \in \mathcal{E}} \zeta(f) \mathbbm{1}_{A(f) > b} \nu(df) < \infty.
\end{equation*}
\begin{proof}
The disintegration described in the Williams decomposition in Theorem \ref{definition Williams decomposition} gives the excursion conditional on $A=w$ as the concatenation of two $\uparrow$-diffusions, the second one being time-reversed so that the concatenated pair form a non-negative excursion of amplitude $w$, see the left graph in Figure \ref{figure diffusion hit level c}. Therefore for $w \in (0,c)$ the expectation under for the lifetime under $\nu(\cdot|A=w)$ is given by
\begin{equation*}
    \int_{f \in \mathcal{E}} \zeta(f) \nu(df|A=w) = 2\int_{f \in \mathcal{E}} R(f) \nu(df|A=w) = 2 \mathbb{E}_{\mathbb{Q}^\uparrow_0} (T_w). 
\end{equation*}
We have calculated the expectation of $T_w$ under $\mathbb{Q}^\uparrow_0$ in \eqref{equation V final}, and so
\begin{equation*}
    \int_{f \in \mathcal{E}} \zeta(f) \nu(df|A=w) = 4 \int_{v=0}^w \left(\frac{1}{s(v)} - \frac{1}{s(w)}\right) s(v)^2 M(dv).
\end{equation*}
We integrate over the amplitude $A=w \in (b,c)$ and take the upper bound of $1/s(v)$ for $1/s(v)-1/s(w)$ to obtain
\begin{equation}\label{equation zeta(f) integral}
    \begin{split}
        \int_{f \in \mathcal{E}} \zeta(f) \mathbbm{1}_{A(f) \in (b,c)} \nu(df) & = 4\int_{w=b}^c \int_{v = 0}^w \frac{s(v)^2}{s(w)^2} \left(\frac{1}{s(v)} - \frac{1}{s(w)}\right) M(dv)s(dw) \\
        & \leq 4\int_{w=b}^c \int_{v = 0}^w s(v) M(dv) \frac{1}{s(w)^2} s(dw). \\
    \end{split}
\end{equation}
Finally, we split the integral for $v \in [0,b]$ and for $v \in (b,w)$ to show the claim. 
\begin{equation*}
    \int_{w=b}^c \int_{v = 0}^b s(v) M(dv) \frac{1}{s(w)^2} s(dw) \leq \frac{1}{s(b)}\int_{v = 0}^b  s(v) M(dv),
\end{equation*}
which is finite by $A1$ and 
\begin{equation*}
    \int_{w=b}^c \int_{v = b}^w s(v)  M(dv)\frac{1}{s(w)^2} s(dw) = \int_{v = b}^c \int_{w=v}^c \frac{1}{s(w)^2} s(dw) s(v)  M(dv) \leq M((b,c))
\end{equation*}
which is finite by $A2$.
\par
The final part we need to consider are excursions that hit level $c$, and this occurs when $s(c)$ is finite (equivalently when $c$ is a regular boundary point of the diffusion). In this case, as stated in the Williams decomposition in Theorem \ref{definition Williams decomposition}, we have that conditional on $\{A=c\}$ the excursion follows a probability distribution that consists of an $\uparrow$-diffusion started at $0$ and stopped at $T_c$, the hitting time of $c$, followed by an independent $0$-diffusion started at $c$, i.e. a diffusion started at $c$ and stopped at $T_0$.
\par
We calculate the contributions. The first term is given for $w=c$ in \eqref{equation V final} and the second contribution is given in Lemma \ref{theorem diffusion started at c}. The combination is
\begin{equation}\label{equation excursion hits regular boundary c}
    \int_{f \in \mathcal{E}} \zeta(f) \nu(df,A=c) = \nu(A=c)\Big(\mathbb{E}_{\mathbb{Q}^\uparrow_0}(T_c)+\mathbb{E}_c(T_0)\Big) = \frac{1}{s(c)} \int_{y \in (0,c]} \Big(2 - \frac{s(y)}{s(c)}\Big)s(y) M(dy).
\end{equation}
This is less than $(2/s(c))\int_{y=0}^c s(y)M(dy)$. The result follows as $\int_{y=0}^b s(y)M(dy)$ is finite by $A1$ and we have that $\int_{y=b}^c s(y)M(dy)$ is finite because $s(c)$ is finite and as $M((b,c))$ is finite by $A2$.
\end{proof}
\end{lemma}
\begin{lemma}\label{lemma levy upper}
Consider a diffusion that satisfies Assumption \ref{assumption a} and let $\nu$ be the Pitman-Yor measure. Let $b \in (0,c)$ be as in $A4$. Then we have
\begin{equation}\label{equation lemma levy upper}
    \int_{f \in \mathcal{E}} \zeta(f)^2 \mathbbm{1}_{A(f) \leq b} \nu (df) \leq 12 \int_{v=0}^b \int_{z=0}^v  s(z) m(z) dz m(v) dv.
\end{equation}
\begin{proof}
As we assume $A3$ in this lemma we can write $M(dv)=m(v)dv$ and $s(dw)=s'(w)dw$, i.e. $M$ has a density $m$ with respect to the Lebesgue measure and $s$ is differentiable. As in the previous lemma, we use from Theorem \ref{definition Williams decomposition} that, for $w \in (0,b)$, $R$ and $\zeta-R$ are identically distributed under $\nu(\cdot|A=w)$ with the distribution of $T_w$ under $\mathbb{Q}^\uparrow_0$. By expanding $\zeta^2 = (R + (\zeta-R))^2$ we obtain the upper bound
\begin{equation*}
    \int_{f \in \mathcal{E}} \big(R(f)^2 + 2 R(f)\big(\zeta(f) - R(f)\big) + (\zeta(f) - R(f))^2\big) \nu (df|A=w) \leq 6 \int_{f \in \mathcal{E}} R(f)^2 \nu (df|A=w).
\end{equation*}
To find an upper bound for the right-hand-side we consider a different diffusion $\Tilde{\mathbf{Z}}$ with infinitesimal drift $\Tilde{\mu}$ and infinitesimal drift $\Tilde{\sigma}^2$ given by
\begin{equation*}
    \Tilde{\mu}(w) \coloneqq \mu(w), \hspace{0.25cm} \Tilde{\sigma}^2(w) \coloneqq \sigma^2(w) \text {   for   } w \in (0,b); \hspace{0.5cm} \Tilde{\mu}(w) \coloneqq \mu(b), \hspace{0.25cm} \Tilde{\sigma}^2(w) \coloneqq \sigma^2(b) \text{   for   } w \geq b.
\end{equation*}
Let $\Tilde{s}$, $\Tilde{m}$ and $\Tilde{\nu}$ be the scale function, speed measure density and Pitman-Yor excursion measure for $\Tilde{\mathbf{Z}}$. Then $\Tilde{\mathbf{Z}}$ is a diffusion that satisfies $A1$ and $A3$ and $\Tilde{\mu}$ and $\Tilde{\sigma}^2$ are continuous on $E$ and so we can apply Lemma \ref{lemma solution green differential} and, substituting the expectation for $T_w^2$, we obtain
\begin{equation*}
    \int_{f \in \mathcal{E}} R(f)^2 \Tilde{\nu} (df|A=w) = 2 \int_{y=0}^w \int_{z=0}^y \int_{v=z}^w \left(\frac{1}{\Tilde{s}(v)} - \frac{1}{\Tilde{s}(w)}\right) \Tilde{s}(v)^2 \Tilde{m}(v) dv \Tilde{s}(z)^2 \Tilde{m}(z) dz \frac{\Tilde{s}'(y)}{\Tilde{s}(y)^2} dy \text{   for   } w \in (0,b).
\end{equation*}
However as, for $w \in (0,b)$, we have $\Tilde{s}(w)=s(w)$, $\Tilde{m}(w)=m(w)$ and also $\nu(\cdot|A=w)=\Tilde{\nu}(\cdot|A=w)$ by the Williams decomposition in Theorem \ref{definition Williams decomposition}, we can substitute the terms in the equation above. We then integrate over $A=w \in (0,b)$ and take the same upper bound of $1/s(v)$ for $1/s(v) - 1/s(w)$ as in the previous lemma to get
\begin{equation*}
    \int_{f \in \mathcal{E}} \zeta(f)^2 \nu (df|A=w) \leq 12 \int_{w=0}^b \int_{y=0}^w \int_{z=0}^y \int_{v=z}^w s(v) m(v)dv  s(z)^2 m(z)dz \frac{s'(y)}{s(y)^2}dy \frac{s'(w)}{s(w)^2} dw \text{   for   } w \in (0,b).
\end{equation*}
We swap the order of integration and then extend the range of the now inner $y$ and $w$ integrals to upper bound the integral above by 
\begin{equation*}
    \int_{z=0}^b \int_{v=z}^b \int_{w=v}^c \int_{y=z}^c  \frac{s'(y)}{s(y)^2}dy \frac{s'(w)}{s(w)^2} dw s(v) m(v)dv s(z) m(z) dz \leq \int_{z=0}^b \int_{v=z}^b  m(v) dv s(z)m(z)dz,
\end{equation*}
where we again use the same upper bound of $1/s(v)$ for $1/s(v) - 1/s(w)$ and also in $z$. Finally we swap the order of integration of $z$ and $v$ to give the upper bound in the claim.
\end{proof}
\end{lemma}
\begin{proof}[Proof of Theorem \ref{theorem l\'evy measure condition}]
Consider a diffusion under $A1$ and $A2$ and for which \eqref{equation 0 is regular} holds (i.e. $0$ is a regular point for the diffusion). Then as in \eqref{equation zeta(f) integral}, but now without the amplitude constraint of $A(f)>b$, we note that
\begin{equation*}
    \int_{f \in \mathcal{E}} \zeta(f) \nu(df) = 4\int_{w=0}^c \int_{v = 0}^w \frac{s(v)^2}{s(w)^2} \left(\frac{1}{s(v)} - \frac{1}{s(w)}\right) M(dv) s(dw)
\end{equation*}
which we can rearrange to give
\begin{equation*}
    4 \int_{v=0}^c \int_{w = v}^c \Big(\frac{s(v)}{s(w)^2} - \frac{s(v)^2}{s(w)^3}\Big) s(dw)M(dv).
\end{equation*}
We split into two cases. When $s(c)$ is infinite the above term is equal to $2M((0,c))$ and therefore we have the Pitman-Yor excursion measure satisfies \eqref{equation levy measure bounded variation} if and only if \eqref{equation 0 is regular} holds. When $s(c)$ is finite we have the above integral is equal to
\begin{equation*}
    2 M((0,c)) - \frac{4}{s(c)}\int_{v=0}^c s(v)M(dv) + \frac{2}{s(c)^2}\int_{v=0}^cs(v)^2M(dv).
\end{equation*}
To complete the proof in this case it suffices to show that $\int_{v=0}^c s(v)M(dv)$ is finite. We have $\int_{v=0}^b s(v)M(dv)$ is finite by $A1$ and furthermore we have that $\int_{v=b}^c s(v)M(dv)$ is finite as $M((b,c))$ and $s(c)$ are both finite.
\par
For the second setting described in Theorem \ref{theorem l\'evy measure condition} (where the diffusion satisfies all of Assumption \ref{assumption a}) the result follows immediately from Lemmas \ref{lemma recurrent upper} and \ref{lemma levy upper}. 
\end{proof}
\subsection{Proof of Theorem \ref{theorem aggregate mass summability}}\label{subsection general rsplps}
Throughout this section we consider a diffusion that satisfies $A1$ and $A2$ and has a Pitman-Yor excursion measure $\nu$ that satisfies \eqref{equation levy measure unbounded variation}, which we recall means that for a scaffolding-and-spindles pair in this setting, the scaffolding $\mathbf{X}$ is an SPLP of unbounded variation with no Gaussian component. With this in mind, we begin with a lemma that applies to such L\'evy processes.
\begin{lemma}[Excursion undershoot-overshoot measure]\label{lemma excursion overshoot-undershoot measure}
Let $\mathbf{X}$ be a recurrent SPLP of unbounded variation with L\'evy measure $\Lambda$ and with no Gaussian component. Let $n$ be the excursion measure for the law of $\mathbf{X}$. For some constant $k_0>0$, the undershoot and overshoot joint measure under the excursion measure $n$ is given by the following measure on $(0,\infty)^2$
\begin{equation}\label{equation levy excursion limit}
    n(-g(T^+_0-) \in dy, g(T^+_0) \in dz) = k_0\Lambda(dz + y)dy.
\end{equation}
\begin{proof}
We substitute the expression for the excursion measure of $\mathbf{S}-\mathbf{X}$ from \eqref{equation overline n} into \eqref{equation PPR} to give
\begin{equation*}
    n\left(-g(T^+_0-) \in dy, g(T^+_0) \in dz, T^+_0>t\right) = k_1 \left(\int_{x=0}^\infty \left(\hat{W}(x)-\hat{W}(x+y)\right) \frac{k_2}{\hat{W}(x)} \hat{\mathbb{P}}^\uparrow(X_t \in dx)\right)\Lambda(dz+y) dy. 
\end{equation*}
For $x>0$ and $y>0$, $\left(\hat{W}(x)-\hat{W}(x-y)\right)/\hat{W}(x)\leq 1$. Also, $\lim_{x \downarrow 0}\left(\hat{W}(x)-\hat{W}(x-y)\right)/\hat{W}(x) = 1$ by monotonicity. Finally, for $X$ with law $\mathbb{P}^\uparrow$, we have that $X_t \rightarrow 0$ as $t \downarrow 0$ a.s. Therefore, as $n(T^+_0=0)=0$,
\begin{equation*}
    \begin{split}
        n\left(-g(T^+_0-) \in dy, g(T^+_0) \in dz\right) & = \lim_{t \downarrow 0} n\left(-g(T^+_0-) \in dy, g(T^+_0) \in dz, T^+_0>t\right) \\
        & = k_1 k_2 \left(\lim_{t \downarrow 0}\int_{x=0}^\infty \frac{\hat{W}(x)-\hat{W}(x-y)}{\hat{W}(x)} \hat{\mathbb{P}}^\uparrow(X_t \in dx)\right)\Lambda(dz+y) dy \\
        & = k_0 \Lambda(dz+y) dy \text{   where   } k_0\coloneqq k_1k_2,
    \end{split}
\end{equation*}
and where the limit of measures is a vague limit.
\end{proof}
\end{lemma}
\begin{definition}\label{definition spaces of N}
    For a complete and separable metric space $S$, let $\mathcal{N}(S)$ be the set of counting measures on $S$ that are boundedly finite. We equip $\mathcal{N}(S)$ with the $\sigma$-algebra $\Sigma(\mathcal{N}(S))$ generated by the evaluation maps $N \mapsto N(B)$ for Borel subsets $B \subseteq S$. Let $N \in \mathcal{N}([0,\infty) \times \mathcal{E})$ and define
    \begin{enumerate}
        \item $\text{len}(N) \coloneqq  \inf \{t>0:N([t,\infty) \times \mathcal{E})=0 \} \in [0,\infty]$;
        \item $N|_{[a,b]}$ is $N$ restricted to the interval $[a,b]$: $N|_{[a,b]} \coloneqq  \sum\limits_{t \in [a,b], (t,f_t) \text{   atom of   } N} \delta_{(t,f_t)}$.
        \item $N|^\leftarrow_{[a,b]}$ is $N$ with the time index shifted back by $a$: $N|^\leftarrow_{[a,b]} \coloneqq  \sum\limits_{t \in [a,b], (t,f_t) \text{   atom of   } N} \delta_{(t-a,f_t)}$;
        \item for an atom $(t,f_t)$ of $N$ we define, for $y \in (\xi_N(t-),\xi_N(t))$, the decomposition of the spindle $f_t$ between the component of $f_t$ that falls below level $y$, denoted as $\check{f}_t \in \mathcal{E}$, and the component of $f_t$ that is above level $y$, denoted as $\hat{f}_t \in \mathcal{E}$:
        \begin{equation*}
            \Check{f}^y_t(z)\coloneqq  f_t(z) \mathbbm{1}_{z \in [0,y-\xi_N(t-))}, \hspace{0.5cm} \Hat{f}^y_t(z)\coloneqq  f_t(y - \xi_N(t-) + z) \mathbbm{1}_{z \in [0,\xi_N(t)-y]},
        \end{equation*}
        and we write $\Check{f}_t\coloneqq \Check{f}^0_t$ and $\Hat{f}_t\coloneqq \Hat{f}^0_t$;
        \item we define a subspace of $\mathcal{N}([0,\infty) \times \mathcal{E})$. Let $\mathcal{N}^{\text{sp}}_{\text{fin}}$ be the set of counting measures with the additional properties:
            \begin{enumerate}
                \item $N(\{t\} \times \mathcal{E}) \leq 1$ for every $t \in [0,\infty)$,
                \item $N([0,t] \times \{f \in \mathcal{E}: \zeta(f)>z \})<\infty$ for every $t$, $z>0$, 
                \item $\text{len}(N)<\infty$ and limit $\xi_N(t)$ exists for each $t \in [0,\text{len}(N)]$;
            \end{enumerate}
        \item we define the subspace $\mathcal{N}^{\text{sp}} \subset \mathcal{N}$ by saying that $N \in \mathcal{N}^{\text{sp}}$ if and only if $N|_{[0,t]} \in \mathcal{N}^{\text{sp}}_{\text{fin}}$ for each $t \geq 0$.
    \end{enumerate} 
\end{definition}
With the above definitions, we can define the following sets that contain measures associated with c\`adl\`ag excursions away from $0$. We define
\begin{equation*}
    \mathcal{N}^{\text{sp}}_{\pm \text{cld}} \coloneqq \{N \in \mathcal{N}^{\text{sp}}_{\text{fin}}: \xi(N) \in \mathcal{D}_{\text{exc}}\},
\end{equation*}
\begin{equation*}
     \mathcal{N}^{\text{sp}}_{+ \text{cld}} \coloneqq \{N \in \mathcal{N}^{\text{sp}}_{\pm \text{cld}}: \inf_t \xi_N(t)=0\}, \hspace{0.5cm} \mathcal{N}^{\text{sp}}_{- \text{cld}} \coloneqq \{N \in \mathcal{N}^{\text{sp}}_{\pm \text{cld}}: \sup_t \xi_N(t)=0\},
\end{equation*}
to classify the possible excursions of $\mathbf{X}$ marked with the points of $\mathbf{N}$. We split $N \in \mathcal{N}^{\text{sp}}_{\pm \text{cld}} \backslash (\mathcal{N}^{\text{sp}}_{- \text{cld}} \cup \mathcal{N}^{\text{sp}}_{+ \text{cld}})$ into what we call its clade component $N^+ \in \mathcal{N}^{\text{sp}}_{+ \text{cld}}$ and its anti-clade component $N^- \in \mathcal{N}^{\text{sp}}_{- \text{cld}}$ by
\begin{equation}\label{equation N+ definition}
    N^-\coloneqq  N|_{[0,T^+_0)} + \delta(T^+_0, \Check{f}^0_{T^+_0}), \hspace{0.5cm} N^+ \coloneqq  \delta(0, \Hat{f}^0_{T^+_0}) + N|^\leftarrow_{(T^+_0,\infty)},
\end{equation}
(and we define $(N^-,N^+)$ to be $(N,0)$ for $N \in \mathcal{N}^{\text{sp}}_{- \text{cld}}$ and similarly for $N \in \mathcal{N}^{\text{sp}}_{+ \text{cld}}$). The terminology is inherited from our main reference \cite{construction}. The component $N^+$ represents all the descendents of an individual, and this is therefore a clade of the total population. The term anti-clade was then coined due to it being the inverse to the clade. Finally, the term bi-clade is the combination of the clade and the anti-clade. We also refer to elements of $\mathcal{N}^{\text{sp}}_{\pm \text{cld}}$, $\mathcal{N}^{\text{sp}}_{+ \text{cld}}$ and $\mathcal{N}^{\text{sp}}_{- \text{cld}}$ as bi-clades, clades and anti-clades respectively.
\par
We now define $V^y$ to be the set of all complete excursions and define $V^y_0$ to be $V^y$ along with the possibly incomplete first and last excursions in the following way from \cite[Definition A.1]{construction}.
\begin{definition}\label{definition set of excursion intervals}
    For a c\`adl\`ag function $g \in \mathcal{D}$ and $y \in \mathbb{R}$ we define
    \begin{equation*}
        V^y(g) \coloneqq \left\{ [a,b] \subset [0,\text{len}(g)] \left|
            \begin{matrix}
                a<b<\infty; g(t-) \neq y \neq g(t) \text{   for   } t \in (a,b); \\
                g(a-)=y \text{   or   } g(a)=y; g(b-)=y \text{   or   } g(b)=y  \\
            \end{matrix}
        \right.\right\}.  
    \end{equation*}
 The set $V^y_0(g) \supset V^y(g)$ then includes any incomplete first or last excursions. To define this rigourously, let 
 \begin{equation*}
    \begin{split}
        & T^y(g) \coloneqq \inf ( \{ t \in [0,\text{len}(g)]: g(t)=y \text{   or   } g(t-)=y \} \cup \{\text{len}(g)\}), \\
        & T^y_*(g) \coloneqq  \sup ( \{ t \in [0,\text{len}(g)]: g(t) = y \text{   or   } g(t-)=y\} \cup \{0\}). \\
    \end{split}
 \end{equation*}
 Then, for $y \neq 0$, we include $[0,T^y(g)] \cap [0,\infty)$ in $V^y_0(g)$, even if $T^y(g)$=0. Also, if $T^y_*(g) < \text{len}(g)$ or $g(\text{len}(g)) \neq y$, then we include $[T^y_*(g), \text{len}(g)] \cap [0,\infty)$ in $V^y_0(g)$. For $N \in \mathcal{N}^{\text{sp}}$ we define $V^y(N) \coloneqq V^y(\xi(N))$ and similarly for $V^y_0$.  
 \par
 Finally, we need to take care with possible non-degenerate excursions and so we carefully define intervals $I^y_N(a,b)$ to be one of $[a,b]$, $(a,b]$, $[a,b)$, $(a,b)$ as follows: we exclude the point $a$ iff $a<b$ and $g(a-)<y=g(a)$; we exclude $b$ iff both $a<b$ and $g(b-)=y<g(b)$. 
 \end{definition}
With this definition we can consider the measure on excursions of the process $\mathbf{X}\coloneqq \xi(\mathbf{N})$ where the jumps are marked by the spindles that generate them. For $N \in \mathcal{N}^{\text{sp}}_{\text{fin}}$, let $V^y(N)$ for $y \in \mathbb{R}$, $N \in \mathcal{N}^{\text{sp}}$ be the set of intervals of complete excursions of $\mathbf{X}$, so that $[a,b] \in V^y(N)$ means that the path $(\xi_N(a+t), 0\leq t \leq b-a)$ is an excursion about level $y$ -- we define this rigourously in Proposition \ref{proposition properties of levy process of unbounded variation} in the Appendix. As $\mathbf{X}$ is of unbounded variation we have that, for each level $y \in \mathbb{R}$, there are no degenerate excursions of $\mathbf{X}$ a.s. and $\mathbf{X}$ has local times. Furthermore, the local times are injective a.s. for each $y$ in the sense that for $[a,b]$, $[c,d] \in V^y$, $\ell^y(a) \neq \ell^y(c)$ unless $[a,b]=[c,d]$; this is shown in Millar \cite{millar1973exit} which summarises whether L\'evy processes either first jump across or first hit fixed levels. While it holds at any fixed level a.s., across all levels simulataneously this would be false. We define $F^0$ for $N \in \mathcal{N}^{\text{sp}}$ and $\mathbf{F}^0$ for $\mathbf{N}$ by
\begin{equation}\label{equation F^0 definition}
    F^0(N) \coloneqq  \sum\limits_{[a,b] \in V^0(N)} \delta_{\ell^0(a), N|^\leftarrow_{[a,b]}}, \hspace{0.5cm} \mathbf{F}^0\coloneqq F^0(\mathbf{N}),
\end{equation}
and we note that $\mathbf{F}^0$ is a Poisson random measure. We define the measure $\nu_{\text{cld}}$ on $\mathcal{N}^{\text{sp}}_{\pm \text{cld}}$ by 
\begin{equation*}
    \nu_{\text{cld}}(A)\coloneqq  \mathbf{E}(\mathbf{F}^0([0,1] \times A)) \text{   for a measurable set   } A \subset \mathcal{N}^{\text{sp}}_{\pm \text{cld}}.
\end{equation*}
The bi-clades here will have a unique jump that crosses zero, that is to say that the set of all bi-clades with no unique crossing jump has measure $0$ under $\nu_{\text{cld}}$, and the unique jump will have associated with it an excursion of our diffusion. Therefore we can discuss the crossing width $m^0$ of the unique crossing jump:
\begin{equation*}
    m^0(N) \coloneqq  \sum_{(r,f_r) \text{   atom of   } N} \max\{f_r((-\xi_N(r-))-), f_r(-\xi_N(r-))\}; \hspace{0.25cm} N \in \mathcal{N}_{\text{cld}}^{\text{sp}}.
\end{equation*}
The reason that we define $m^0$ in this way is that it makes sense for $N \in \mathcal{N}^{\text{sp}}_{- \text{cld}}\cup \mathcal{N}^{\text{sp}}_{+ \text{cld}}$. Indeed, we have
\begin{equation*}
    m^0(N) = f_{T^+_0(\xi_N)}(-\xi_N(T^+_0(\xi_N)-)) \text{   for   } N \in \mathcal{N}^{\text{sp}}_{\pm \text{cld}} \backslash (\mathcal{N}^{\text{sp}}_{- \text{cld}} \cup \mathcal{N}^{\text{sp}}_{+ \text{cld}}).
\end{equation*}
This second formulation for the crossing width is more intuitive. We use the following notation, similar in style to that used in Section \ref{subsection ito excursion theory}: take as an example, for a measurable set $A$ of $[0,\infty) \times (-\infty, 0] \times [0,\infty)$,
\begin{equation*}
    \nu_{\text{cld}}((m^0, \xi_N(T^+_0-), \xi_N(T^+_0)) \in A) \coloneqq  \nu_{\text{cld}}(\{N \in \mathcal{N}^{\text{sp}}_{\pm \text{cld}}: (m^0(N), \xi_N(T^+_0(\xi_N)-), \xi_N(T^+_0(\xi_N))) \in A\}).
\end{equation*}
\begin{proposition}\label{proposition crossing width density}
Consider a diffusion that satisfies $A1$ and $A2$ with Pitman-Yor excursion measure $\nu$ that satisfies \eqref{equation levy measure unbounded variation}. Let $(\mathbf{N},\mathbf{X})$ be a scaffolding-and-spindles pair for this diffusion. Then, in the setting described above, $\nu_{\text{cld}}(m^0 \in dx)=k_0M(dx)$ where $k_0$ is the constant in Lemma \ref{lemma excursion overshoot-undershoot measure} for the SPLP $\mathbf{X}$. 
\begin{proof}
We rearrange the result in Lemma \ref{lemma excursion overshoot-undershoot measure} to obtain
\begin{equation*}
    \nu_{\text{cld}}(-\xi_N(T^+_0-) \in dv, \xi_N(T^+_0)-\xi_N(T^+_0-) \in du) = k_0 dv \nu(\zeta \in du) \text{   for   } k_0>0, \text{   for   } 0<v<u,
\end{equation*}
where $k$ is the constant stated in Lemma \ref{lemma excursion overshoot-undershoot measure}. Recall that $\Lambda$ is the jump measure of $\mathbf{X}$ and therefore $\nu(\zeta \in \cdot)=\Lambda(\cdot)$ in this context. Let $H$ be a positive measurable function $H:[0,\infty) \rightarrow [0,\infty)$. We claim that
\begin{equation*}
    \nu_{\text{cld}}\big(H(m^0)\big) = \int_{N \in \mathcal{N}^{\text{sp}}_{\pm \text{cld}}} H(m^0(N)) \nu_{\text{cld}}(dN) = k_0 \int_{y \in E} H(y)M(dy).
\end{equation*}
First note that
\begin{equation*}
    \int\limits_{N \in \mathcal{N}^{\text{sp}}_{\pm \text{cld}}} H(m^0(N)) \nu_{\text{cld}}(dN) = \int_{u=0}^\infty \int_{v=0}^u \int_{x \in E} H(x) \nu_{\text{cld}}(m^0 \in dx, -\xi_N(T^+_0-) \in dv, \xi_N(T^+_0)-\xi_N(T^+_0-) \in du)
\end{equation*}
and in the next step we disintegrate $\nu_{\text{cld}}(m^0 \in \cdot)$ with respect to $(-\xi_N(T^+_0-), \xi_N(T^+_0))$. This is possible as the bi-clade measure $\nu_{\text{cld}}$ can be obtained from marking the jumps of an excursion of a c\`adl\`ag process under the excursion measure of the L\'evy process $\mathbf{X}$, \cite[Proposition 4.9]{construction}. We obtain
\begin{equation*}
    \begin{split}
        k_0 \int_{u=0}^\infty \int_{v=0}^u \int_{x \in E} & H(x) \nu_{\text{cld}}(m^0 \in dx| -\xi_N(T^+_0-) = v, \xi_N(T^+_0)-\xi_N(T^+_0-) = u) dv \nu(\zeta \in du) \\
        & = k_0 \int_{u=0}^\infty \int_{v=0}^u \int_{f \in \mathcal{E}} H(f(v)) \nu(df|\zeta = u) dv \nu(\zeta \in du), \\
    \end{split}
\end{equation*}
which we can rewrite as
\begin{equation*}
    k_0 \int_{u=0}^\infty \int_{f \in \mathcal{E}} \int_{v=0}^{\zeta(f)} H(f(v)) dv \nu(df|\zeta = u) \nu(\zeta \in du),
\end{equation*}
where $f$ represents the excursion associated with the unique crossing jump and $\nu$ is the Pitman-Yor excursion measure of our diffusion. We can express this in terms of an integral against the Pitman-Yor excursion measure as
\begin{equation*}
    k_0 \int_{f \in \mathcal{E}} \Big(\int_{v=0}^{\zeta(f)} H(f(v)) dv \Big) \nu(df).
\end{equation*}
From Theorem \ref{definition Williams decomposition} we see this disintegrates with respect to the amplitude $A$ to give 
\begin{equation}\label{equation williams split}
    k_0 \int_{w \in E \backslash \{0\}} \int_{f \in \mathcal{E}} \Big(\int_{v=0}^{\zeta(f)} H(f(v))dv \Big)\nu(df|A=w) \nu(A \in dw).
\end{equation}
As shown in the left graph of Figure \ref{figure diffusion hit level c}, for $w \in (0,c)$ the Williams decomposition result gives a disintegration with the conditional distribution $\nu(\cdot|A=w)$ being the concatenation of two $\uparrow$-diffusions both stopped at $T_w$. Therefore for the $w \in (0,c)$ in \eqref{equation williams split} we have that
\begin{equation*}
    \int_{w=0}^c \int_{f \in \mathcal{E}} \Big(\int_{v=0}^{\zeta(f)} H(f(v))dv \Big)\nu(df|A=w) \nu(A \in dw) = \int_{w=0}^c 2\mathbb{E}_{\mathbb{Q}^\uparrow_0}\Big(\int_{t=0}^{T_w}H(Y_t)dt\Big) \frac{s(dw)}{s(w)^2}.
\end{equation*}
which is equal to
\begin{equation*}
     \int_{w=0}^c \int_{y=0}^w 2 H(y) (s^\uparrow(w)-s^\uparrow(v))  M^\uparrow(dy) \frac{s(dw)}{s(w)^2}
\end{equation*}
by Theorem \ref{theorem Green's function general} (where we have taken monotonic limits $a \downarrow 0$ and $x \downarrow 0$ as we did in Lemma \ref{theorem diffusion started at c}). With the equations for $s^\uparrow$ and $M^\uparrow$ measure from \eqref{equation s',s,m uparrow}, we substitute these terms and swap the order of integration to obtain
\begin{equation*}
    2 \int_{y=0}^c \int_{w=y}^c \left(\frac{1}{s(y)} - \frac{1}{s(w)}\right) \frac{1}{s(w)^2} s(dw) H(y) s(y)^2 M(dy).
\end{equation*}
We split into cases depending on whether or not $s(c)$ is finite. If $s(c)$ is infinite then $E=[0,c)$ and the above integral is equal to $\int_{y \in E} H(y)M(dy)$. When $s(c)$ is finite we instead have that $E=[0,c]$ and the above integral is equal to
\begin{equation*}
    \int_{y \in E} \Big(1-2\frac{s(y)}{s(c)} + \frac{s(y)^2}{s(c)^2}\Big) H(y) M(dy).
\end{equation*}
Furthermore when $s(c)$ is finite we have that the measure $\nu(A \in \cdot)$ has an atom at $c$ and so there is a positive contribution to the integral in \eqref{equation williams split} at $w=c$. We can calculate this from the Williams decomposition in Theorem \ref{definition Williams decomposition} to obtain
\begin{equation*}
    \nu(A=c)\Big(\mathbb{E}_{\mathbb{Q}^\uparrow_0}\Big(\int_{y=0}^{T_c} H(Y_s)ds\Big) + \mathbb{E}_{\mathbb{Q}_c}\Big( \int_{s=0}^{T_0} H(Y_s) ds \Big)\Big)\Big).
\end{equation*}
We have that $\nu(A=c)$ is equal to $1/s(c)$ and 
\begin{equation*}
    \begin{split}
        \mathbb{E}_{\mathbb{Q}^\uparrow_0}\Big(\int_{y=0}^{T_c} H(Y_s)ds\Big) + \mathbb{E}_{\mathbb{Q}_c}\Big(\int_{s=0}^{T_0} H(Y_s) ds \Big) & = \int_{y \in E} H(y) \Big(\frac{1}{s(y)}-\frac{1}{s(c)}\Big)s(y)^2 M(dy) + \int_{y \in E} H(y) s(y) M(dy) \\
        & = \int_{y \in E} H(y) \Big(2s(y) - \frac{s(y)^2}{s(c)}\Big) M(dy).
    \end{split}
\end{equation*}
Therefore we have \eqref{equation williams split} is equal to $\int_{y \in E} H(y)M(dy)$ in both cases as required.
\end{proof}
\end{proposition}
With Proposition \ref{proposition crossing width density}, we can prove Theorem \ref{theorem aggregate mass summability}. 
\begin{proof}[Proof of Theorem \ref{theorem aggregate mass summability}]
Let $(\mathbf{N},\mathbf{X})$ be a scaffolding-and-spindles pair stopped at an a.s. finite $T$. In the setting where $\nu_\mathbf{Z}$ satisfies \eqref{equation levy measure bounded variation} for $\nu=\nu_\mathbf{Z}$, the scaffolding $\mathbf{X}$ is an SPLP of bounded variation. Therefore the set $\{r \in [0,t]:\mathbf{X}_{r-}<0\leq \mathbf{X}_r\}$ is finite a.s., i.e. there are only finitely many crossing widths to sum; therefore, the aggregate health $M^0_{\mathbf{N},\mathbf{X}}(T)$ at level $0$ is finite as it is a finite sum of positive values. 
\par
In the setting where $\nu_\mathbf{Z}$ satisfies \eqref{equation levy measure unbounded variation} for $\nu=\nu_\mathbf{Z}$ we can apply Proposition \ref{proposition crossing width density}. This proposition states that, for fixed $t>0$, the aggregate fitness process $M^0_{\mathbf{N},\mathbf{X}}(t)$ is a sum over a Poisson random measure with intensity $\text{Leb} \otimes M_\mathbf{Z}$ on $[0,t] \times E \backslash \{0\}$. Let $b\coloneqq 1\wedge c/2 \in (0,c)$. We deal with the sum of crossing widths less than $b$ and the crossing widths of size greater than $b$ separately. The crossing widths of size less than $b$ form a Poisson random measure with intensity $\text{Leb} \otimes M_\mathbf{Z}$ on $[0,t] \times (0,b)$. Therefore, by Campbell's Theorem, see \cite[Section 3.2]{kingman1992poisson}, these are summable if and only if
\begin{equation*}
    \int_{[0,t] \times (0,b)} z dt M_\mathbf{Z}(dz) = t\int_{z=0}^bzM_\mathbf{Z}(dz)<\infty
\end{equation*}
and this is finite by the Assumption in the theorem. The crossing widths of size greater than $b$ are from a Poisson random measure with intensity $\text{Leb} \otimes M_\mathbf{Z}$ on $[0,t] \times [b,c]$. Note that
\begin{equation*}
    \int_{[0,t] \times [b,c]} dt M_\mathbf{Z}(dz) = tM_\mathbf{Z}([b,c])
\end{equation*}
and this is finite by $A2$ and therefore there are only finitely many such crossings so the sum is clearly finite. The result then follows for the aggregate fitness process evaluated at the a.s. finite time $T$.
\par
To see the result in terms of $\mathbf{Y}$, apply Proposition \ref{proposition crossing width density} to diffusion $\mathbf{Y}$ to obtain that the crossing widths in this case are a Poisson random measure with intensity $\text{Leb} \otimes M_\mathbf{Y}$. We apply the transformation $g$ to these crossing widths and this is equal in distribution to the crossing widths from the diffusion $\mathbf{Z}$. The condition for summability of the small crossing widths then becomes the condition in \eqref{equation IP existence condition}.
\end{proof}
\subsection{Disintegration with respect to the crossing width}\label{subsection disintegration proof}
In this section we disintegrate the bi-clade measure $\nu_{\text{cld}}$ with respect to the crossing width $m^0$. As in the previous section, we assume that our diffusion satisfies $A1$ and $A2$ and that the Pitman-Yor excursion measure $\nu$ satisfies \eqref{equation levy measure unbounded variation}. We define some operators that reverses in time excursions and point measures.
\begin{definition}
Let $f \in \mathcal{E}$ and $N \in \mathcal{N}^{\text{sp}}_{\text{fin}}$. Then
\begin{enumerate}
    \item We define the time-reversal for $f$ by
\begin{equation}\label{equation spindle reversal}
    \mathcal{R}_{\text{spdl}} (f) \coloneqq  (f((\zeta(f)-y)-), y \in \mathbb{R}),
\end{equation}
following the definition given in \cite[Equation (2.16)]{construction}. 
\item We also define the time-reversal for $N$ by
\begin{equation}\label{equation clade reversal}
    \mathcal{R}_{\text{cld}} (N) \coloneqq  \sum_{(t,f) \text{   atom of   } N} \delta(\text{len}(N) - t, \mathcal{R}_{\text{spdl}}(f)),
\end{equation}
following the definition given in \cite[Equation (4.19)]{construction}.
\end{enumerate}
\end{definition}
The reversibility of the excursions under the Pitman-Yor excursion measure $\nu$, along with the independence between excursions away from $0$ for the L\'evy process $\mathbf{X}$, means that 
\begin{equation}\label{equation reversal under nu clade}
    \nu_{\text{cld}}(A) = \nu_{\text{cld}}(\mathcal{R}_{\text{cld}}(A)) \text{   for a measurable   } A \subset \mathcal{N}^{\text{sp}}_{\pm \text{cld}}.
\end{equation}
The next result is a generalisation of \cite[Lemma 4.13]{construction}, and it proves conditional independence between the clade and the anti-clade given the crossing width.
\begin{lemma}[Mid-spindle Markov property]\label{lemma mid-spindle}
Consider a diffusion that satisfies $A1$ and $A2$ with Pitman-Yor excursion measure $\nu$ that satisfies \eqref{equation levy measure unbounded variation}, and let $(\mathbf{N},\mathbf{X})$ be a scaffolding-and-spindles construction for this diffusion. Let $T$ be a stopping time either of the form $T^{\geq y}$ for some $y>0$ or of the form
\begin{equation}\label{equation mid-spindle T}
    T^y_{n,j} \coloneqq  \inf \Big\{ t>0 \Big| \int_{[0,t]\times \mathcal{E}} \mathbbm{1}_{f(y-\mathbf{X}(s-))>1/n} d\mathbf{N}(s,f) \geq j \Big\} \text{   for   } y \in \mathbb{R}, n \in \mathbb{N},
\end{equation}
and let $f_T$ denote spindle of $\mathbf{N}$ at this stopping time $T$. Then, given $f_T(-\mathbf{X}(T-))=x>1/n$, we have that 
\begin{equation*}
    (\mathbf{N}|_{[0,T)}, \Check{f}_T) \text{   is conditionally independent of   } (\mathbf{N}|^\leftarrow_{(T,\infty)}, \Hat{f}_T)
\end{equation*}
(see Definition \ref{definition spaces of N} for $\check{f}_T$ and $\hat{f}_T$). Under this conditional law, $\mathbf{N}|^\leftarrow_{(T,\infty)}$ is a Poisson random measure with intensity measure $\text{Leb} \otimes \nu$ independent of $\Hat{f}_T$, which is a $0$-diffusion started at $x$. 
\begin{proof}
The proof given in \cite[Lemma 4.13]{construction} applies in our setting as we have the following conditions: the Pitman-Yor description of an excursion that hits level $1/n$ is given by a $\uparrow$-diffusion started at $0$ and stopped at $1/n$, followed by an independent $0$-diffusion started at $1/n$; the Poisson random measure $\mathbf{N}$ has a strong Markov property; the scaffolding process $\mathbf{X}$ is of unbounded variation and so we have $T^{\geq y} = \inf_{n \geq 1} T^y_{n,1}$. 
\end{proof}
\end{lemma}
For $x>0$, let $\mathbb{Q}^0_x$ be the probability measure on $\mathcal{D}$ for the $0$-diffusion started from $x$.
\begin{lemma}\label{lemma joint probability undershoot overshoot}
    In the setting of Lemma \ref{lemma mid-spindle}, for $n>1$ and $x>1/n$ we have that
\begin{equation*}
    \begin{split}
        & \mathbf{P}\Big(\hat{f}_{T^0_{n,1}} \in A, \check{f}_{T^0_{n,1}} \in B, \mathbf{N}|_{(\tau^0(\ell^0(T^0_{n,1})-),T^0_{n,1})} \in C, \mathbf{N}|^\leftarrow_{(T^0_{n,1},\tau^0(\ell^0(T^0_{n,1})))} \in D\Big| f_{T^0_{n,1}}(-\mathbf{X}(T^0_{n,1}-))=x\Big) \\
        & = \int_{f \in A} \mathbf{P}(\mathbf{N}|_{[0,T^{-\zeta(f))})} \in D) \mathbb{Q}^0_x(df) \int_{f' \in \mathcal{R}_{\text{spdl}}(B)} \mathbf{P}(\mathbf{N}|_{[0,T^{-\zeta(f'))})} \in \mathcal{R}_{\text{cld}}(C)) \mathbb{Q}^0_x(df'). \\
    \end{split}
\end{equation*}
    \begin{proof}
First note by the conditional independence property given in Lemma \ref{lemma mid-spindle} we have, for $x>0$, that
\begin{equation*}
    \mathbf{P}\Big(\hat{f}_{T^0_{n,1}} \in A, \check{f}_{T^0_{n,1}} \in B\Big| f_{T^0_{n,1}}(-\mathbf{X}(T^0_{n,1}-))=x\Big) = \mathbb{Q}_x^0(A)\mathbf{P}\Big(\check{f}_{T^0_{n,1}} \in B\Big| f_{T^0_{n,1}}(-\mathbf{X}(T^0_{n,1}-))=x\Big).
\end{equation*}
Let $\lambda > 0$ and let $e_\lambda$ be an exponential random variable with rate $\lambda$. The excursions of $\mathbf{X}$ away from $0$ are indexed by its local time. Therefore $\mathbf{N}|_{[0,\tau^0(e_\lambda)]}$ consists of a collection of complete bi-clades. Therefore, by reversibility of the clade measure, noted in \eqref{equation reversal under nu clade}, we have
\begin{equation}\label{equation equality in distributions}
    \mathbf{N}|_{[0,\tau^0(e_\lambda)]} \overset{d}{=} \mathcal{R}_{cld}\Big(\mathbf{N}|_{[0,\tau^0(e_\lambda)]}\Big).
\end{equation}
Now, on the event $\{T^0_{n,1} < \tau^0(e_\lambda)\}$ there exists $j \in \mathbb{N}$ such that 
\begin{equation}\label{equation events relabel}
    T^0_{n,1}(\mathbf{N}|_{[0,e_\lambda]}) = \tau^0(e_\lambda) - T^0_{n,j}(\mathcal{R}_{\text{cld}}(\mathbf{N}|_{[0,e_\lambda]})).
\end{equation}
Furthermore, we have that
\begin{equation*}
    \mathbf{P}\Big(\check{f}_{T^0_{n,1}} \in B\Big| f_{T^0_{n,1}}(-\mathbf{X}(T^0_{n,1}-))=x\Big) = \lim_{\lambda \downarrow 0} \sum_{j=1}^\infty \mathbf{P}\Big(\check{f}_{T^0_{n,1}} \in B, T^0_{n,j}<\tau^0(e_\lambda)<T^0_{n,j+1}\Big| f_{T^0_{n,1}}(-\mathbf{X}(T^0_{n,1}-))=x\Big)
\end{equation*}
as $\lim_{\lambda \downarrow 0} \mathbf{P}(T^0_{n,1} < \tau^0(e_\lambda)) = 1$. We consider each term in the sum. We apply the time reversal and the mid-spindle Markov property (Lemma \ref{lemma mid-spindle}) to obtain
\begin{equation*}
    \mathbf{P}\Big(\check{f}_{T^0_{n,1}} \in B\Big| T^0_{n,j}<\tau^0(e_\lambda)<T^0_{n,j+1}, f_{T^0_{n,1}}(-\mathbf{X}(T^0_{n,1}-))=x\Big) = \mathbf{P}\Big(\hat{f}_{T^0_{n,j}} \in \mathcal{R}_{\text{spdl}}(B)\Big| f_{T^0_{n,j}}(-\mathbf{X}(T^0_{n,j}-))=x\Big)
\end{equation*}
where the conditioning $T^0_{n,j}<\tau^0(e_\lambda)<T^0_{n,j+1}$ was able to be removed due to the conditional independence given $f_{T^0_{n,1}}(-\mathbf{X}(T^0_{n,1}-))=x$ from the mid-spindle Markov property. Let $p_\lambda=\mathbf{P}(\tau^0(e_\lambda)>T^0_{n,1})$. Then
\begin{equation*}
    \mathbf{P}\Big(T^0_{n,j}<\tau^0(e_\lambda)<T^0_{n,j+1}\Big| f_{T^0_{n,1}}(-\mathbf{X}(T^0_{n,1}-))=x\Big) = \mathbf{P}\Big(T^0_{n,j}<\tau^0(e_\lambda)<T^0_{n,j+1}\Big) = p_\lambda^j(1-p_\lambda).
\end{equation*}
Pulling all this together we obtain
\begin{equation*}
    \lim_{\lambda \downarrow 0} \mathbf{P}\Big(\check{f}_{T^0_{n,1}} \in B, T^0_{n,j}<\tau^0(e_\lambda)<T^0_{n,j+1}\Big| f_{T^0_{n,1}}(-\mathbf{X}(T^0_{n,1}-))=x\Big) = \lim_{\lambda \downarrow 0} \sum_{j=1}^\infty p_\lambda^j(1-p_\lambda)\mathbb{Q}^0_x(\mathcal{R}_{\text{spdl}}(B)),
\end{equation*}
which equals $\mathbb{Q}^0_x(\mathcal{R}_{\text{spdl}}(B))$ because $p_\lambda \rightarrow 1$ as $\lambda \rightarrow 0$. The generalisation to the statement in the proposition follows from the Lemma \ref{lemma mid-spindle}; the strong Markov property of $\mathbf{N}$ at $T^0_{n,1}$ leads to the independence of the Poisson random measures $\mathbf{N}|_{[0,T^0_{n,1})}$ and $\mathbf{N}|^\leftarrow_{(T^0_{n,1},\infty)}$.
\end{proof}
\end{lemma}
We now describe the regular conditional distribution of $\nu_{\text{cld}}$ disintegrated with respect to $m^0$. Let $x>0$ denote the crossing width. Let $\hat{f}$ and $\check{f}$ be independent $0$-diffusions started at $x$. Let $\check{\mathbf{N}}$ and $\hat{\mathbf{N}}$ be independent Poisson random measures with intensity $\text{Leb} \otimes \nu$.
\par
We define the idea of concatenating point measures. Let $N_i \in \mathcal{N}^{\text{sp}}_{\text{fin}}$ be a collection of point measures where $i \in I$ is an ordered index set. We define $N\coloneqq \bigstar_{i \in I} N_i$ to be the concatenation with respect to the ordering of $I$ in the following way. On the condition that for all $i \in I$ we have $S(j-)\coloneqq \sum_{j \in I: j<i} \text{len}(N_j)<\infty$ we take $N$ to be the point measure with length $\sum_{i \in I}\text{len}(N_i) \leq \infty$ such that $N|^\leftarrow_{(S(i-), S(i-)+\text{len}(N_i))}=N_i$. If the condition does not hold we take $N=\emptyset$. When $I=\{1,2, \dots, n\}$ is a finite set we write $N_1 \star \dots \star N_n$ for $\bigstar_{i \in I} N_i$.
\par
We define a function to stick spindles together: for $f_1$, $f_2$ excursions we define $f_{1 \rightarrow 2}$ by
\begin{equation*}
    f_{1 \rightarrow 2}(y) \coloneqq  \min \left\{
        \begin{matrix}
            & f_1(\zeta(f_1)-y) & \text{   for   } & 0<y<\zeta(f_1) \\
            & f_2(y-\zeta(f_1)) & \text{   for   } & \zeta(f_1) \leq y \leq \zeta(f_1) + \zeta(f_2)  \\
            & 0 & \text{   otherwise.} \\
        \end{matrix}
    \right.
\end{equation*}
Define the map for paths $f_1$, $f_2$ and for $N_1$, $N_2 \in \mathcal{N}^{\text{sp}}$
\begin{equation*}
    H(f_1, f_2, N_1, N_2)\coloneqq \mathcal{R}_{\text{cld}}(N_1|_{[0,T^{-\zeta(f_1)})}) \star \delta_{(0,f_{1 \rightarrow 2})} \star N_2|_{[0,T^{-\zeta(f_2)})}
\end{equation*}
and define the measure $\mathbf{M}$ by
\begin{equation*}
    \mathbf{M}\coloneqq  H(\check{f}, \hat{f}, \check{\mathbf{N}}, \hat{\mathbf{N}}),
\end{equation*}
and let $\nu_x$ be the distribution of $\mathbf{M}$ on $\mathcal{N}^{\text{sp}}_{\text{cld}}$. We note from the independence in the construction that for $x>0$ and $G_1 \subseteq \mathcal{N}^{\text{sp}}_{+\text{cld}}$, $G_2 \subseteq \mathcal{N}^{\text{sp}}_{-\text{cld}}$ that 
\begin{equation*}
    \nu_x(N^+ \in G_1, N^- \in G_2) = \nu_x(N^+ \in G_1)\nu_x(N^- \in G_2).
\end{equation*}
Furthermore, let $\hat{f}$ be a $0$-diffusion started at $x$, independent of $\mathbf{N}$. Let $\hat{T}^0\coloneqq \inf\{t>0:\xi_\mathbf{N}(t) = -\zeta(\hat{f})\}$. Then
\begin{equation*}
    \overline{N}^+_x\coloneqq \delta(0,\hat{f})+\mathbf{N}|_{[0,\hat{T}^0)}
\end{equation*}
has the law $\nu_x(N^+ \in \cdot)$, and $\mathcal{R}_{\text{cld}}(\overline{N}^+_x)$ has the law $\nu_x(N^- \in \cdot)$.
\begin{proposition}\label{proposition disintegration}
Consider a diffusion that satisfies $A1$ and $A2$ with Pitman-Yor excursion measure $\nu$ that satisfies \eqref{equation levy measure unbounded variation}. Then $\nu_x$, $x \in E$, is a disintegration of $\nu_{\text{cld}}$ with respect to $m^0$. 
    \begin{proof}
        This proof concerns kernels as defined in Kallenberg \cite[Chapter I.3]{Kallenberg_2021}. Let $\mu_x$ be the distribution of $(\check{f}, \hat{f}, \check{\mathbf{N}}, \hat{\mathbf{N}})$. We first show that $x \mapsto \mu_x$ is a kernel. First note that $x \mapsto \mathbf{P}(\mathbf{N} \in \cdot)$ is a kernel and also have that $x \mapsto \mathbb{Q}^0_x$ is a kernel. The first follows as the image is independent of $x$ and the second follows as the measure $\mathbb{Q}^0_x$ is continuous in the weak sense and is therefore measurable. It therefore follows that $x \mapsto \mu_x$ is a kernel as we have that $\mu_x=\mathbb{Q}^0_x \otimes \mathbb{Q}^0_x \otimes \mathbf{P}(\mathbf{N} \in \cdot) \otimes \mathbf{P}(\mathbf{N} \in \cdot)$, see \cite[Lemma I.3.3]{Kallenberg_2021}. The result then follows as $H$ is a measurable map and $\nu_x$ is $\mu_x(H^{-1}(\cdot))$, and therefore map $x \mapsto \nu_x$ is a measurable map and therefore is a kernel, see \cite[Lemma I.3.1]{Kallenberg_2021}. Therefore the map on $\mathbb{R}_+ \rightarrow \mathcal{N}_{\text{cld}}^{\text{sp}}$ given by $x \mapsto \nu_x$ is a kernel because the map $H$ from $(\check{f}, \hat{f}, \check{\mathbf{N}}, \hat{\mathbf{N}})$ to $\mathbf{M}$ is measurable.
        \par
        We now verify that $(\nu_x, x > 0)$ do integrate to the bi-clade measure $\nu_{\text{cld}}$. Let $A$, $B \subset \mathcal{E}$ measurable and let $C$, $D \subset \mathcal{N}^{\text{sp}}_{\text{fin}}$ measurable. By definition of $\nu_{\text{cld}}$ we have that $\nu_{\text{cld}}(H(A \times B \times C \times D) \cap \{m^0>1/n\})$ equals
        \begin{equation*}
            \begin{split}
                & \nu_{\text{cld}}(\hat{f}_{T^+_0} \in A, \check{f}_{T^+_0} \in B, N^-|_{[0,\text{len}(N^-))} \in C, N^+|_{[0,\text{len}(N^+))} \in D, m^0>1/n) \\
                & = k_0 M(E \backslash [0,1/n]) \mathbf{P}(\hat{f}_{T^0_{n,1}} \in A, \check{f}_{T^0_{n,1}} \in B, \mathbf{N}|_{[\tau^0(\ell^0(T^0_{n,1})-),T^0_{n,1})} \in C, \mathbf{N}|^\leftarrow_{(T^0_{n,1},\tau^0(\ell^0(T^0_{n,1})))} \in D). \\
            \end{split}
        \end{equation*}
        With the result from Lemma \ref{lemma joint probability undershoot overshoot} along with the crossing width intensity measure from Proposition \ref{proposition crossing width density} this is equal to
        \begin{equation*}
            \int_{x \in E \backslash [0,1/n]} \int_{f \in \mathcal{R}_{\text{spdl}}(B)} \mathbf{P}(\mathbf{N}|_{[0,T^{-\zeta(f))})} \in \mathcal{R}_{\text{cld}}(C)) \mathbb{Q}^0_x(df) \int_{f' \in A} \mathbf{P}(\mathbf{N}|_{[0,T^{-\zeta(f')})} \in D) \mathbb{Q}^0_x(df')k_0M(dx).
        \end{equation*}
        We take the limit as $n$ tends to infinity to obtain
        \begin{equation*}
            \begin{split}
                & \nu_{\text{cld}}(H(A \times B \times C \times D)) = \nu_{\text{cld}}(\hat{f}_{T^+_0} \in A, \check{f}_{T^+_0} \in B, N^-|_{[0,\text{len}(N^-))} \in C, N^+|_{[0,\text{len}(N^+))} \in D) \\
                & = \int_{x \in E \backslash \{0\}} \int_{f \in \mathcal{R}_{\text{spdl}}(B)} \mathbf{P}(\mathbf{N}|_{[0,T^{-\zeta(f))})} \in \mathcal{R}_{\text{cld}}(C)) \mathbb{Q}^0_x(df) \int_{f' \in A} \mathbf{P}(\mathbf{N}|_{[0,T^{-\zeta(f'))})} \in D) \mathbb{Q}^0_x(df')2k_0M(dx). \\
            \end{split}
        \end{equation*}
        By the definition of $\nu_x$ and $\mu_x$ this is equal to 
        \begin{equation*}
            \int_{x \in E \backslash \{0\}} \mu_x(A \times B \times C \times D) \nu_{\text{cld}}(m^0 \in dx) = \int_{x \in E \backslash \{0\}} \nu_x(H(A \times B \times C \times D)) \nu_{\text{cld}}(m^0 \in dx).
        \end{equation*}
        The generality of $A$, $B \in \mathcal{E}$, and $C$, $D \in \mathcal{N}^{\text{sp}}_{\text{fin}}$ is sufficient here to prove the proposition by the monotone class theorem.
    \end{proof}
\end{proposition}
We also note from the proof of the proposition above that $\mu_x$, $x \in E \backslash \{0\}$, is in fact a disintegration of $\nu_{\text{cld}}(\check{f}_{T^+_0} \in \cdot, \hat{f}_{T^+_0} \in \cdot, N^-|_{[0,\text{len}(N^-))} \in \cdot, N^+|_{[0,\text{len}(N^+))} \in \cdot)$.
\subsection{L\'evy process scaffolding from fixed initial interval partition}\label{section space of interval partitions}
In this section we give the construction of scaffolding-and-spindles pairs started from certain interval partitions $\beta \in \mathcal{I}_H$, which we denote as $(\mathbf{N}_\beta, \mathbf{X}_\beta)$. We give the construction of $\mathbf{N}_\beta$ in Definition \ref{definition starting interval partition} and we prove that this construction is possible under certain conditions in Proposition \ref{proposition starting interval partition}. We define random measures $\mathbf{N}_\beta$ on the space $\mathcal{N}^{\text{sp}}_{\text{fin}}$. 
\begin{definition}\label{definition starting interval partition}
    Let $\mathcal{I}_H(E) \coloneqq \{\gamma \in \mathcal{I}_H: \text{Leb}(U) \in E \hspace{0.125cm} \forall U \in \gamma\}$ and let $\beta \in \mathcal{I}_H(E)$. If $\beta = \varnothing$ then define $\mathbf{N}_\beta$ to be the zero measure. Otherwise, for each interval $U$ of $\beta$ we construct each $\mathbf{N}_U$ independently in the following way. First let $f_U$ have law $\mathbb{Q}^0_{\text{Leb}(U)}$ and let $\mathbf{N}$ be a Poisson random measure with intensity $\text{Leb} \otimes \nu$. Let $T\coloneqq  \inf\{t>0: \xi_\mathbf{N}(t) = -\zeta(f_U)\}$ and let $\mathbf{N}_U \coloneqq  \delta(0,f_U) + \mathbf{N}|_{[0,T]}$.
    \par
    Under the condition that $\sum_{U \in \beta} \text{len}(\mathbf{N}_U)$ is finite, let $\mathbf{N}_\beta \coloneqq \bigstar_{U \in \beta} \mathbf{N}_U$ (and otherwise take $\mathbf{N}_\beta$ to be the zero measure). Let $\mathbf{X}_\beta \coloneqq \xi(\mathbf{N}_\beta)$ be the corresponding scaffolding constructed from $\mathbf{N}_\beta$. We write $\mathbf{P}_\beta$ to denote the law of $\mathbf{N}_\beta$ on $\mathcal{N}^{\text{sp}}_{\text{fin}}$. 
\end{definition}
In all the results that follow we consider intervals partitions $\beta \in \mathcal{I}_H(E)$. We give sufficient conditions on an interval partition $\beta$ for concatenating clades $(\mathbf{N}_U, U \in \beta) \subset \mathcal{N}^{\text{sp}}_{+ \text{cld}}$ as in Definition \ref{definition starting interval partition}.
\begin{proposition}\label{proposition starting interval partition}
    Consider a diffusion that satisfies $A1$, $A2$ with a Pitman-Yor excursion measure $\nu$ that satisfies \eqref{equation levy measure unbounded variation}. Further suppose that
    \begin{equation}\label{equation starting IP condition}
        \limsup_{x \downarrow 0} \frac{1}{x} \Big(\int_{v \in E} s(v \wedge x) M(dv)\Big)=:k_3<\infty.
    \end{equation}
    Then, in the setting of Definition \ref{definition starting interval partition}, for any interval partition $\beta \in \mathcal{I}_H(E)$, we have that $\sum_{U \in \beta} \text{len}(\mathbf{N}_U)$ is finite a.s. and hence $(\mathbf{N}_\beta, \mathbf{X}_\beta)$ is a non-trival scaffolding-and-spindles pair (unless $\beta=\emptyset$).
    \begin{proof}
    We calculate the expected lifetime of an interval with initial width $x$, for $x\in E$, and we define the function $r(x)\coloneqq \mathbb{E}_{\mathbb{Q}^0_x}(\zeta)$ to be this expectation. We first assume that $s(c)$ is finite and we calculate $\mathbb{E}_{\mathbb{Q}^0_x}(T_0 \wedge T_c)$. For $0<a<x<w<c$ we use $\mathbb{E}_{\mathbb{Q}^0_x}(T_a \wedge T_w)$ from \eqref{equation expected hitting time} and we can take the monotonic limits of $a \downarrow 0$ and $w \uparrow c$ to find 
    \begin{equation*}
        \mathbb{E}_{\mathbb{Q}_x}(T_0 \wedge T_c) = \int_{v \in E} \Big(\frac{s(x)}{s(c)} \Big(s(c)-s(v)\Big) \mathbbm{1}_{x<v} + \frac{s(v)}{s(c)} \Big(s(c)-s(x)\Big) \mathbbm{1}_{v \leq x}\Big) M(dv).
    \end{equation*}
    To account for the path from $c$ to $0$ on the event that $T_c<T_0$ we have
    \begin{equation*}
        \mathbb{E}_{\mathbb{Q}_x}\Big(T_0-T_c:T_c<T_0\Big) = \mathbb{Q}_x(T_c<T_0) \mathbb{E}_{\mathbb{Q}_c}(T_0).
    \end{equation*}
    Firstly we have $\mathbb{Q}_x(T_c<T_0)=s(x)/s(c)$ (see \eqref{equation scale function property}). The second term can be calculated from Lemma \ref{theorem diffusion started at c} to give
    \begin{equation*}
        \mathbb{E}_{\mathbb{Q}_x}\Big(T_0-T_c:T_c<T_0\Big) = \frac{s(x)}{s(c)} \int_{v \in E} s(v)M(dv).
    \end{equation*}
    We add the two contributions to obtain
    \begin{equation}\label{equation r definition}
        r(x) = \int_{v \in E} s(v \wedge x) M(dv).
    \end{equation}
    In the case where $s(c)$ is infinite we use $\mathbb{E}_{\mathbb{Q}^0_x}(T_a \wedge T_w)$ from \eqref{equation expected hitting time} and we can take the monotonic limits of $a \downarrow 0$ and $w \uparrow c$ and this equals $\mathbb{E}_{\mathbb{Q}^0_x}(T_0)$ and this limit is as in \eqref{equation r definition}.
    \par
    Under the assumption given in \eqref{equation starting IP condition}, there exists $\varepsilon'>0$ such that for all $x<\varepsilon'$ we have that $r(x) \leq (1+k_3) x$. Therefore, for an interval partition $\beta \in \mathcal{I}_H(E)$ we have
    \begin{equation}\label{equation sum of starting diffusion lifetimes}
        \mathbb{E}\Big(\sum_{\substack{U \in \beta \\ \text{Leb}(U) < \varepsilon'}} \zeta(f_U)\Big) = \sum_{\substack{U \in \beta \\ \text{Leb}(U) < \varepsilon'}} r(\text{Leb}(U)) \leq \sum_{\substack{U \in \beta \\ \text{Leb}(U) < \varepsilon'}} (1+k_3) \text{Leb}(U) = (1+k_3) \norm{\beta} < \infty.
    \end{equation}
    As there are only finitely many intervals $U$ of the interval partition $\beta$ that have length greater than $\varepsilon'$ the result follows.
    \end{proof}
\end{proposition}
We show in Corollary \ref{corollary assumption b' starting interval partition} that \eqref{equation starting IP condition} is satisfied for a large class of diffusions and in Remark \ref{remark squared bessel dim 0} we highlight a class of diffusions that satisfy $A1$, $A2$ and \eqref{equation levy measure unbounded variation} but that do not satisfy \eqref{equation starting IP condition}.
\begin{corollary}\label{corollary assumption b' starting interval partition}
    Consider a diffusion that satisfies $A1$, $A2$, $A3$ with Pitman-Yor excursion measure that satisfies \eqref{equation levy measure unbounded variation} and suppose that there exists $k_4<0$ such that 
    \begin{equation*}
        \limsup_{y \downarrow 0} \mu(y)= k_4<0.
    \end{equation*}
    Then the diffusion satisfies \eqref{equation starting IP condition} and therefore for any interval partition $\beta \in \mathcal{I}_H(E)$ the construction of $(\mathbf{N}_\beta,\mathbf{X}_\beta)$ as in Definition \ref{definition starting interval partition} exists and $\text{len}(\mathbf{N}_\beta)$ is finite a.s. Furthermore we have that the function $r$ defined in Proposition \ref{proposition starting interval partition} has a continuous derivative $r'$ on $(0,c)$ and $\limsup_{x \rightarrow 0}r'(x)<\infty$.
    \begin{proof}
        As our diffusion satisfies $A3$ we can use the equations for the scale function derivative and speed measure density from \eqref{equation scale and speed measure}. This means that we can calculate $r'(x)$ and note that the condition in \eqref{equation starting IP condition} is satisfied if $\limsup_{x \rightarrow 0} r'(x)$ is finite because this is greater than or equal to $\limsup_{x \rightarrow 0} r(x)/x$. We obtain
        \begin{equation*}
            r'(x) = s'(x) \int_{v=x}^c m(v)dv
        \end{equation*}
    and this function is continuous on $(0,c)$. By the condition given in the corollary, there exists $\varepsilon>0$ such that for $y<\varepsilon$ we have $\mu(y)<k_4/2$. Now, for $x<\varepsilon$, we have that
    \begin{equation*}
        s'(x) \int_{v=x}^\varepsilon m(v) dv = \int_{v=x}^\varepsilon \frac{2}{\sigma^2(v)} \exp \Big( \int_{z=x}^v \frac{2\mu(z)}{\sigma^2(z)} dz \Big) dv \leq \int_{v=x}^\varepsilon \frac{2}{\sigma^2(v)} \exp \Big( \int_{z=x}^v \frac{k_4}{\sigma^2(z)} dz \Big) dv
    \end{equation*}
    and therefore
    \begin{equation}\label{equation r' upper bound}
        \begin{split}
            \limsup_{x \downarrow 0} s'(x) \int_{v=x}^\varepsilon m(v) dv
            & \leq \limsup_{x \downarrow 0}     \int_{v=x}^\varepsilon \frac{2}{\sigma^2(v)} \exp \Big( \int_{z=x}^v \frac{k_4}{\sigma^2(z)} dz \Big) dv \\
            & = \limsup_{x \downarrow 0} \frac{-2}{k_4} \Big( 1 - \exp \Big( \int_{z=x}^\varepsilon \frac{k_4}{\sigma^2(z)} dz \Big) \Big) = \frac{-2}{k_4} < \infty. \\
        \end{split}
    \end{equation}
    It remains to show that $\limsup_{x \downarrow 0} s'(x) \int_{v=\varepsilon}^c m(v) dv$ is finite. We have that $\int_{v=\varepsilon}^c m(v) dv$ is finite by $A2$. Also as $\mu$ is negative locally at $0$ we have that the limit $\lim_{x \downarrow 0} s'(x)$ exists and equals
    \begin{equation*}
        s'(0)\coloneqq \exp\left(\int_{z=0}^b \frac{2 \mu (z)}{\sigma^2(z)}dz\right)
    \end{equation*}
    where $b \in (0,c)$ is from \eqref{equation scale and speed measure} as the function $s'(x)$ is monotonic decreasing for $x\in (0,\varepsilon)$. Therefore we have an upper bound for $\limsup_{x \rightarrow 0} r(x)/x$ and for $\limsup_{x \rightarrow 0}r'(x)$ of $-2/k_4 + s'(0)\int_{v=\varepsilon}^c m(v)dv$ and from this we see that the diffusion satisfies \eqref{equation starting IP condition} by Proposition \ref{proposition starting interval partition}.
    \end{proof}
\end{corollary}
We note that the untransformed diffusions under Assumption \ref{assumption b} satisfy the condition in Corollary \ref{corollary assumption b' starting interval partition} and so satisfy \eqref{equation starting IP condition}. This does not extend to all diffusions under Assumption \ref{assumption b} as \eqref{equation starting IP condition} in not invariant under spatial transformations $g$ (unlike our previous conditions with the speed measure and scale function in \eqref{equation 0 is regular} and \eqref{equation x^2 bound}).
\begin{remark}\label{remark squared bessel dim 0}
    The untransformed diffusions in Example \ref{example squared bessel dim 0} satisfy $A1$, $A2$ and \eqref{equation levy measure unbounded variation} but do not satisfy \eqref{equation starting IP condition}: using \eqref{equation speed scale squared bessel dim 0} for  the density $m$ of the speed measure $M$ and for the scale function $s$ we have, for $x \in (0,\varepsilon_2)$, 
    \begin{equation*}
        \frac{1}{x} \Big(\int_{v \in E} s(v \wedge x) m(v) dv \Big) = \frac{1}{x} \Big( \int_{v=0}^x \frac{1}{2} dv + \int_{v=x}^{\varepsilon_2} x \frac{1}{2v} dv + xM((\varepsilon_2, c]) \Big) = \frac{-\log (x)+ 1 + \log(\varepsilon_2)}{2} + M((\varepsilon_2, c])
        \end{equation*}
    and this goes to infinity as $x$ goes to zero.
\end{remark}
The purpose of Proposition \ref{proposition starting interval partition} and Corollary \ref{corollary assumption b' starting interval partition} is to give a framework from which we can develop a simple Markov property for the branching IP-evolutions for Theorem \ref{theorem diffusion}.
For the remainder of propositions and lemmas in this section we assume \eqref{equation starting IP condition} and so that for any interval partition $\beta \in \mathcal{I}_H(E)$ we have a non-trivial scaffolding-and-spindles pair $(\mathbf{N}_\beta, \mathbf{X}_\beta)$. We now define several functions that will allow us to construct suitable filtrations in space and in time, and then prove the simple Markov property in the setting of the Poisson random measure $\mathbf{N}$ on the probability space $(\Omega, \mathcal{A}, \mathbf{P})$. We omit the exact proofs, rather, we give references to results from \cite{construction} and explain how they generalise to our setting. First we define $M$-diversity which is crucial for the arguments.
\begin{definition}\label{definition m-diversity}
Let $\beta$ be an interval partition and $t \geq 0$. Then we define the $M$-diversity of $\beta$ at $t$ as the following limit when it exists:
\begin{equation}\label{equation definition m-diversity}
    \mathcal{D}^M_\beta(t) \coloneqq  \lim_{x \rightarrow 0} \frac{1}{M(E \backslash [0,x])} \Big| \Big\{ (a,b) \in \beta: |b-a| > x, b \leq t \Big\}\Big| \in [0,\infty),
\end{equation}
and if the limit does not exist then we say that $\beta$ at $t$ does not have an $M$-diversity. Furthermore, for an interval $U = (a,b) \in \beta$ we define the limit $\mathcal{D}^M_\beta(U) \coloneqq  \mathcal{D}^M_\beta(b)$, when this limit exists. Finally, we also define the limsup-$M$-diversity by 
\begin{equation}\label{equation definition limsup m-diversity}
    \overline{\mathcal{D}}^M_\beta(t) \coloneqq  \limsup_{x \rightarrow 0} \frac{1}{M(E \backslash [0,x])} \Big| \Big\{ (a,b) \in \beta: |b-a| > x, b \leq t \Big\}\Big| \in [0,\infty],
\end{equation}
and the limit always exists for any interval partition $\beta$, and for an interval $U \in \beta$ we define $\overline{\mathcal{D}}^M_\beta(U)\coloneqq \overline{\mathcal{D}}^M_\beta(b)$. Let $\mathcal{I}_M \subset \overline{\mathcal{I}}_M$ be
\begin{equation}\label{equation interval partition class with m diversity}
    \mathcal{I}_M \coloneqq  \{ \beta \in \mathcal{I}_H: \forall t>0, \mathcal{D}^M_\beta(t) \text{   exists}\}; \hspace{0.5cm} \overline{\mathcal{I}}_M \coloneqq  \{ \beta \in \mathcal{I}_H: \forall t>0, \overline{\mathcal{D}}^M_\beta(t) <\infty \},
\end{equation}
and let $\mathcal{I}_M(E) = \mathcal{I}_H(E) \cap \mathcal{I}_M$ for a measurable subset $E \subseteq [0,\infty)$ and similarly for $\overline{\mathcal{I}}_M(E)$. We substitute the $M$ for $m$ and write $\mathcal{D}^m_\beta$ for $\mathcal{D}^M_\beta$ and similarly for the other terms when the speed measure $M$ has a density $m$, and we call replace the $M$- for $m$- in their names.
\end{definition}
\begin{lemma}\label{lemma diversity is local time}
    Consider a diffusion that satisfies $A1$ and $A2$ of Assumption \ref{assumption a} with speed measure and scale function that satisfies \eqref{equation starting IP condition}, and for which the Pitman-Yor excursion measure $\nu$ satisfies \eqref{equation levy measure unbounded variation}. Let $(\mathbf{N},\mathbf{X})$ be a stopped scaffolding-and-spindles pair for this diffusion, defined on $(\Omega, \mathcal{A}, \mathbf{P})$. Then, for $y \in \mathbb{R}$ and $t>0$, $\beta^y\coloneqq \text{skewer}(y, \mathbf{N}_\beta,\mathbf{X}_\beta)$ has an $M$-diversity a.s. We further have that for all intervals $U \in \beta^y$ simultaneously that $\mathcal{D}^M_{\beta^y}(U)=\ell^y(U)$ a.s. 
    \begin{proof}   
        Let $x>0$. Then, under $\mathbf{P}(\cdot|\mathbf{X})$ (see \cite[Proposition 4.9]{construction} for a description), the number of intervals before $t$ with width greater than $x$ in the interval partition $\beta^y$ is a Poisson random variable with parameter $\ell^y(t)M(E \backslash [0,x])$. Therefore, as $M(E\backslash \{0\})=\infty$ and by the strong law of large numbers, we have the limit $\mathcal{D}^M_\beta(t)$ exists and equals $\ell^y(t)$ a.s.
    \end{proof}
\end{lemma}
\begin{lemma}\label{lemma injective local times}
     Consider a diffusion that satisfies $A1$ and $A2$ of Assumption \ref{assumption a} with speed measure and scale function that satisfies \eqref{equation starting IP condition}, and for which the Pitman-Yor excursion measure $\nu$ satisfies \eqref{equation levy measure unbounded variation}. Let $\beta \in \mathcal{I}_H(E)$ and let $(\mathbf{N}_\beta, \mathbf{X}_\beta)$ be a scaffolding-and-spindles pair started from $\beta$ for this diffusion as in Definition \ref{definition starting interval partition}. For $y>0$ fixed, it is a.s. the case that the local times of the excursions of $\xi(\mathbf{N}_\beta)$ are distinct about level $y$, the crossing widths are summable, and the $M$-diversities exist for all $U \in \beta^y\coloneqq \text{skewer}(y, \mathbf{N}_\beta,\mathbf{X}_\beta)$ and are distinct (in the sense that for $U \neq V \in \beta^y$ we have that $\mathcal{D}^M_{\beta^y}(U) \neq \mathcal{D}^M_{\beta^y}(V)$). 
    \begin{proof} 
        First note that $\xi_{\mathbf{N}_U}$ is a L\'evy process of unbounded variation started from a random level $\zeta(f_U)$, where $f_U$ is the incomplete spindle with initial value $\text{Leb}(U)$, so that $\xi_{\mathbf{N}_U}(0) \sim \mathbb{Q}^0_{\text{Leb}(U)}(\zeta \in \cdot)$. For $y>0$ at most finitely many of the $\xi_{\mathbf{N}_U}$ from $\mathbf{N}_\beta$ intersect level $y$. Therefore, the local times claim follows from Millar \cite{millar1973exit}; see Proposition \ref{proposition properties of levy process of unbounded variation} in the Appendix for a full description. 
        \par
        The summability of the crossing widths at level $y$ follows from Theorem \ref{theorem aggregate mass summability} along with the fact that only finitely many of the $\mathbf{N}_U$ intersect level $y$. The $M$-diversities claim follows from the previous two results of this proposition along with Lemma \ref{lemma diversity is local time}.
    \end{proof}
\end{lemma}
We now define cutoff functions which are used to define the filtrations from \cite{construction}. The idea is to have a filtration indexed by level $y \in \mathbb{R}$ that can be used to define a Markov property of the skewer process. 
\begin{equation*}
    \sigma^y_N(t)\coloneqq  \text{Leb}\{u \leq t: \xi_N(u) \leq y \},
\end{equation*}
\begin{equation}\label{equation cutoff one}
    \begin{split}
        & \text{cutoff}^{\leq y}_{\xi(N)}(s)\coloneqq \xi_N(\sup \{t \geq 0: \sigma^y_N(t) \leq s\}) - \min\{y,0\}, \\
        & \text{cutoff}^{\geq y}_{\xi(N)}(s)\coloneqq \xi_N(\sup \{t \geq 0: t - \sigma^y_N(t) \leq s\}) - \max\{y,0\}, \\
    \end{split}
\end{equation}
\begin{equation}\label{equation cutoff two}
    \begin{split}
        & \text{cutoff}^{\leq y}_N\coloneqq  \sum\limits_{(t,f) \in N} \mathbbm{1}_{\xi_N(t-)<y<\xi_N(t)} \delta(\sigma^y_N(t),\Check{f}^y_t) + \mathbbm{1}_{\xi_N(t) \leq y}\delta(\sigma^y_N(t), f_t), \\
        & \text{cutoff}^{\geq y}_N\coloneqq  \sum\limits_{(t,f) \in N} \mathbbm{1}_{\xi_N(t-)<y<\xi_N(t)} \delta(t-\sigma^y_N(t),\Check{f}^y_t) + \mathbbm{1}_{\xi_N(t-) \geq y}\delta(t-\sigma^y_N(t), f_t), \\
    \end{split}
\end{equation}
\begin{definition}[Filtrations]\label{definition filtrations}
We define a filtration indexed by level $y \geq 0$, $(\overline{\mathcal{F}^y}, y \geq 0)$, and a filtration, indexed by time, $(\overline{\mathcal{F}_t}, t \geq 0)$.
    \begin{enumerate}
        \item We define $(\mathcal{F}^y, y \geq 0)$ to be the least right-continuous filtration under which the maps $N \mapsto \text{cutoff}^{\leq y}_{\xi(N)}$ and $N \mapsto \text{cutoff}^{\leq y}_N$ are $\mathcal{F}^y$-measurable for $y \geq 0$.
        \item We define $(\mathcal{F}_t, t \geq 0)$ to be the least right-continuous filtration on $\mathcal{N}^{\text{sp}}$ under which $N \mapsto N|_{[0,t]}$ is $\mathcal{F}_t$-measurable for $t \geq 0$.
    \end{enumerate}
    We define $(\overline{\mathcal{F}}_t, t \geq 0)$ and $(\overline{\mathcal{F}}^y, y \geq 0)$ to be the $\mathbf{P}$-completions of $(\mathcal{F}_t, t \geq 0)$ and $(\mathcal{F}^y, y \geq 0)$ respectively on $\mathcal{N}^{\text{sp}}_{\text{fin}}$, where $(\Omega, \mathcal{A}, \mathbf{P})$ is the probability space on which the PRM $\mathbf{N}$ is defined. 
\end{definition}
For $\beta \in \mathcal{I}_H(E)$ recall $\mathbf{N}_\beta$, the $\mathbf{N}_U$ for $U \in \beta$ and the $\mathbf{P}_\beta$ probability measure from Definition \ref{definition starting interval partition}. We define the collection of excursion intervals about a fixed level $y \geq 0$ for $N \in \mathcal{N}^{\text{sp}}$ by
\begin{equation*}
    F^{\geq y}_0(N) = \sum_{[a,b] \in V^y_0(N)} \delta\Big(\ell^y_N(a),(N|^\leftarrow_{[a,b]})^+\Big).
\end{equation*}
By a slight abuse of notation with the above definition for $\beta \in \mathcal{I}_M(E)$ we define for a measure $\mathbf{N}_\beta$ as in Definition \ref{definition starting interval partition}  
\begin{equation}\label{equation PRM formulation above y}
    F^{\geq 0}_0(\mathbf{N}_\beta) \coloneqq  \sum_{U \in \beta} \delta(\mathcal{D}^M_\beta(U), \mathbf{N}_U),
\end{equation}
and we set $F^{\geq 0}_0(\mathbf{N}_\beta)$ to be the zero measure when $\beta \in \mathcal{I}_H(E) \backslash \mathcal{I}_M(E)$. We state results from \cite[Section 5]{construction} that generalise to our setting.
\begin{proposition}\label{proposition clade poisson random measure}
Consider a diffusion that satisfies $A1$ and $A2$ of Assumption \ref{assumption a} with speed measure and scale function that satisfies \eqref{equation starting IP condition}, and for which the Pitman-Yor excursion measure $\nu$ satisfies \eqref{equation levy measure unbounded variation}. Let $(\mathbf{N},\mathbf{X})$ be a scaffolding-and-spindles pair where $\mathbf{N}$ has intensity $\text{Leb} \otimes \nu$ stopped at $T$, an a.s. finite time $T$ with the following properties:
\begin{enumerate}
    \item $S^0 \coloneqq  \ell^0(T)$ is measurable in $\overline{\mathcal{F}}^0$;
    \item $\mathbf{X} < 0$ on the time interval $(\tau^0(S^0-), T)$.
\end{enumerate}
Then for each $y \geq 0$, the measure $F^{\geq y}_0(\mathbf{N})$ is conditionally independent of $\overline{\mathcal{F}}^y$ given $\beta^y\coloneqq \text{skewer}(\mathbf{N}, \mathbf{X}, y)$, with the regular conditional distribution $\mathbf{P}_{\beta^y}(F^{\geq 0}_0 \in \cdot)$.
\begin{proof}
    The proof of \cite[Proposition 5.6]{construction} applies in our setting as we have the following:
    \begin{enumerate}
        \item The scaffolding $\mathbf{X}$ is a L\'evy process of unbounded variation and therefore we can decompose $\mathbf{N}$ into biclades at every fixed level $y$ of $\mathbf{X}$: $(\mathbf{N}|^\leftarrow_I, I \in V^y_0(\mathbf{N}))$.
        \item We have a mid-spindle Markov property (Lemma \ref{lemma mid-spindle}).
        \item We have a disintegration for the bi-clade excursion measure $\nu_{\text{cld}}$ with respect to the crossing width $m^0$ (Proposition \ref{proposition disintegration}).
        \item Local times of the excursions at a fixed positive level $y>0$ are equal to the $M$-diversities of the interval partitions a.s. (Lemma \ref{lemma diversity is local time}). These are measurable functions of the skewer function at that level. 
    \end{enumerate}
\end{proof}
\end{proposition}
This result is sufficient to prove the simple Markov property with respect to the filtration generated by the skewer process. For a probability measure $\mu$ on the space of interval partitions $\mathcal{I}_H(E)$, we define the probability measure $\mathbf{P}_\mu$ on the space $\mathcal{N}^{\text{sp}}_{\text{fin}}$ by $\mathbf{P}_\mu(\cdot)\coloneqq \int \mathbf{P}_\beta(\cdot) \mu(d\beta)$.
\begin{proposition}[Simple Markov property of the skewer process under $\mathbf{P}_\mu$]\label{proposition simple markov skewer}
Consider a diffusion that satisfies $A1$ and $A2$ of Assumption \ref{assumption a} and either with speed measure and scale function that satisfies \eqref{equation starting IP condition} and for which the Pitman-Yor excursion measure $\nu$ satisfies \eqref{equation levy measure unbounded variation}, or with Pitman-Yor excursion measure that satisfies \eqref{equation levy measure bounded variation}. Let $\mu$ be a probability measure on $\mathcal{I}_H(E)$. Let $z>0$ and $0 \leq y_1< \dots < y_n$. Let $\eta: \mathcal{N}^{\text{sp}}_{\text{fin}} \rightarrow [0,\infty)$ be $\overline{\mathcal{F}}^z$ measurable, and let $f:(\mathcal{I}_H(E))^n \rightarrow [0,\infty)$ be measurable. Then
\begin{equation*}
    \mathbf{P}_\mu[\eta f(\text{skewer}(z+y_j, \cdot, \xi(\cdot)), j \in [n])] = \int \eta(N) \mathbf{P}_{\text{skewer}(z,N)} [f(\text{skewer}(y_j, \cdot, \xi(\cdot)), j \in [n])]\mathbf{P}_\mu(dN). 
\end{equation*}
   \begin{proof}   
        We split the proof into two cases. If the Pitman-Yor excursion measure satisfies \eqref{equation levy measure bounded variation}, the SPLP $\mathbf{X}$ is of bounded variation a.s. and therefore there are only finitely many crossings at level $y$ a.s. Therefore the simple Markov property follows in this case as the concatenation of the simple Markov properties for each incomplete spindle and its SPLP extension added to it. 
        \par
        For the other case, the proof follows from the previous results and it uses the same argument as given in \cite[Corollary 6.7]{construction} with the following conditions:
        \begin{enumerate}
            \item For an interval partition $\beta \in \mathcal{I}_H(E)$, and fixed $y>0$, the local times of the excursions of $\xi(\mathbf{N}_\beta)$ about level $y$ are distinct a.s. (Lemma \ref{lemma injective local times}).
            \item We have conditional independence of $F^{\geq y}_0(\mathbf{N})$ and $\overline{\mathcal{F}}^y$ given $\text{skewer}(\mathbf{N}, \xi(\mathbf{N}), y)$, for certain stopping times $T$ (Proposition \ref{proposition clade poisson random measure}). In particular, for an interval partition $\beta \in \mathcal{I}_H(E)$, and for $U \in \beta$, it is shown in \cite[Proposition 5.10]{construction} that $F_0^{\geq y}(\mathbf{N}_U)$ is conditionally independent of $\overline{\mathcal{F}}^y$ given $\text{skewer}(y,\mathbf{N}_U,\xi_{\mathbf{N}_U})$. As only finitely many of the $\mathbf{N}_U$, $U \in \beta$, will contribute to the skewer process at levels above $y$, the result follows. 
        \end{enumerate}
    \end{proof}
\end{proposition}

%% file: continuity.tex
\section{Continuity and Markov property -- Assumptions \ref{assumption b} and \ref{assumption c}}\label{section continuity}
In this section we prove Theorems \ref{theorem diffusion} and \ref{theorem interval partitions regular}. The majority of the section is dedicated to the proof of Theorem \ref{theorem diffusion} which is completed in Section \ref{subsection strong markov}. The proof of Theorem \ref{theorem interval partitions regular} is simpler and is completed at the end of Section \ref{subsection continuity in aggregate mass}.
\subsection{Scale function and speed measure under Assumptions \ref{assumption B'} and \ref{assumption c'}}\label{section model assumptions}
We now define a modification of Assumption \ref{assumption b} which we will use to prove intermediate results in this paper.
\begin{customass}{B'}\label{assumption B'}
A diffusion $\mathbf{Y}$ with state space $E \subseteq [0,\infty)$ of the form $E=[0,c]$ or $E=[0,c)$ for some $c \leq \infty$ satisfies \emph{Assumption \ref{assumption B'}} if it satisfies $A1$, $A2$, $A3$, $B1$, $B2$, $B3$.
\end{customass}
Note that this is Assumption \ref{assumption b} to the cases where the transformation $g$ is the identity map. We define Assumption \ref{assumption c'} in the same way for Assumption \ref{assumption c}.
\begin{customass}{C'}\label{assumption c'}
A diffusion $\mathbf{Y}$ with state space $E \subseteq [0,\infty)$ of the form $E=[0,c]$ or $E=[0,c)$ for some $c \leq \infty$ satisfies \emph{Assumption \ref{assumption c'}} if it satisfies $A1$, $A2$, $A3$, $C1$, $C2$. 
\end{customass}
We use the definition of the derivative $s'$ of a scale function $s$ from \eqref{equation scale and speed measure} and we take $b$ in \eqref{equation scale and speed measure} to be $\varepsilon_0$ when under Assumption \ref{assumption B'} and to be $\varepsilon_1$ when under Assumption \ref{assumption c'}. We obtain upper and lower bounds for the speed measure and scale function from the corresponding quantities in the squared Bessel case given above. We have
\begin{equation}
    \varepsilon_0^{-\alpha^+}x^{\alpha^+} = \exp\left(\int_{z=x}^{\varepsilon_0} \frac{-\alpha^+}{z}dz\right) \leq s'(x) =  \exp\left(\int_{z=x}^{\varepsilon_0} \frac{\mu(z)}{2z}dz\right) \leq \varepsilon_0^{\alpha^-} x^{\alpha^-} , x \in (0,\varepsilon_0) \text{   under \ref{assumption B'}}
\end{equation}
and
\begin{equation}\label{equation s' bounds}
    \varepsilon_1^{\beta^-} x^{-\beta^-} = \exp\left(\int_{z=x}^{\varepsilon_1} \frac{\beta^-}{z}dz\right) \leq s'(x) =  \exp\left(\int_{z=x}^{\varepsilon_1} \frac{\mu(z)}{2z}dz\right) \leq \varepsilon_1^{\beta^+} x^{-\beta^+} , x \in (0,\varepsilon_0) \text{   under \ref{assumption c'}}.
\end{equation}
For the density of the speed measure this results in 
\begin{equation}\label{equation m bounds B}
    \frac{1}{2} \varepsilon_0^{\alpha^-} x^{-1-\alpha^-} \leq m(x) \leq \frac{1}{2} \varepsilon_0^{\alpha^+} x^{-1-\alpha^+} , x \in (0,\varepsilon_0) \text{   under \ref{assumption B'}}
\end{equation}
and
\begin{equation}\label{equation m bounds C}
    \frac{1}{2} \varepsilon_1^{-\beta^+} x^{-1+\beta^+} \leq m(x) \leq \frac{1}{2} \varepsilon_1^{-\beta^-} x^{-1+\beta^-} , x \in (0,\varepsilon_0) \text{   under \ref{assumption c'}}.
\end{equation}
For $s(x)$ we have
\begin{equation}\label{equation scale function}
    \frac{\varepsilon_0^{-\alpha^+}}{1+\alpha^+}  x^{1+\alpha^+} = \varepsilon_0^{-\alpha^+} \int_{y=0}^x y^{\alpha^+} dy \leq s(x) \leq \varepsilon_0^{-\alpha^-} \int_{y=0}^x y^{\alpha^-} dy \leq \frac{\varepsilon_0^{-\alpha^-}}{1+\alpha^-}  x^{1+\alpha^-}, x \in (0,\varepsilon_0) \text{   under \ref{assumption B'}}.
\end{equation}
It immediately follows that diffusions under Assumption \ref{assumption c} have scaffolding of bounded variation because $0$ is a regular point for these diffusions. As for diffusions under Assumption \ref{assumption b}, if we further assume $A4$ then $B1$ is immediate and is not required as a condition. We include this in the result below for completeness.
\begin{proposition}\label{corollary levy measure assumption b}
    Consider a diffusion that satisfies Assumption \ref{assumption b} and let $(\mathbf{N},\mathbf{X})$ be a scaffolding-and-spindles pair for this diffusion which is stopped at some a.s. finite random time $T$. Then we have the aggregate health process $M^0_{\mathbf{N}, \mathbf{X}}(T)$ is finite a.s.
    \par
    Furthermore, if instead we have a diffusion that satisfies all of Assumption \ref{assumption a} and $B2$, $B3$ then the diffusion also satisfies $B1$.
\begin{proof} 
For the first result we note that the untransformed diffusion $\mathbf{Y}$ satisfies Assumption \ref{assumption B'} and so we can use the bound in \eqref{equation m bounds B} for its speed measure $m_{\mathbf{Y}}$. Let $q \in (\alpha^+, q_0)$ and let $C_q$ be the H\"older constant for $g$ as in from Assumption $B4$. Then we have
\begin{equation*}
    \int_{y=0}^{\varepsilon_0} g(y) m_{\mathbf{Y}}(y) dy \leq \frac{1}{2} C_q \varepsilon_0^{\alpha^+} \int_{y=0}^{\varepsilon_0} y^{-1-\alpha^++q} dy = \frac{C_q \varepsilon_0^{\alpha^+}}{2(q-\alpha^+)}<\infty.
\end{equation*}
For the second result we begin by noting that on $(0,\varepsilon_0)$
\begin{equation*}
    s'(x) = \exp\left(\int_{z=x}^{\varepsilon_0} \frac{\mu(z)}{2z} dz\right),
\end{equation*}
and so as $\mu$ is negative on $(0,\varepsilon_0)$ by $B3$ we have $s'$ is increasing on $(0,\varepsilon_0)$, and so $s(x) \leq x s'(x)$ for $x \in (0,\varepsilon_0)$. Recall from \eqref{equation scale and speed measure} that $m(x)=2/(\sigma^2(x)s'(x))$. Then, as $\sigma^2(x)=4x$ by $B2$ we have 
\begin{equation*}
    \int_{x=0}^{\varepsilon_0} m(x) \int_{z=0}^x s(z)m(z) dz dx \leq \frac{1}{2}\int_{x=0}^{\varepsilon_0} m(x) \int_{z=0}^x dz dx
\end{equation*}
which by the bounds for $m$ given in \eqref{equation m bounds B} (that follow from $B2$ and $B3$) is less than or equal to
\begin{equation*}
    \frac{\varepsilon_0^{\alpha^+}}{4}\int_{x=0}^{\varepsilon_0} x^{-\alpha^+} dx
\end{equation*}
which is finite, and so the condition in \eqref{equation x^2 bound} is satisfied. Also, by the lower bound for $m(x)$ given in \eqref{equation m bounds B} we have
\begin{equation*}
    M((0,\varepsilon_0)) \geq \frac{\varepsilon_0^{\alpha^-}}{2} \int_{x=0}^{\varepsilon_0} x^{-1-\alpha^-} dx = \infty,
\end{equation*}
and so the condition in \eqref{equation 0 is regular} is not satisfied and therefore by Theorem \ref{theorem l\'evy measure condition} we have that the Pitman-Yor excursion measure satisfies \eqref{equation levy measure unbounded variation} as required.
\end{proof}
\end{proposition}
\subsection{Continuity of local time}\label{subsection continuity local time}
In our proof Theorem \ref{theorem diffusion} we use the continuity of the local time in space and time of the scaffolding process $\mathbf{X}$. To show continuity of the local time, we first prove some technical lemmas and we use them to show a necessary and sufficient condition for the continuity of local time on Laplace exponents of SPLPs with infinite variation, from Lambert and Simatos \cite[Lemma 2.2]{lambert2015asymptotic}. We note that necessary and sufficient conditions for the space-time continuity of the local time of L\'evy processes more generally had been determined by Barlow \cite{barlow1988} but the conditions in \cite{lambert2015asymptotic} are more accessible for us to verify in our setting. To begin, we define the measure of the overshoot (or equivalently of the undershoot) $\chi$ from the L\'evy measure $\nu(\zeta \in \cdot)$ as
\begin{equation*}
    \chi(\cdot) \coloneqq \nu_{\text{cld}}(\zeta(\hat{f}_{T^+_0}) \in \cdot)= \nu_{\text{cld}}(\zeta(\check{f}_{T^+_0}) \in \cdot),
\end{equation*}
where the second equality is an immediate consequence of \eqref{equation reversal under nu clade}. We also have from Lemma \ref{lemma excursion overshoot-undershoot measure} that
\begin{equation*}
    \chi(dz) = \nu_{\text{cld}}(\zeta(\check{f}_{T^+_0}) \in dz) = \int_{y=0}^\infty \nu(\zeta \in dz+y) dy = \nu(\zeta \in (z,\infty))dz.
\end{equation*}
\begin{lemma}\label{lemma chi lower bound}
Consider a diffusion that satisfies Assumption \ref{assumption B'}. Let $\nu$ be its Pitman-Yor excursion measure and let $\chi$ be defined from $\nu$ as above. Then there exist $\varepsilon_2 >0$ and $k_5>0$ such that
\begin{equation*}
    \chi((x,\infty)) \geq k_5 x^{-\alpha^-} \text{   for all   } x \in (0,\varepsilon_2).
\end{equation*}
\begin{proof}
First note that we have
\begin{equation}\label{equation chi definition}
    \chi((x,\infty)) = \nu_{\text{cld}}(\zeta(\hat{f}_{T^+_0}) >x) = k_0 \int_{u=0}^\infty m(u) \mathbb{Q}^0_u(\zeta >x) du,
\end{equation}
where $k_0$ is the constant given in Lemma \ref{lemma excursion overshoot-undershoot measure}. We can bound the density of the speed measure $m$ below by \eqref{equation m bounds B}. Let $u \in (0,1)$, and let $B$ be a standard Brownian motion on some probability space $(\Omega, \mathcal{F}, \mathbb{P})$. Define the stochastic processes $Y$, $Y^+$, $Y^-$ by the stochastic differential equations,
\begin{equation}\label{equation coupling of diffusions}
    dY_t = u + \mu(Y_t)dt + 2\sqrt{|Y_t|}dB_t; \hspace{0.25cm} dY^+_t = u - 2\alpha^+dt + 2\sqrt{|Y^+_t|}dB_t; \hspace{0.25cm} dY^-_t = u - 2\alpha^-dt + 2\sqrt{|Y^-_t|}dB_t,
\end{equation}
which are stopped when they hit $0$. We have that
\begin{equation*}
    \mathbb{P}(Y \in \cdot) = \mathbb{Q}_u(\cdot); \hspace{0.25cm} \mathbb{P}(Y^+ \in \cdot) = {^{-2\alpha^+}}\mathbb{Q}_u(\cdot); \hspace{0.25cm} \mathbb{P}(Y^- \in \cdot) = {^{-2\alpha^-}}\mathbb{Q}_u(\cdot),
\end{equation*}
and furthermore that
\begin{equation*}
    Y^+_t \leq Y_t \leq Y^-_t \text{   for   } t \in (0,T_0(Y^+) \wedge T_{\varepsilon_0}(Y^-)) \hspace{0.25cm} \mathbb{P}\text{-a.s.},
\end{equation*}
where the $\varepsilon_0$ is as in Assumption \ref{assumption B'}. Therefore, by using these inequalities along with the marginal distributions of $Y$, $Y^+$, $Y^-$ under $\mathbb{P}$, we have that
\begin{equation*}
    \begin{split}
    \mathbb{Q}^0_u(\zeta >x) = \mathbb{P}(T_0(Y) >x) & \geq \mathbb{P}\Big(T_0(Y^+) >x, \sup_{0 \leq t \leq T_0(Y^-)}\{Y^-_t\}<\varepsilon_0 \Big) \\
    & \geq {^{-2\alpha^+}}\mathbb{Q}^0_u(\zeta >x) - {^{-2\alpha^-}}\mathbb{Q}^0_u(A\geq \varepsilon_0). \\
    \end{split}
\end{equation*}
G\"oing-Jaeschke and Yor \cite[Equation (13)]{going2003survey} state $\zeta$ under ${^{-2\alpha}}\mathbb{Q}^0_u$ has distribution $\text{InverseGamma}(1+\alpha,u/2)$, i.e.
\begin{equation}\label{equation zeta distribution}
    {^{-2\alpha}}\mathbb{Q}^0_u(\zeta \in dx) = \frac{1}{\Gamma(1+\alpha)} \Big(\frac{u}{2}\Big)^{1+\alpha} x^{-2-\alpha} \exp\Big(-\frac{u}{2x}\Big).
\end{equation}
By the two-sided exit problem and \eqref{equation scale function} we have ${^{-2\alpha^-}}\mathbb{Q}^0_u(A \geq \varepsilon_0) = (u/\varepsilon_0)^{1+\alpha^-}$. For $x<\varepsilon_0/4$, we now state an upper bound for $\chi((x,\infty))$ by using the upper bound for $\mathbb{Q}^0_u(\zeta >x)$ given above along with the upper bound for $m$ given in \eqref{equation m bounds B}. We also upper bound $\mathbb{Q}^0_u(\zeta>x)$ by $\mathbb{Q}^0_u(\zeta \in (x, \varepsilon_0))$ and obtain
\begin{equation*}
    \chi((x,\infty)) \geq \int_{u=0}^{\varepsilon_0} \frac{k_0\varepsilon_0^{\alpha^-}}{2} u^{-1-\alpha^-} \Big(\int_{y=x}^{\varepsilon_0/2} \Big(\Big(\frac{1}{\Gamma(1+\alpha^+)} \frac{u^{1+\alpha^+}}{y^{2+\alpha^+}} \exp\Big(-\frac{u}{2y}\Big) dy - \Big(\frac{u}{\varepsilon_0}\Big)^{1+\alpha^-}\Big)\wedge 0\Big) \Big)du
\end{equation*}
which is greater than or equal to
\begin{equation*}
    \int_{u=0}^{\varepsilon_0} \frac{k_0\varepsilon_0^{\alpha^-}}{2} u^{-1-\alpha^-} \int_{y=x}^{\varepsilon_0/2} \frac{1}{\Gamma(1+\alpha^+)} \frac{u^{1+\alpha^+}}{y^{2+\alpha^+}} \exp\Big(-\frac{u}{2y}\Big) dy du -\int_{u=0}^{\varepsilon_0} \frac{k_0\varepsilon_0^{\alpha^-}}{2} u^{-1-\alpha^-} \int_{y=x}^{\varepsilon_0/2}\Big(\frac{u}{\varepsilon_0}\Big)^{1+\alpha^-}dydu.
\end{equation*}
We consider the two double integrals separately. For the first double integral we have swap the integrals to find and restrict $u \in (0,\varepsilon_0)$ to $u \in (0,2y)$ as $y \in (x,\varepsilon_0/2)$. We get a lower bound of
\begin{equation*}
    \frac{k_0\varepsilon_0^{\alpha^-}}{2\Gamma(1+\alpha^+)} \int_{y=x}^{\varepsilon_0/2} y^{-2-\alpha^+} \int_{u=0}^{2y} u^{\alpha^+-\alpha^-} \exp\Big(-\frac{u}{2y}\Big) du dy
\end{equation*}
and we can bound $\exp(-u/2y)$ from below by $1/e$ to get a lower bound of
\begin{equation*}
    \frac{k_0\varepsilon_0^{\alpha^-}}{2\Gamma(1+\alpha^+)e} \int_{y=x}^{\varepsilon_0/2} y^{-2-\alpha^+} \int_{u=0}^{2y} u^{\alpha^+-\alpha^-} du dy = \frac{2^{\alpha^+-\alpha^-}k_0\varepsilon_0^{\alpha^-}}{\Gamma(1+\alpha^+)e(1-\alpha^++\alpha^-)} \int_{y=x}^{\varepsilon_0/2} y^{-1-\alpha^-} dy
\end{equation*}
which is equal to 
\begin{equation*}
    \frac{2^{\alpha^+-\alpha^-}k_0\varepsilon_0^{\alpha^-}}{\Gamma(1+\alpha^+)e(1-\alpha^++\alpha^-)\alpha^-} \Big( x^{-\alpha^-} - \Big(\frac{\varepsilon_0}{2}\Big)^{-\alpha^-}\Big).
\end{equation*}
For the second double integral we see the upper bound given by
\begin{equation*}
    \int_{u=0}^{\varepsilon_0} \frac{k_0\varepsilon_0^{\alpha^-}}{2} u^{-1-\alpha^-} \int_{y=0}^{\varepsilon_0/2}\Big(\frac{u}{\varepsilon_0}\Big)^{1+\alpha^-}dydu = \frac{k_0}{2\varepsilon_0} \int_{u=0}^{\varepsilon_0} \int_{y=0}^{\varepsilon_0/2} dy du = \frac{k_0\varepsilon_0}{4}. 
\end{equation*}
Therefore the claim in the lemma is satisfied for some $k_5$ and $\varepsilon_2$.
\end{proof}
\end{lemma}
\begin{lemma}\label{lemma chi upper bound}
In the same setting as in Lemma \ref{lemma chi lower bound}, we have that, for some positive constant $k_6$,
\begin{equation*}
    \chi((x,\infty)) \leq k_6 (x^{-\alpha^+} + 1) \text{   for all  } x \in (0,\infty),
\end{equation*}
\begin{proof}
We use the formulation of $\chi((x,\infty))$ as given in \eqref{equation chi definition}, as well as the coupled processes $Y$, $Y^+$ and $Y^-$ as in Lemma \ref{lemma chi lower bound}. We note that
\begin{equation*}
    \begin{split}
        \mathbb{Q}^0_u(\zeta >x) = \mathbb{P}(T_0(Y) >x) & = \mathbb{P}\Big(T_0(Y) >x, \sup_{0 \leq t \leq T_0(Y^-)} \{Y^-_t\} \geq \varepsilon_0\Big) + \mathbb{P}\Big(T_0(Y) >x, \sup_{0 \leq t \leq T_0(Y^-)} \{Y^-_t\} < \varepsilon_0\Big) \\
        & \leq \mathbb{P}(A(Y^-)\geq \varepsilon_0) + \mathbb{P}(\zeta(Y^-) >x) = {^{-2\alpha^-}}\mathbb{Q}_u^0(A\geq \varepsilon_0) + {^{-2\alpha^-}}\mathbb{Q}^0_u(\zeta >x), \\
    \end{split}
\end{equation*}
where the $\varepsilon_0$ is as in Assumption \ref{assumption B'}. Therefore we have that
\begin{equation*}
    \chi((x,\infty)) \leq \frac{k_0\varepsilon_0^{\alpha^+}}{2}\int_{u=0}^{\varepsilon_0} u^{-1-\alpha^+} \Big({^{-2\alpha^-}}\mathbb{Q}_u^0(A \geq \varepsilon_0) + {^{-2\alpha^-}}\mathbb{Q}^0_u(\zeta >x)\Big) du + k_0\int_{u=\varepsilon_0}^c \mathbb{Q}^0_u(\zeta > x) m(u) du
\end{equation*}
and we bound this above by
\begin{equation*}
    \frac{k_0\varepsilon_0^{\alpha^+}}{2} \int_{u=0}^{\varepsilon_0} u^{-1-\alpha^+} \Big(u^{1+\alpha^-} + \int_{y=x}^\infty \frac{1}{\Gamma(1+\alpha^-)} \frac{u^{1+\alpha^-}}{y^{2+\alpha^-}} \exp\Big(-\frac{u}{2y}\Big) dy \Big) du + \frac{k_0}{2} M(E \backslash [0,\varepsilon_0]).
\end{equation*}
We deal with the two parts of the integral separately. We have
\begin{equation*}
    \frac{k_0\varepsilon_0^{\alpha^+}}{2}\int_{u=0}^{\varepsilon_0} u^{-1-\alpha^+} \int_{y=x}^\infty \frac{1}{\Gamma(1+\alpha^-)} \frac{u^{1+\alpha^-}}{y^{2+\alpha^-}} \exp\Big(-\frac{u}{2y}\Big) dy du \leq \frac{k_0\varepsilon_0^{\alpha^+}}{2\Gamma(1+\alpha^-)}\int_{y=x}^\infty y^{-2-\alpha^-} \int_{u=0}^{\varepsilon_0} u^{\alpha^- - \alpha^+} du dy
\end{equation*}
which equals
\begin{equation*}
    \frac{k_0\varepsilon_0^{1+\alpha^-}}{2\Gamma(1+\alpha^-)(1+\alpha^--\alpha^+)} \int_{y=x}^\infty y^{-1-\alpha^+} dy = \frac{k_0\varepsilon_0^{1+\alpha^-}}{2\Gamma(1+\alpha^-)(1+\alpha^--\alpha^+)\alpha^+} x^{-\alpha^+}. 
\end{equation*}
For the second part of the integral we have that
\begin{equation*}
    \frac{k_0\varepsilon_0^{\alpha^+}}{2} \int_{u=0}^{\varepsilon_0} u^{\alpha^--\alpha^+} du = \frac{k_0\varepsilon_0^{1+\alpha^-}}{2(1+\alpha^--\alpha^+)},
\end{equation*}
and so
\begin{equation*}
    \chi((x,\infty)) \leq \frac{k_0\varepsilon_0^{1+\alpha^-}}{2\Gamma(1+\alpha^-)(1+\alpha^--\alpha^+)\alpha^+} x^{-\alpha^+} + \frac{k_0\varepsilon_0^{1+\alpha^-}}{2(1+\alpha^--\alpha^+)} + \frac{k_0}{2} M(E \backslash [0,\varepsilon_0])
\end{equation*}
and the result follows.
\end{proof}
\end{lemma}
\begin{lemma}\label{lemma integration by parts}
Consider a diffusion that satisfies Assumption \ref{assumption B'}, and let $\psi$ be the Laplace exponent for the scaffolding process $\mathbf{X}$, see \eqref{equation laplace exponent}. Then we have for $\lambda \in (0,\infty)$ that
\begin{equation}\label{equation improper integral representation}
    \psi(\lambda) = \int_{x=0}^\infty \lambda^2 e^{-\lambda x} \chi((x,\infty))dx.
\end{equation}
\begin{proof}
    This follows by integration-by-parts applied twice to $\psi(\lambda) = \int_{x=0}^\infty (e^{-\lambda x} - 1 + \lambda x) \Lambda(dx)$ and we require the upper bound for $\chi((x,\infty))$ in Lemma \ref{lemma chi upper bound} to do this.
\end{proof}
\end{lemma}
We now use a necessary and sufficient condition from Lambert \cite{lambert2015asymptotic} to determine that the local time is a.s. jointly continuous in time and space.
\begin{proposition}\label{corollary bi-continuous}
    Consider a diffusion under Assumption \ref{assumption b} and let $(\mathbf{N}, \mathbf{X})$ be a scaffolding-and-spindles construction. Then the local time $(\ell^u_t: t , u \in \mathbb{R})$ is jointly continuous a.s.
    \begin{proof}
        By Lambert \cite[Lemma 2.2]{lambert2015asymptotic}, as $\mathbf{X}$ is a SPLP that satisfies \eqref{equation levy measure unbounded variation}, a necessary and sufficient condition for the local time to be jointly continuous a.s. is 
        \begin{equation*}
            \int^{\infty-} \frac{d\lambda}{\psi(\lambda) \sqrt{\log(\lambda)}} < \infty.
        \end{equation*}
        Now, with the lower bound for $\chi((x,\infty))$ from Lemma \ref{lemma chi lower bound} we have that
        \begin{equation*}
            \psi(\lambda) \geq k_5 \int_{x=0}^\delta \lambda^2 e^{-\lambda x} x^{-\alpha^-}dx.
        \end{equation*}
        We let $\lambda>1/\delta$ and consider the Taylor expansion of $e^{-\lambda x}$ to obtain that
        \begin{equation*}
            \psi(\lambda) \geq k_5 \sum_{n=0}^\infty \frac{(-1)^n \lambda^{n+2}}{n!} \int_{x=0}^{1/\lambda} x^{n-\alpha^-} dx = k_5 \Big( \sum_{n=0}^\infty \frac{(-1)^n}{n!(n+1-\alpha^-)} \Big) \lambda^{1+\alpha^-}
        \end{equation*}
        and therefore
        \begin{equation*}
            \int_{\lambda = 1/\delta}^\infty \frac{d\lambda}{\psi(\lambda) \sqrt{\log(\lambda)}} \leq  \Big( k_5\sum_{n=0}^\infty \frac{(-1)^n}{n!(n+1-\alpha^-)} \Big)^{-1} \int_{\lambda = 1/\delta}^\infty \frac{d\lambda}{\lambda^{1+\alpha^-} \sqrt{\log(\lambda)}} < \infty.
        \end{equation*}
    \end{proof}
\end{proposition}
\subsection{Some further calculations}\label{subsection diffusion deviations}
In this section we give upper bounds of certain probabilities and expectations for diffusions with a state space $E=[0,\infty)$ and we use these results in our proof of Theorem \ref{theorem diffusion}. In Section \ref{subsection conditional moments} we prove results for a general such diffusions conditioned to not hit $0$ and we reference the literature on quasi-stationary distributions, M\'el\'eard and Villemonais \cite{meleard2012quasi}. In Section \ref{section squared bessel deviations} we consider squared Bessel cases of negative dimension.
\subsubsection{Expectation for a diffusion conditioned not to hit level $0$}\label{subsection conditional moments}
In this section we calculate conditional moments for the deviation of a diffusion with state space $E \subseteq [0,\infty)$. As before, for $s>0$, $\mathbb{Q}_s$ denotes the law of the diffusion started from $s$ and $\mathbb{Q}^0_s$ denotes the law of the diffusion started from $s$ and stopped at $0$.
\begin{lemma}\label{lemma conditional expectation lower bound}
Consider a diffusion that satisfies $A1$ of Assumption \ref{assumption a}. Then for $r$, $s>0$, $0<a<c$ and $0 < t \leq r$ we have
\begin{equation*}
    \mathbb{Q}^0_s(Y_t \in E \backslash [0,a)) \leq \mathbb{Q}^0_s(Y_t \in E \backslash [0,a) | \zeta > r).
\end{equation*}
\begin{proof}
Let $t \in (0,r)$. We first note that
\begin{equation*}
    b\coloneqq\mathbb{Q}^0_s(\zeta>r | Y_t \in [0,a)) \leq \mathbb{Q}^0_a(\zeta > r-t) \leq \mathbb{Q}^0_s(\zeta>r | Y_t \in E \backslash [0,a))=:c,
\end{equation*}
by the Markov property and almost-sure continuity of the sample paths. Let $p = \mathbb{Q}^0_s(Y_t \in E \backslash [0,a))$. We then have that
\begin{equation*}
    \mathbb{Q}^0_s(Y_t \in E \backslash [0,a))\mathbb{Q}^0_s(\zeta>r) = p (cp + b(1-p)) \leq cp = \mathbb{Q}^0_s(Y_t \in [a, \infty), \zeta>r)
\end{equation*}
which implies that
\begin{equation*}
    \mathbb{Q}^0_s(Y_t \in E \backslash [0,a)) \leq \mathbb{Q}^0_s(Y_t \in [a, \infty)| \zeta>r).
\end{equation*}
Taking $t \uparrow r$  we have by a.s. continuity of the sample paths that
\begin{equation*}
    \mathbb{Q}^0_s(Y_r \in E \backslash [0,a)) \leq \mathbb{Q}^0_s(Y_r \in [a, \infty)| \zeta>r).
\end{equation*}
\end{proof}
\end{lemma}
\begin{lemma}\label{lemma conditional expectation upper bound}
In the same setting as Lemma \ref{lemma conditional expectation lower bound}, for $r$, $s$, $u>0$ and $0<a<c$ we have
\begin{equation*}
    \mathbb{Q}^0_s(Y_r \in E \backslash [0,a)|\zeta > r) \leq \mathbb{Q}^0_s(Y_r \in E \backslash [0,a) | \zeta > r+u).
\end{equation*}
\begin{proof}
In a very similar fashion to the previous lemma we first note that
\begin{equation*}
    b_u\coloneqq \mathbb{Q}^0_s(\zeta > r+u | Y_r \in (0,a)) \leq \mathbb{Q}^0_a(\zeta>u) \leq \mathbb{Q}^0_s(\zeta > r+u | Y_r \in E \backslash [0,a)) =:c_u,
\end{equation*}
and let $p\coloneqq \mathbb{Q}^0_s(Y_r \in (a,\infty)) \leq \mathbb{Q}^0_s(\zeta>y)=:q$. Then we have that
\begin{equation*}
    \mathbb{Q}^0_s(\zeta > r+u, Y_r \in E \backslash [0,a)) = c_up.
\end{equation*}
Furthermore note that
\begin{equation*}
    \mathbb{Q}^0_s(\zeta > r+u | \zeta > r) = c_u \frac{p}{q} + b_u \frac{q-p}{q} \leq c_u, \hspace{0.25cm} \mathbb{Q}^0_s(Y_r \in E \backslash [0,a)| \zeta > r) = \frac{p}{q},
\end{equation*}
and therefore
\begin{equation*}
    \mathbb{Q}^0_s(Y_r \in E \backslash [0,a), \zeta>r+u|\zeta > r) = \frac{pc}{q} \geq \mathbb{Q}^0_s(Y_r \in E \backslash [0,a)|\zeta>r) \mathbb{Q}^0_s(\zeta > r+u|\zeta>r),
\end{equation*}
and the claim follows.
\end{proof}
\end{lemma}
This final lemma, Lemma \ref{lemma expectation increasing functional}, uses the following convergence result for the $\uparrow$-diffusion from \cite[Theorem 27]{meleard2012quasi}.
\begin{lemma}\label{lemma never hit zero}
In the same setting as the previous lemmas, for $r$, $s$, $a>0$ we have
\begin{equation*}
    \lim_{u \rightarrow \infty} \mathbb{Q}^0_s(Y_r \in  E \backslash [0,a) | \zeta > r+u) = \mathbb{Q}^\uparrow_s(Y_r \in (a,\infty) E \backslash [0,a)).
\end{equation*}
\end{lemma}
\begin{lemma}\label{lemma expectation increasing functional}
Let $F:[0,\infty) \rightarrow \mathbb{R}$ be a continuous increasing function with derivative $F'$ Lebesgue-almost everywhere. In the same setting as the previous lemmas, for $r$, $s>0$ we have
\begin{equation*}
    \mathbb{E}_{\mathbb{Q}^0_s}(F(Y_r)) \leq \mathbb{E}_{\mathbb{Q}^0_s}(F(Y_r) | \zeta>r) \leq \mathbb{E}_{\mathbb{Q}^\uparrow_s} (F(Y_r)). 
\end{equation*}
\begin{proof}
We first note that
\begin{equation*}
    \mathbb{E}_{\mathbb{Q}^0_s}(F(Y_r)) = \int_{z=0}^\infty \mathbb{Q}^0_s(Y_r \in  E \backslash [0,z)) F'(z) dz,
\end{equation*}
and similarly for $\mathbb{E}_{\mathbb{Q}^0_s}(F(Y_r) | \zeta>r)$ and $\mathbb{E}_{\mathbb{Q}^\uparrow_s}(F(Y_r))$. By Lemmas \ref{lemma conditional expectation lower bound}, \ref{lemma conditional expectation upper bound} and \ref{lemma never hit zero} we have that for any $a \in (0,\infty)$
\begin{equation*}
    \mathbb{Q}^0_s(Y_r \in E \backslash [0,a)) \leq \mathbb{Q}^0_s(Y_r \in E \backslash [0,a)|\zeta>r) \leq \lim_{u \rightarrow \infty} \mathbb{Q}^0_s(Y_r \in E \backslash [0,a) | \zeta > r+u) = \mathbb{Q}^\uparrow_s(Y_r \in E \backslash [0,a)),
\end{equation*}
and the result follows. 
\end{proof}
\end{lemma}
\subsubsection{Certain calculations concerning squared Bessel processes with negative dimension}\label{section squared bessel deviations}
In this section we work with squared Bessel processes with dimension $-2\alpha$, where $\alpha \in (0,1)$. We make a slight abuse of notation and write ${^{-2\alpha}}\mathbb{Q}^0_s$ in place of $\mathbb{E}_{{^{-2\alpha}}\mathbb{Q}^0_s}$ for readability.
\begin{lemma}\label{lemma conditional moments}
Let $r$, $q>0$, $0<\alpha,s<1$ and $n \in \mathbb{N}$. Then we have
\begin{equation*}
    {^{-2\alpha}}\mathbb{Q}^0_s(|Y_r - s|^{qn} | \zeta >r ) \leq {^{-2\alpha}}\mathbb{Q}^0_s(|Y_r - s|^{qn}) + {^{4+2\alpha}}\mathbb{Q}^0_s (|Y_r - s|^{qn}).
\end{equation*}
\begin{proof}
We note that the functions
\begin{equation*}
    |z - s|^{qn} \mathbbm{1}_{z \geq s}; \hspace{0.5cm} -|z - s|^{qn} \mathbbm{1}_{z \leq s},
\end{equation*}
are increasing and Lebesgue almost everywhere differentiable. Therefore, by Lemma \ref{lemma expectation increasing functional} we have that
\begin{equation*}
    \begin{split}
        {^{-2\alpha}}\mathbb{Q}^0_s(|Y_r - s|^{qn}| \zeta >r ) & \leq {^{-2\alpha}}\mathbb{Q}^\uparrow_s(|Y_r - s|^{qn} \mathbbm{1}_{z \geq s} | \zeta >r ) - {^{-2\alpha}}\mathbb{Q}^0_s(-|Y_r - s|^{qn} \mathbbm{1}_{z \leq s} | \zeta >r ) \\
        & \leq {^{4+2\alpha}}\mathbb{Q}^0_s(|Y_r - s|^{qn}) + {^{-2\alpha}}\mathbb{Q}^0_s(|Y_r - s|^{qn}). \\
    \end{split}
\end{equation*}
\end{proof}
\end{lemma}
We establish upper and lower bounds for the probability that the lifetime $\zeta$ is greater than $r$ under ${^{-2\alpha}}\mathbb{Q}^0_s$.
\begin{lemma}
For $r$, $s>0$, $0<\alpha<1$, we have
\begin{equation}\label{equation lifetime upper bound}
    {^{-2\alpha}}\mathbb{Q}^0_s(\zeta \leq r) \leq \frac{\Gamma(1+\alpha+v)2^v}{\Gamma(1+\alpha)} \frac{r^v}{s^v}; \text{ for } v>0 \hspace{0.5cm} {^{-2\alpha}}\mathbb{Q}^0_s(\zeta>r) \leq \frac{\Gamma(1+\alpha - v)}{\Gamma(1+\alpha)2^v} \frac{s^v}{r^v} \text{ for } 0<v<1+\alpha.
\end{equation}
\begin{proof}
For the first inequality let $v>0$ and we apply the Markov inequality and \eqref{equation zeta distribution} to see that
\begin{equation*}
    {^{-2\alpha}}\mathbb{Q}^0_s(\zeta \leq r) = {^{-2\alpha}}\mathbb{Q}^0_s\left(\frac{1}{\zeta^v} \geq \frac{1}{r^v}\right) \leq {^{-2\alpha}}\mathbb{Q}^0_s\left(\frac{1}{\zeta^v}\right)r^v = r^v \int_{x=0}^\infty \frac{s^{1+\alpha}}{2^{1+\alpha}\Gamma(1+\alpha)} \frac{1}{x^{2+\alpha +v}} \exp\left(-\frac{s}{2x}\right) dx
\end{equation*}
and we integrate out the probability density function of an $\text{InverseGamma}(1+\alpha+v, s/2)$ to obtain
\begin{equation*}
    \frac{\Gamma(1+\alpha+v)2^v}{\Gamma(1+\alpha)} \frac{r^v}{s^v} \int_{x=0}^\infty \frac{s^{1+\alpha+v}}{2^{1+\alpha+v}\Gamma(1+\alpha + v)} \frac{1}{x^{2+\alpha +v}} \exp\left(-\frac{s}{2x}\right) dx = \frac{\Gamma(1+\alpha+v)2^v}{\Gamma(1+\alpha)} \frac{r^v}{s^v}.
\end{equation*}
For the second inequality let $0<v<1+\alpha$, we use the same ideas:
\begin{equation*}
    {^{-2\alpha}}\mathbb{Q}^0_s(\zeta>r) \leq \frac{1}{r^v} {^{-2\alpha}}\mathbb{Q}^0_s(\zeta^v) = \frac{1}{r^v} \int_{x=0}^\infty \frac{s^{1+\alpha}}{2^{1+\alpha}\Gamma(1+\alpha)} \frac{1}{x^{2+\alpha - v}} \exp\left(-\frac{s}{2x}\right) dx = \frac{\Gamma(1+\alpha - v)}{\Gamma(1+\alpha)2^v} \frac{s^v}{r^v},
\end{equation*}
where we have integrated out the probability distribution function of an $\text{InverseGamma}(1+\alpha-v, s/2)$.
\end{proof}
\end{lemma}
We also establish a further lemma with some useful probability bounds for squared Bessel processes.
\begin{lemma}\label{lemma bessel bounds}
Let $\gamma \in (-\infty,0)\cup[2,\infty)$, and $d>0$. Then there exists a constant $C=C(\gamma,d)$ such that for any $r$, $s>0$, we have
\begin{equation*}
    {^\gamma}\mathbb{Q}^0_s(|Y_r - s|^{d}) \leq C r^{d/2} (r^{d/2} \vee s^{d/2}).
\end{equation*}
\begin{proof}
We use \cite[Equation (15)]{forman2018uniform} to note that for $d \geq 2$ and $\gamma \in \mathbb{R}$ that we have
\begin{equation}\label{equation squared bessel bound}
    \begin{split}
        \left({^\gamma}\mathbb{Q}_s(|Y_r - s|^{d})\right)^{(1/d)} & \leq |\gamma| r + 2 \sqrt{(d-1)(|s|r + r^2(|\gamma|+(d-2)))} \\
        & \leq \Big(|\gamma| + 2\sqrt{(d-1)(1+(|\gamma| +(d-2)))}\Big) r^{1/2} \Big(r^{1/2}\vee s^{1/2}\Big)
    \end{split}
\end{equation}
Squared Bessel processes with negative dimension never return to the positive half-line after the first hitting time of $0$ a.s. Also, squared Bessel processes with dimension greater than or equal to $2$ that start positive never hit $0$ a.s. Therefore, for $\gamma \in (-\infty,0)\cup [2,\infty)$ we have that
\begin{equation*}
    \left({^\gamma}\mathbb{Q}^0_s(|Y_r - s|^{d})\right)^{(1/d)} \leq  \Big(|\gamma| + 2\sqrt{(d-1)(1+(|\gamma| +(d-2)))}\Big) r^{1/2} \Big(r^{1/2} \vee s^{1/2}\Big) \text{   for   } d \geq 2, \gamma \in (-\infty,0) \cup [2,\infty)
\end{equation*}
by \eqref{equation squared bessel bound}. Now let $d<2$. Then $4/d>2$ and so we can apply Jensen's inequality to obtain that
\begin{equation*}
    {^\gamma}\mathbb{Q}^0_s(|Y_r - s|^{d}) \leq ({^\gamma}\mathbb{Q}^0_s(|Y_r - s|^{4}))^{d/4} \leq C(\gamma,4)^{d/4} r^{d/2} \Big(r^{d/2} \vee s^{d/2}\Big) \text{   for   } 0 < d < 2, \gamma \in (-\infty,0) \cup [2,\infty).
\end{equation*}
\end{proof}
\end{lemma}
\begin{lemma}\label{lemma bound for bessel moments}
Let $n \in \mathbb{N}$, $0<\alpha,\gamma<1$, $0<\kappa<1-\alpha$, and $\alpha+\kappa<q<1$. Then there exists a constant $C=C(n,q,\alpha,\gamma,\kappa)>0$ such that for $0<r\leq 1$ we have
\begin{equation*}
    \int_{s=0}^1 {^{-2\gamma}}\mathbb{Q}^0_s(|Y_r - s|^{qn}) s^{-1-\alpha} ds \leq C r^{n((q-\alpha-\kappa) \wedge q/2)}.
\end{equation*}
\begin{proof}
We can apply Lemma \ref{lemma bessel bounds} as
\begin{equation*}
    {^{-2\gamma}}\mathbb{Q}^0_s(|Y_r^q - s^q|^{n}) \leq {^{-2\gamma}}\mathbb{Q}^0_s(|Y_r - s|^{qn}) \text{   for   } 0< q< 1, n\geq 1.
\end{equation*}
We split the integral. We first consider the expectation on the event $\{\zeta \leq r\}$ and use the upper bound for ${^{-2\gamma}}\mathbb{Q}^0_s(\zeta \leq r)$ in \eqref{equation lifetime upper bound} with $v=qn-\alpha-\kappa$:
\begin{equation*}
    \begin{split}
        \int_{s=0}^1 {^{-2\gamma}}\mathbb{Q}^0_s(|Y_r - s|^{qn} : \zeta \leq r) s^{-1-\alpha} ds & \leq \frac{\Gamma(1+\gamma+qn-\alpha-\kappa)2^{qn-\alpha-\kappa}}{\Gamma(1+\alpha)} \int_{s=0}^1 s^{-1-\alpha} s^{qn} \frac{r^{qn-\alpha-\kappa}}{s^{qn-\alpha-\kappa}} ds \\
        & = \frac{\Gamma(1+\gamma+qn-\alpha-\kappa)2^{qn-\alpha-\kappa}}{\Gamma(1+\alpha)\kappa} r^{qn - \alpha - \kappa}. \\
    \end{split}
\end{equation*}
We now consider the integral on the event $\{\zeta > r\}$ on the range $0<s<r$. By Lemmas \ref{lemma conditional moments} and \ref{lemma bessel bounds} we have that
\begin{equation*}
    {^{-2\gamma}}\mathbb{Q}^0_s(|Y_r - s|^{qn} | \zeta > r ) \leq (C(-2\alpha, qn) + C(4+2\alpha, qn)) r^{qn/2} (s^{qn/2} \vee r^{qn/2}),
\end{equation*}
where $C(-2\alpha, qn)$ and $C(4+2\alpha, qn)$ are from Lemma \ref{lemma bessel bounds}. Let $K=(C(-2\alpha, qn) + C(4+2\alpha, qn))$. We use \eqref{equation lifetime upper bound} with $v=\alpha+\kappa$ to bound above ${^{-2\gamma}}\mathbb{Q}^0_s(\zeta>r)$:
\begin{equation*}
    \begin{split}
        \int_{s=0}^{r \wedge 1} {^{-2\gamma}}\mathbb{Q}^0_s[|Y_r - s|^{qn}:\zeta > r] s^{-1-\alpha} ds
        & \leq \frac{K\Gamma(1+\gamma-\alpha-\kappa)}{\Gamma(1+\gamma)2^{\alpha+\kappa}} \int_{s=0}^{r \wedge 1} \frac{s^{\alpha+\kappa}}{r^{\alpha+\kappa}} s^{-1-\alpha}r^{qn} ds \\
        & \leq \frac{K\Gamma(1+\gamma-\alpha-\kappa)}{\Gamma(1+\gamma)2^{\alpha+\kappa}} r^{qn - \alpha - \kappa} \int_{s=0}^1 s^{-1+\kappa} ds \\
        & = \frac{K\Gamma(1+\gamma-\alpha-\kappa)}{\Gamma(1+\gamma)2^{\alpha+\kappa}\kappa} r^{qn - \alpha}. \\
    \end{split}
\end{equation*}
We note that
\begin{equation*}
    r^{qn-\alpha - \kappa} \leq r^{n((q-\alpha-\kappa) \wedge q/2)} \text{   for   } 0<r<1.
\end{equation*}
Finally, we consider the integral on the event $\{\zeta > r\}$ over the range $r \wedge 1 \leq s \leq 1$. We split this calculation into two parts to deal with the cases $q/2 \geq q-\alpha-\kappa$ and $q/2 < q-\alpha-\kappa$ separately. In the case where $q/2 \geq q-\alpha-\kappa$ we use the bound in \eqref{equation lifetime upper bound} for ${^{-2\gamma}}\mathbb{Q}^0_s(\zeta>r)$ with $0 < v = \alpha + \kappa - q/2 < 1 + \alpha$:
\begin{equation*}
    \begin{split}
        \int_{s=r \wedge 1}^1 {^{-2\gamma}}\mathbb{Q}^0_s(|Y_r - s|^{qn} : \zeta > r) s^{-1-\alpha} ds & \leq \frac{K\Gamma(1+\gamma + q/2 - \alpha-\kappa)}{\Gamma(1+\alpha)2^{\alpha+\kappa-q/2}} \int_{s=r \wedge 1}^1 \frac{s^{\alpha+\kappa-q/2}}{r^{\alpha+\kappa-q/2}} s^{-1-\alpha} r^{qn/2} s^{qn/2} ds \\
        & \leq \frac{K\Gamma(1+\gamma + q/2 - \alpha-\kappa)}{\Gamma(1+\alpha)2^{\alpha+\kappa-q/2}} r^{qn/2 + q/2 - \alpha - \kappa} \int_{s=0}^1 s^{-1+\kappa} ds \\
        & \leq \frac{K\Gamma(1+\gamma + q/2 - \alpha-\kappa)}{\Gamma(1+\alpha)2^{\alpha+\kappa-q/2}\kappa} r^{n(q-\alpha-\kappa)} \text{   for   } q/2 \geq q-\alpha-\kappa. \\
    \end{split}
\end{equation*}
For the case where $q/2 < q-\alpha-\kappa$ we do not use \eqref{equation lifetime upper bound} to bound ${^{-2\gamma}}\mathbb{Q}^0_s(\zeta>r)$ (we instead use the trivial bound of $1$):
\begin{equation*}
    \begin{split}
        \int_{s=r \wedge 1}^1 {^{-2\gamma}}\mathbb{Q}^0_s(|Y_r - s|^{qn} : \zeta > r)s^{-1-\alpha} ds
        & \leq K \int_{s=r \wedge 1}^1 s^{-1-\alpha} r^{qn/2} s^{qn/2} ds \\
        & = K \int_{s=r \wedge 1}^1 s^{-1-\alpha + qn/2} ds r^{qn/2} \\
        & \leq K r^{qn/2} \int_{s=0}^1 s^{-1+\kappa} ds = \frac{K}{\kappa} r^{qn/2} \text{   for   } q/2 < q - \alpha - \kappa. \\
    \end{split}
\end{equation*}
\end{proof}
\end{lemma}
\subsection{Continuity of the skewer process}\label{subsection continuity in aggregate mass}
In this section we prove the continuity results for Theorems \ref{theorem diffusion} and \ref{theorem interval partitions regular}. We begin with some definitions. Consider a diffusion under Assumption \ref{assumption B'} with Pitman-Yor excursion measure $\nu$, and a transformation $g:E \rightarrow g(E)$ (where $E=[0,c)$ or $E=[0,c]$ for some $c \leq \infty$) that satisfies the conditions stated in Assumption \ref{assumption b}. Without loss of generality we can take the $\varepsilon_0$ in Assumption \ref{assumption b} to be less than $1$ and we do this here. As before, let $\mathbf{N}$ be a stopped PRM with intensity $\text{Leb} \otimes \nu$, where the stopping is at a positive random variable $T$. Finally, for $M$, $\varepsilon_0>0$ we define $\mathbf{N}_{M,\varepsilon_0}$ by
\begin{equation}\label{equation n with m and epsilon}
    \mathbf{N}_{M,\varepsilon_0}\coloneqq \mathbbm{1}_{\sup_{u \in \mathbb{R}}\ell^u(T)<M} \sum_{\substack{(t,f) \text{   atom of   } \mathbf{N}: \\ A(f) < \varepsilon_0}} \delta_{(t,f)}.
\end{equation}
\begin{lemma}\label{lemma x up to y}
Let $x$, $y \in \mathbb{R}$, $z \in \{x,y\}$, and let $k \in \mathbb{N}$. Then we have that
\begin{equation}\label{equation deviations bound}
    \begin{split}
        & \mathbf{E}\Bigg(\Bigg( \sum_{\substack{(t,f) \text{   atom of   } \mathbf{N}_{M,\varepsilon_0}: \\ \mathbf{X}_{t-}<z<\mathbf{X}_t, \ell^z(t)>0}} 
        \Big|g\big(f(y-\mathbf{X}_{t-})\big) - g\big(f(x-\mathbf{X}_{t-})\big)\Big|\Bigg)^k\Bigg) \\
        & \leq C^k_q\sum_{\substack{\left(n_1, \dots, n_a\right)\\  \text{composition}\\\text{of   }k}} (k_0M)^a \prod_{i=1}^a \int_{s=0}^{\varepsilon_0} m(s) \Big(\hspace{0.01cm} {^{-2\alpha^+}}\mathbb{Q}^0_s(|Y_{|x-y|}-s|^{qn_i}) + {^{-2\alpha^-}}\mathbb{Q}^0_s(|Y_{|x-y|}-s|^{qn_i})\Big) ds,
    \end{split}
\end{equation}
where $C_q$ is the H\"older constant for $g$ as in Assumption \ref{assumption b} and $k_0$ is as in Proposition \ref{proposition crossing width density}.
\begin{proof}
The bound is trivial for $x = y$ so assume $x \neq y$. First recall from Assumption \ref{assumption b} that there exists $C_q>0$ such that $|g(a)-g(b)| \leq C_q|a-b|^q$ for $0<a$, $b<\varepsilon_0$. We can bound the term on the left-hand-side of \eqref{equation deviations bound} by
\begin{equation}\label{equation bound by tau M}
    C_q^k \mathbf{E}\Bigg(\Bigg(\sum_{\substack{(t,f) \text{   atom of   } \mathbf{N}: A(f)<\varepsilon_0 \\ 0<\ell^z(t)<M, \mathbf{X}_{t-}<z<\mathbf{X}_t}} 
    \Big|f(y-\mathbf{X}_{t-}) - f(x-\mathbf{X}_{t-})\Big|^q\Bigg)^k\Bigg).
\end{equation}
As $T^z$ is a stopping time for $\mathbf{X}$, we have by Lemma \ref{lemma joint probability undershoot overshoot} that the process $\xi(\mathbf{N}|^\leftarrow_{[T^z,\infty)})$ can be split into excursions about level $0$, and all these excursions a.s. include a single jump across level $0$. We have a.s. a one-to-one correspondence between crossing across a level and excursions about the level. Therefore the sum in \eqref{equation bound by tau M} less than or equal to a sum indexed by the Poisson point process of crossing widths at level $z$ and is summing functionals of the clades or anti-clades associated with each of these crossing widths. Specifically, let $\mathbf{N}'$ be the unstopped Poisson random measure, then
\begin{equation*}
    \mathbf{J}\coloneqq\sum_{\substack{t:0<\ell^z(t)<M, \\ A(f)<\varepsilon_0 }} \delta_{(\ell^z(t), \mathbf{N}'|^\leftarrow_{(\ell^z(t)-,\ell^z(t))})}
\end{equation*}
is distributed as a Poisson random measure with intensity 
\begin{equation*}
\text{Leb}|_{[0,M]} \otimes k_0 \int_{s=0}^{\varepsilon_0} m(s)\nu_{\text{cld}}(\mathbbm{1}_{A(f_{T^+_0})<\varepsilon_0}N \in \cdot|f_{T^+_0}(-\xi_N(T^+_0-)) = s)ds,
\end{equation*}
where $k_0$ is as in Proposition \ref{proposition crossing width density} and $\varepsilon_0$ is as in Assumption \ref{assumption b}. Let $w \in \{x,y\}\backslash \{z\}$. Then, we define
\begin{equation*}
    F(N) \coloneqq |f_{T^+_0}(w-z-\xi_N(T^+_0-)) - f_{T^+_0}(-\xi_N(T^+_0-))|^q \text{   for   } N \in \mathcal{N}^{\text{sp}}_{\pm \text{cld}}
\end{equation*}
and therefore the expectation in \eqref{equation bound by tau M} is bounded above by
\begin{equation*}
    \mathbf{E}\Bigg(\Bigg(\sum_{\substack{(s,N) \text{   atom}\\\text{  of   } \mathbf{J}}} F(N)\Bigg)^k\Bigg) = \sum_{\substack{\left(n_1, \dots, n_a\right)\\  \text{composition}\\\text{of   }k}} \mathbf{E}\bigg(\sum_{\substack{(s_1,N_1), \dots, (s_a,N_a) \\ \text{distinct atoms of } \mathbf{J}}} \prod_{i=1}^a F(N_i)^{n_i}\bigg).
\end{equation*}
We apply Campbell's formula for Poisson random measures, see Kingman \cite[Equation (3.14)]{kingman1992poisson}, to obtain, for a composition $(n_1,\dots, n_a)$ of $k$, that the expectation in the right-hand-side of the above equation equals
\begin{equation*}
    (k_0M)^a\prod_{i=1}^a \int_{s=0}^{\varepsilon_0} m(s) \nu_{\text{cld}}\Big(|f_{T^+_0}(w-z - \xi_N(T^+_0-)) - f_{T^+_0}(-\xi_N(T^+_0-))|^{qn_i}\mathbbm{1}_{A(f_{T^+_0})<\varepsilon_0}\Big|f_{T^+_0}(-\xi_N(T^+_0-)) = s\Big) ds.
\end{equation*}
Now, let $\Tilde{\mathbf{Y}}$ be a diffusion with infinitesimal parameters
\begin{equation}\label{equation y tilde}
    \sigma^2_{\Tilde{\mathbf{Y}}}(x) = 4x; \hspace{0.5cm} \mu_{\Tilde{\mathbf{Y}}}(x) = \left\{
        \begin{matrix}
            \hspace{0.25cm} \mu(x) \text{   for   } x \in (0,\varepsilon_0]; \\
            \mu(\varepsilon_0) \text{   for   } x > \varepsilon_0, \\
        \end{matrix}
    \right.
\end{equation}
and let $\Tilde{\nu}_{\text{cld}}$ be the bi-clade excursion measure for $\Tilde{\mathbf{Y}}$. Then we have for any measurable set $B \subseteq \mathcal{C}$ of continuous paths we have for $C \coloneqq \{N \in \mathcal{N}^{\text{sp}}_{\pm \text{cld}}: f_{T^+_0} \in B\}$ that
\begin{equation*}
    \nu_{\text{cld}}\Big(\{A(f_{T^+_0})<\varepsilon_0\} \cap C\Big|f_{T^+_0}(-\xi_N(T^+_0-)) = s\Big) = \Tilde{\nu}_{\text{cld}}\Big(\{A(f_{T^+_0})<\varepsilon_0\} \cap C\Big|f_{T^+_0}(-\xi_N(T^+_0-))) = s\Big)
\end{equation*}
because the two excursion measures differ only on excursions with amplitude greater than $\varepsilon_0$. Indeed, you can construct a Poisson random measure $\mathbf{N}$ with bi-clade measure $\nu_{\text{cld}}$ by first constructing a Poisson random measure $\Tilde{\mathbf{N}}$ with bi-clade measure $\Tilde{\nu}_{\text{cld}}$, removing the finitely many bi-clades with spindles with amplitude greater than $\varepsilon_0$ (up to a finite time $T$) and adding finitely many bi-clades with spindles with amplitude greater than $\varepsilon_0$ since such bi-clades have finite intensity under the bi-clade measures.
\par
In Proposition \ref{proposition disintegration} we established that the bi-clade measure about a fixed level disintegrates with respect to the crossing width. Let $n \in \mathbb{N}$. The terms in the sum in \eqref{equation bound by tau M} are the functionals of the clade (if $z=\min\{x,y\}$) or the anti-clade (if $z>\min\{x,y\}$) of excursions about level $z$. Therefore we have an upper bound given by
\begin{equation*}
    \Tilde{\nu}_{\text{cld}}\Big(|f_{T^+_0}(w-z - \xi_N(T^+_0-)) - f_{T^+_0}(-\xi_N(T^+_0-))|^{qn}\Big|f_{T^+_0}(-\xi_N(T^+_0-))) = s\Big) = \Tilde{\mathbb{Q}}^0_s(|Y_{|x-y|}-s|^{qn}),
\end{equation*}
and as $-2\alpha^+ \leq \mu_{\Tilde{\mathbf{Y}}}(x) \leq -2\alpha^-$ we have, through a coupling argument, that
\begin{equation}\label{equation bound above by bessels}
    \Tilde{\mathbb{Q}}^0_s(|Y_{|x-y|}-s|^{qn}) \leq {^{-2\alpha^+}}\mathbb{Q}^0_s(|Y_{|x-y|}-s|^{qn}) + {^{-2\alpha^-}}\mathbb{Q}^0_s(|Y_{|x-y|}-s|^{qn}),
\end{equation}
which gives the bound stated in this lemma.
\end{proof}
\end{lemma}
\begin{figure}
    \centering
    \includegraphics[width=\textwidth]{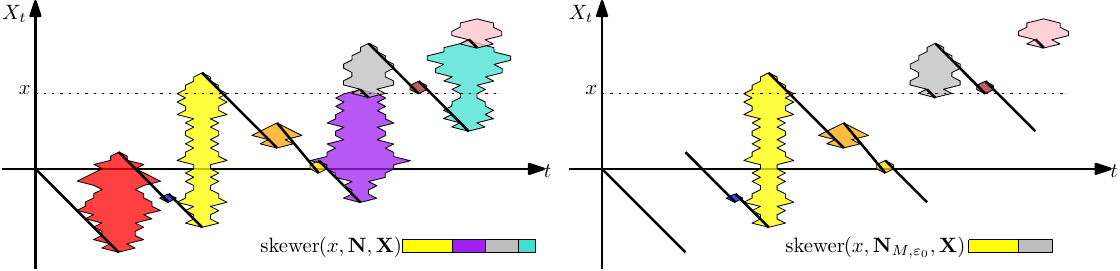}
    \caption{Left: The skewer process at level $x$ constructed from the PRM $\mathbf{N}$ and scaffolding $\mathbf{X}$. Right: The skewer process at level $x$ constructed from the reduced PRM $\mathbf{N}_{M, \varepsilon_0}$ from \eqref{equation n with m and epsilon} and scaffolding $\mathbf{X}$.}
    \label{figure reduced scaffolding and spindles}
\end{figure}
The next lemma deals with the first spindle to cross level $z>0$. In this case, we do not have a complete excursion about level $z$, and moreover the distribution of the crossing width is not clear. However, it is sufficient to know that some probability density function $p_z, z \in [0,\infty)$ gives the distribution for the crossing width where the anti-clade has infimum less than or equal to $-z$, i.e. the distribution of 
\begin{equation*}
    \nu_{\text{cld}}(f_{T^+_0}(-\xi_N(T^+_0-))) \in \cdot | \inf \xi_N \leq -z).
\end{equation*}
\begin{lemma}\label{lemma first crosses level x}
In the setting of Lemma \ref{lemma x up to y} we have that 
\begin{equation}\label{equation initial contributions}
    \begin{split}
        & \mathbf{E}\Bigg(\Bigg( \sum_{\substack{(t,f) \text{   atom of   } \mathbf{N}:\\ t<T, \ell^z(t)=0, \\ A(f) < \varepsilon_0, \mathbf{X}_{t-}<z<\mathbf{X}_t}}\Big|g\big(f(y-\mathbf{X}_{t-})\big) - g\big(f(x-\mathbf{X}_{t-})\big)\Big|\Bigg)^k\Bigg) \\
        & \leq C_q^k \int_{s=0}^{\varepsilon_0} p_z(s) \Big(\hspace{0.01cm} {^{-2\alpha^+}}\mathbb{Q}^0_s(|Y_{|x-y|}-s|^{qk}) + {^{-2\alpha^-}}\mathbb{Q}^0_s(|Y_{|x-y|}-s|^{qk})\Big) ds \\
    \end{split}
\end{equation}
where $p_z$ is the probability density function on $[0,\infty)$ that gives the distribution for 
$\nu_{\text{cld}}(m^0 \in \cdot | \zeta^+(N^-)>z)$ (i.e. the distribution for the crossing width under $\nu_{\text{cld}}$ conditional on the anti-clade having depth at least $-z$) and where $C_q$ is the H\"older constant given in Assumption \ref{assumption b}.
\begin{proof}
    We use the same method as in the previous lemma. The sum on the left-hand-side consists of the first spindle to cross level $z$ a.s. if it exists and has amplitude less than $\varepsilon_0$. We can then apply Lemma \ref{lemma mid-spindle} with stopping time $T^{\geq z}$. This tells us that the spindle above level $z$ evolves diffusively. We compare this to the diffusion $\Tilde{\mathbf{Y}}$ to upper bound this expectation in the same way as in Lemma \ref{lemma x up to y}.
\end{proof}
\end{lemma}
\begin{theorem}\label{theorem continuity}
    Consider a diffusion that satisfies Assumption \ref{assumption b}. Let $(\mathbf{N}, \mathbf{X})$ be a scaffolding-and-spindles construction for this diffusion on some probability space $(\Omega,\mathcal{A},\mathbf{P})$. Let $T$ be any $(0,\infty)$-valued $\mathbf{P}$-random variable. Then there exists a continuous IP-evolution that is a modification of $(\text{skewer}(y, \mathbf{N}, \mathbf{X}), y \in \mathbb{R})$ which takes values on $(\mathcal{I}_H(g(E)), d'_H)$.
    \begin{proof}
    First assume $g(E)$ is closed and $x$, $y \geq 0$. Let $\kappa>0$ be small enough so that $\alpha^++\kappa<q$. We first show that there exists $A \geq 0$ such that
    \begin{equation}\label{equation sum of differences x and y}
        \mathbf{E}\Bigg(\Bigg(\sum_{\substack{(t,f) \text{   atom of   } \mathbf{N}_{M,\varepsilon_0}: \\ \mathbf{X}_{t-}<x<\mathbf{X}_t \\ \text{   or   } \mathbf{X}_{t-}<y<\mathbf{X}_t}} \Big|g\big(f(y-\mathbf{X}_{t-})\big) - g\big(f(x-\mathbf{X}_{t-})\big)\Big|\Bigg)^k \Bigg) \leq A |x-y|^{k(q/2 \wedge (q-\alpha^+-\kappa))}.
    \end{equation}
    First we upper bound sum in the equation above by
    \begin{equation*}
        \Bigg( \sum_{\substack{(t,f) \text{   atom of   } \\ \mathbf{N}_{M,\varepsilon_0}: t < T, \\ \mathbf{X}_{t-}<x<\mathbf{X}_t, \\ \ell^x(t)=0}} + \sum_{\substack{(t,f) \text{   atom of   } \\ \mathbf{N}_{M,\varepsilon_0}: t < T, \\ \mathbf{X}_{t-}<y<\mathbf{X}_t, \\ \ell^y(t)=0}} + \sum_{\substack{(t,f) \text{   atom of   } \\ \mathbf{N}_{M,\varepsilon_0}: t < T, \\ \mathbf{X}_{t-}<x<\mathbf{X}_t, \\ \ell^x(t)>0}} + \sum_{\substack{(t,f) \text{   atom of   } \\ \mathbf{N}_{M,\varepsilon_0}: t < T, \\ \mathbf{X}_{t-}<y<\mathbf{X}_t, \\ \ell^y(t)>0}}\Bigg)
        \Big|g\big(f(y-\mathbf{X}_{t-})\big) - g\big(f(x-\mathbf{X}_{t-})\big)\Big|.
    \end{equation*}
    Now note that for any $a$, $b$, $c$, $d \geq 0$ we have that $(a+b+c+d)^k \leq 4^{k-1}(a^k+b^k+c^k+d^k)$ and therefore we can deal with each of the sums in the term above individually. We apply Lemmas \ref{lemma bound for bessel moments} and \ref{lemma x up to y} and the upper bound for $m$ given in \eqref{equation m bounds B} to obtain an upper bound for the third sum in the equation above:
    \begin{equation*}
        \mathbf{E}\Bigg( \Bigg(\sum_{\substack{(t,f) \text{   atom of   } \mathbf{N}_{M,\varepsilon_0}: \\ t < T, \mathbf{X}_{t-}<x<\mathbf{X}_t, \ell^x(t)>0}}  \Big|g\big(f(y-\mathbf{X}_{t-})\big) - g\big(f(x-\mathbf{X}_{t-})\big)\Big|\Bigg)^k \Bigg) \leq A' |x-y|^{k(q/2 \wedge (q-\alpha^+-\kappa))},
    \end{equation*}
    for some $A' \geq 0$ and similarly for the fourth sum. We deal with the first two sums in the same way, using Lemma \ref{lemma first crosses level x} and therefore obtain the upper bound given in \eqref{equation sum of differences x and y}. We now consider the $\mathbf{P}$-random times 
    \begin{equation}\label{equation stopping times Sn}
        S_0\coloneqq 0; \hspace{0.5cm} S_n\coloneqq \inf\{ t>S_{n-1}: (t,f) \text{   atom of   } \mathbf{N}, \hspace{0.25cm} A(f) \geq \varepsilon_0 \} \text{   for   } n \geq 1.
    \end{equation}
    As $\mathbf{N}$ is a stopped PRM stopped at random time $T$ there exists $P \in \mathbb{N}$ such that $S_{P-1} \leq T < S_P$. For $1 \leq n \leq P$ we define the random measures
    \begin{equation*}
        \mathbf{N}_{M,\varepsilon_0,n}\coloneqq \sum_{\substack{(t,f) \text{   atom of   } \mathbf{N}_{M,\varepsilon_0}: \\ t\in [0,T] \cap (S_n, S_{n+1})}} \delta_{(t,f)} \text{   for   } n \geq 0
    \end{equation*}
    so that $\mathbf{N}_{M,\varepsilon_0} = \sum_{n=0}^P \mathbf{N}_{M,\varepsilon_0,n}$. Let $n \in \mathbb{N}$. Consider the correspondence between $\text{skewer}(x,\mathbf{N}_{M,\varepsilon_0,n},\mathbf{X})$ and $\text{skewer}(y,\mathbf{N}_{M,\varepsilon_0,n},\mathbf{X})$ given by pairing all the intervals generated by spindles that cross both levels $x$ and $y$. This is a correspondence as we only have finitely many spindles cross level $x$ and level $y$ $\mathbf{P}$-a.s and has distortion
    \begin{equation*}
        \sum_{\substack{(t,f) \text{   atom  } \\ \text{   of   } \mathbf{N}_{M,\varepsilon_0,n}: \\ \mathbf{X}_{t-}<x \leq y<\mathbf{X}_t }} \Big|g\big(f(y-\mathbf{X}_{t-})\big) - g\big(f(x-\mathbf{X}_{t-})\big)\Big|+ \max\Big\{\sum_{\substack{(t,f) \text{   atom   } \\ \text{   of   } \mathbf{N}_{M,\varepsilon_0,n}: \\ x\leq \mathbf{X}_{t-}<y<\mathbf{X}_t }} g\big(f(y-\mathbf{X}_{t-})\big), \sum_{\substack{(t,f) \text{   atom   } \\ \text{   of   } \mathbf{N}_{M,\varepsilon_0}: \\ \mathbf{X}_{t-}<x<\mathbf{X}_t\leq y } } g\big(f(x-\mathbf{X}_{t-})\big)\Big\}.
    \end{equation*}
    This sum is less than or equal to the sum in \eqref{equation sum of differences x and y}. Therefore we have $\overline{\text{skewer}}(\mathbf{N}_{M,\varepsilon_0,n}, \mathbf{X})$ satisfies the Kolmogorov-\v{C}entsov continuity criterion and so there exists a continuous modification $(\gamma^y_{M,n}, y \geq 0)$. We can take $M$ to infinity as by Proposition \ref{corollary bi-continuous} we have that $\sup_{u \in \mathbb{R}}\ell^u(T)$ is finite $\mathbf{P}$-a.s. and so we have a continuous modification $(\gamma_n^y, y \geq 0)$ of $\overline{\text{skewer}}(\mathbf{N}|_{(S_n,S_{n+1})},\mathbf{X})$. 
    \par
    From here we can construct a continuous modification of $\overline{\text{skewer}}(\mathbf{N},\mathbf{X})$. For $n \in \mathbb{N}$ let $f_n$ be the $n^{\text{th}}$ spindle of $\mathbf{N}$ with amplitude greater than or equal to $\varepsilon_0$, i.e. so that $\delta_{(S_n,f_n)}$ is an atom of $\mathbf{N}$. Then we define the IP-evolution $(\gamma^y, y \geq 0)$ by
    \begin{equation}\label{equation concatenate}
        \gamma^y \coloneqq \bigstar_{n=1}^P \big(\gamma^y_n \star \{(0,f_n(y-\mathbf{X}_{S_n-}))\}\big).
    \end{equation}
    Note there are only finitely many spindles with amplitude greater than or equal to $\varepsilon_0$ on $[0,T]$ $\mathbf{P}$-a.s. Therefore we have that $(\gamma^y, y \geq 0)$ is the concatenation of finitely many $d'_H$-continuous IP-evolutions (the IP-evolutions $(\{(0,f_n(y-\mathbf{X}_{S_n}))\}, y \geq 0)$ are clearly $d'_H$-continuous) and therefore $(\gamma^y, y \geq 0)$ is $d'_H$-continuous. Finally we note that concatenating modificiations gives modifications and so $(\gamma^y, y \geq 0)$ is a modification of $\overline{\text{skewer}}(\mathbf{N}, \mathbf{X})$ on $(\mathcal{I}_H, d'_H)$. The subspace $(\mathcal{I}_H(g(E)), d'_H)$ is closed and so this modification exists on the appropriate subspace. To finish the proof assume $g(E)$ is open. This only occurs if $g(c)$ is finite and $g(E)=[0,g(c))$.
    In this case however we still obtain our result as we do not take a modification of the diffusions with amplitude greater than $\varepsilon_0$ and as each of these diffusions with amplitude greater than $\varepsilon_0$ do not hit $g(c)$ then $(\gamma^y, y \geq 0)$ does take values on the subspace $(\mathcal{I}_H(g(E)), d'_H)$. 
    \par
    Finally, assume $x$, $y \leq 0$. Then this case follows ultimately by \eqref{equation reversal under nu clade}. If $x \leq 0 \leq y$ then we have the result by the previous two cases and the triangle inequality. 
\end{proof}
\end{theorem}
We conclude this section with the analogous result for diffusions under Assumption \ref{assumption c}.
\begin{proposition}\label{proposition continuity assumption c}
    Let $\mathbf{Z}$ be a diffusion that satisfies Assumption \ref{assumption c}, and let $(\mathbf{N},\mathbf{X})$ be a scaffolding-and-spindles pair stopped at an a.s. finite stopping time. Then $(\overline{\text{skewer}}(\mathbf{N},\mathbf{X})$ is a.s. continuous and takes values on $(\mathcal{I}_H(g(E)), d'_H)$.
    \begin{proof}
    Let $\mathbf{Y}$ be an untransformed diffusion that satisfies Assumption \ref{assumption c'} and $g$ be a transformation that satisfies $C3$ of Assumption \ref{assumption c} such that $\mathbf{Z} = g(\mathbf{Y})$. 
    \par
    By the bounds in \eqref{equation m bounds C} we have that $M((0,\varepsilon_0))$ is finite. Therefore by Theorem \ref{theorem l\'evy measure condition} the Pitman-Yor excursion measure for $\mathbf{Y}$ (equivalently $\mathbf{Z}$) satisfies \eqref{equation levy measure bounded variation}. For the diffusion $\mathbf{Z}$ let $(\mathbf{N},\mathbf{X})$ be a scaffolding-and-spindles pair stopped at an a.s. finite stopping time. 
    \par
    To show continuity of the IP-evolution we first show that the sum of the amplitudes of all of the spindles is finite. From Theorem \ref{definition Williams decomposition} we know that the intensity of spindles of $\mathbf{N}$ with amplitude greater than $x$ is $1/s(x)$. Therefore, the amplitudes of the spindles is a Poisson point process with intensity $s'(x)/s(x)^2$ which, when stopped at an a.s. finite random time, is summable iff
    is finite . With the substitution $y=s(x)$ we have that it is equal to
    \begin{equation}\label{equation needs to be summable}
        \int_{x=0}^b g(x) \frac{s'(x)}{s(x)^2}dx = \int_{y=0}^{s(b)} \frac{g(s^{-1}(y))}{y^2} dy <\infty
    \end{equation}
    for some (equivalently all) $b \in (0,c)$ (we have included the integral obtained with the substitution $y=s(x)$ here as we use it below). By the bound in \eqref{equation s' bounds} of $s'$ we have
    \begin{equation*}
        \frac{\varepsilon_1^{\beta^-}}{1-\beta^-} x^{1-\beta^-} \leq s(x) \text{   for   } x \in (0,\varepsilon_1).
    \end{equation*}
    Therefore for $y=\big(\varepsilon_1^{-\beta^-}(1-\beta^-)\big)^{1/(1-\beta^-)}x^{1/(1-\beta^-)}$ we have 
    \begin{equation*}
        x = \frac{\varepsilon_1^{\beta^-}}{1-\beta^-} y^{(1-\beta^-)} \leq s(y) \text{   for   } y \in (0,\varepsilon_1).
    \end{equation*}
    This means that
    \begin{equation*}
        s^{-1}(x) \leq \big(\varepsilon_1^{-\beta^-}(1-\beta^-)\big)^{1/(1-\beta^-)}x^{1/(1-\beta^-)} \hspace{0.5cm} \text{ for } x < s(\varepsilon_1)
    \end{equation*}
    and so the second integral in \eqref{equation needs to be summable} is less than or equal to
    \begin{equation*}
        \int_{y=0}^{s(b)} \frac{g(\big(\varepsilon_1^{-\beta^-}(1-\beta^-)\big)^{1/(1-\beta^-)}y^{1/(1-\beta^-)})}{y^2} dy \hspace{0.5cm} \text{ provided that } b < \varepsilon_1.
    \end{equation*}
    To put this into the form of $C3$ requires a linear transformation $z=ay$ for $a=\big(\varepsilon_1^{-\beta^-}(1-\beta^-)\big)$. Therefore, the integral in \eqref{equation needs to be summable} is finite by $C3$ and the summability of the amplitudes follows.
    \par
    From here we can show that $\overline{\text{skewer}}(\mathbf{N},\mathbf{X})$ is continuous a.s. with respect to the metric $d'_H$. Let $I$ be an index set for the spindles of $\mathbf{N}$, $\varepsilon>0$ and $C=\sum_{i \in I}A(f_i)<\infty$ be the sum of amplitudes Let $(t_1,f_1), \dots, (t_K, f_K)$ be atoms of $\mathbf{N}$ such that $\sum_{i=1}^K A(f_i)>C-\varepsilon/2$, and let $\Tilde{\mathbf{N}}$ be the Poisson random measure $\mathbf{N}$ with the atoms $(t_1,f_1), \dots, (t_K, f_K)$ removed. The spindles $f_i$ are continuous functions with finite lifetimes so are uniformly continuous. Let $\delta>0$ such that for any $1 \leq i \leq K$ and any $s$, $t$ with $|s-t|<\delta$ we have $|f_i(s)-f_i(t)|<\varepsilon/2K$. Let $x$, $y \in \mathbb{R}$ with $|x-y|<\delta$. Then, for the interval partitions $\text{skewer}(x,\mathbf{N},\mathbf{X})$ and $\text{skewer}(y,\mathbf{N},\mathbf{X})$ we can define the correspondence of size at most $K$ that matches every interval generated by spindle $f_i$ in $\text{skewer}(x,\mathbf{N},\mathbf{X})$ to the corresponding interval in $\text{skewer}(y,\mathbf{N},\mathbf{X})$. The remaining intervals in $\text{skewer}(x,\mathbf{N},\mathbf{X})$ or $\text{skewer}(y,\mathbf{N},\mathbf{X})$ have sum less than $\varepsilon/2$ and therefore the distortion is equal to 
    \begin{equation*}
        \sum_{j=1}^K |f_i(x-\mathbf{X}_{t_i-})-f_i(y-\mathbf{X}_{t_i-})| + \max\Big\{\sum_{(t,f) \text{   atom of   } \Tilde{\mathbf{N}}} f(x-\mathbf{X}_{t-}), \sum_{(t,f) \text{   atom of   } \Tilde{\mathbf{N}}} f(y-\mathbf{X}_{t-})\Big\} \leq K \frac{\varepsilon}{2K} + \frac{\varepsilon}{2}
    \end{equation*}
    and therefore $d'_H(\text{skewer}(x,\mathbf{N},\mathbf{X}), \text{skewer}(y,\mathbf{N},\mathbf{X}))\leq \varepsilon$. As in the previous result the process takes values on the subspace $(\mathcal{I}_H(g(E)), d'_H)$ as required.
    \end{proof}
\end{proposition}
\subsection{Simple Markov property and proof of Theorem \ref{theorem interval partitions regular}}\label{section simple markov}
For a measurable map $g:E \rightarrow g(E)$ (where as before $E=[0,c)$ or $E=[0,c]$ for some $c \leq \infty$) we define $g^*:\mathcal{I}_H(E) \rightarrow g^*(\mathcal{I}_H(E))$ in the following way: for an interval partition $\beta=(U_i, i \in I)$ such that $\sum_{i \in I} g(\text{Leb}(U_i))$ is finite we define $g^*(\beta)$ to be the interval partition $\gamma\coloneqq(V_i, i \in I)$ where we have $\text{Leb}(V_i)=g(\text{Leb}(U_i))$ and the order of the intervals in $\beta$ and $\gamma$ is consistent; when $\sum_{i \in I} g(\text{Leb}(U_i))$ is infinite we take $g^*(\beta)\coloneqq\emptyset$, the empty interval partition. Note that, $g^*(\mathcal{I}_H(E)) \subseteq \mathcal{I}_H(g(E))$. Also, as $g$ is an injection, $g^*$ restricted to $\{\beta \in \mathcal{I}_H(E): \sum_{U \in \beta} g(\text{Leb}(U))<\infty\}$ has an inverse on $\{\mathcal{I} \hspace{0.125cm} | \hspace{0.125cm} \mathcal{I} \subset g^*(\mathcal{I}_H)\}$ and we have that $(g^*)^{-1}(\mathcal{I}) = (g^{-1})^*(\mathcal{I})$. We also note that $g^*$ is a measurable function. To see this, first define the functions $(g_n, n \in \mathbb{N})$ by $g_n(x)=g(x)$ if $x>1/n$ and $g_n(x)=x$ otherwise and observe that $g^*$ is the limit of the sequence of the functions $(g_n^*, n \in \mathbb{N})$. The measurability of $g^*$ then follows by the fact that for all $n \in \mathbb{N}$ the functions $g_n^*$ are continuous. The proof of continuity of the functions $(g_n^*, n \in \mathbb{N})$ follows from the continuity of $g$ and that for any interval partition $\beta$ the interval partition differs from its image $g_n^*(\beta)$ on only finitely many interval widths (see also Lemma \ref{lemma g star is continuous}).
\par
The next step is to develop a simple Markov property. In Corollary \ref{corollary assumption b' starting interval partition} we have that for an interval partition $\beta \in \mathcal{I}_H(E)$ and a diffusion that satisfies Assumption \ref{assumption B'}, we can construct a scaffolding-and-spindles pair $(\mathbf{N}_\beta, \mathbf{X}_\beta)$ as in Definition \ref{definition starting interval partition}. We state a corollary here for diffusions under Assumption \ref{assumption b}.
\begin{corollary}\label{corollary starting g^*(I)}
    Consider a diffusion $\mathbf{Z}$ that satisfies Assumption \ref{assumption b} and let $\gamma \in g^*(\mathcal{I}_H(E))$. Consider a scaffolding-and-spindles pair started from $\gamma$, denoted as $(\mathbf{N}_\gamma, \mathbf{X}_\gamma)$ as in Definition \ref{definition starting interval partition}. Then we have that $\text{len}(\mathbf{N}_\gamma)=\sum_{U \in \gamma} \text{len}(\mathbf{N}_U)$ is finite a.s.
    \begin{proof}
        Let $\mathbf{Y}$ be the untransformed diffusion that satisfies Assumption \ref{assumption B'}, let $\gamma \in g^*(\mathcal{I}_H(E))$ and then let $\beta \in \mathcal{I}_H(E)$ such that $g^*(\beta)=\gamma$. As mentioned above, by Corollary \ref{corollary assumption b' starting interval partition} we can construct a scaffolding-and-spindles pair $(\mathbf{N}_\beta^{\mathbf{Y}}, \mathbf{X}_\beta^{\mathbf{Y}})$ for $\mathbf{Y}$ started from $\beta$ and $\text{len}(\mathbf{N}_\beta)$ is finite a.s. Therefore we can define a scaffolding-and-spindles pair for diffusion $\mathbf{Z}$ and started from $\gamma$ by
        \begin{equation*}
            \mathbf{N}_\gamma^{\mathbf{Z}} = \sum_{(t,f) \text{   atom of   } \mathbf{N}_\beta^{\mathbf{Y}}} \delta_{(t,g(f))},
        \end{equation*}
        where for a function $(f(x), x \in \mathbb{R})$ we write $g(f)$ to be the function $(g(f(x)), x \in \mathbb{R})$. Then define $\mathbf{X}_\gamma^{\mathbf{Z}}=\xi(\mathbf{N}_\gamma^{\mathbf{Z}})$ which is equal to $\mathbf{X}_\beta^{\mathbf{Y}}$ as the spatial transformation of the spindles does not change the scaffolding.
    \end{proof}
\end{corollary}
We now state a continuity result for an IP-evolution with a fixed initial state. For $N \in \mathcal{N}^{\text{sp}}_{\text{fin}}$, let $\zeta^+(N)\coloneqq\sup_{t \in \mathbb{R}}\xi_N(t)$.
\begin{lemma}\label{lemma continuity away from zero}
    Consider a diffusion $\mathbf{Z}$ that satisfies Assumption \ref{assumption b} and let $\gamma \in g^*(\mathcal{I}_H(E))$. Let $(\mathbf{N}_\gamma, \mathbf{X}_\gamma)$ be a scaffolding-and-spindles pair $\gamma$ as in Definition \ref{definition starting interval partition}. Then there exists a modification of $\overline{\text{skewer}}(\mathbf{N}_\gamma, \mathbf{X}_\gamma)$ that we call $(\gamma^y, y \geq 0)$ that is continuous for $y \in (0,\infty)$. Furthermore, for $y \in (0,\infty)$ fixed, we have $\gamma^y \in g^*(\mathcal{I}_H(E))$ a.s.
    \begin{proof}
        Let $n \in \mathbb{N}$ and $\gamma \in g^*(\mathcal{I}_H(E))$. From Corollary \ref{corollary starting g^*(I)} we can construct $(\mathbf{N}_U, U \in \gamma)$ as in Definition \ref{definition starting interval partition}. Recall that $\mathbf{N}_U \coloneqq \delta_{(0,f_U)} + \mathbf{N}|_{[0,T^{-\zeta(f_U)}]}$ where $\mathbf{N}$ is a Poisson random measure. The process $(\text{skewer}(y, \mathbf{N}|_{[0,T^{-\zeta(f_U)}]}, \xi(\mathbf{N}|_{[0,T^{-\zeta(f_U)}]})), y \in \mathbb{R})$ has a continuous modification by Theorem \ref{theorem continuity} and it follows that $\overline{\text{skewer}}(\mathbf{N}_U, \mathbf{X}_U)$ has a continuous modification. Also, from the fact that $\sum_{U \in \gamma} \text{len}(\mathbf{N}_U)$ is finite a.s., from Corollary \ref{corollary starting g^*(I)}, we have that for only finitely many $U \in \gamma$ that $\zeta^+(\mathbf{N}_U)>1/n$. Therefore we can conclude that there exists a continuous modification of $\overline{\text{skewer}}(\mathbf{N}_\gamma, \mathbf{X}_\gamma)$ on $(1/n,\infty)$. As this is true for all $n \in \mathbb{N}$, we have a continuous modification on $(0,\infty)$.
        \par
        It remains to show that for $y \in (0,\infty)$ fixed, we have $\gamma^y \in g^*(\mathcal{I}_H(E))$ a.s. Let $y>0$. For each $U \in \gamma$ we have that $\mathbf{N}_U$ is constructed by an initial incomplete spindle $f_U$ and then a stopped Poisson random measure $\mathbf{N}$ with stopping time $T^{-\zeta(f_U)}$. We claim that $\text{skewer}(y, \mathbf{N}_U, \mathbf{X}_U) \in g^*(\mathcal{I}_H(E))$ a.s. By the spatial homogeneity of the scaffolding-and-spindles pair and the independence of $f_U$ and $\mathbf{N}$ it is sufficient to show that $\text{skewer}(0, \mathbf{N}, \mathbf{X}) \in g^*(\mathcal{I}_H(E))$ a.s. We have that $\text{skewer}(0, \mathbf{N}, \mathbf{X})$ is in $\mathcal{I}_H(E)$ by the first part of Proposition \ref{corollary levy measure assumption b} and we can conclude that it is in fact in $g^*(\mathcal{I}_H(E))$ by again using the first part of Proposition \ref{corollary levy measure assumption b} with the identity transformation $g(x)=x$ (this tells us the untransformed spindles produce an interval partition and therefore the transformed spindles produce an interval partition in $g^*(\mathcal{I}_H(E))$ a.s.).
        \par
        We concatenate the interval partitions $\text{skewer}(y, \mathbf{N}_U, \mathbf{X}_U)$ for the finitely many of the $U \in \gamma$ with $\zeta^+(\mathbf{N}_U)>1/n$ to form $\text{skewer}(y,\mathbf{N}_\gamma, \mathbf{X}_\gamma)$ and conclude that $\text{skewer}(y,\mathbf{N}_\gamma, \mathbf{X}_\gamma) \in g^*(\mathcal{I}_H(E))$ a.s. and therefore the modification $\gamma^y \in g^*(\mathcal{I}_H(E))$ a.s.
    \end{proof}
\end{lemma}
We define the filtration $(\mathcal{F}_{\mathcal{I}}^y, y \geq 0)$ to be the right-continuous filtration generated by $(\gamma^y, y \in (0,\infty))$ on $\Omega\coloneqq\mathcal{C}((0,\infty), \mathcal{I}_H)$, the class of all paths in $\mathcal{I}_H$ continuous on $(0,\infty)$. We define a filtered probability space $(\Omega, \mathcal{F}_{\mathcal{I}}, (\mathcal{F}_{\mathcal{I}}^y, y \geq 0), \mathbb{P}_\gamma)$ where $\mathcal{F}_{\mathcal{I}}\coloneqq\bigcup_{y \geq 0} \mathcal{F}_{\mathcal{I}}^y$, and $\mathbb{P}_\gamma(\cdot)$ is the probability measure 
\begin{equation*}
    \mathbb{P}_\gamma(\cdot) = \mathbf{P}((\gamma^y, y \geq 0) \in \cdot)
\end{equation*}
where $(\gamma^y, y \geq 0)$ is the modification of $\overline{\text{skewer}}(\mathbf{N}_\gamma, \mathbf{X}_\gamma)$ that is continuous on $(0,\infty)$ that is given in the lemma above.
\begin{corollary}[Simple Markov property under Assumption \ref{assumption b}]\label{corollary simple markov property}
    Consider a diffusion under Assumption \ref{assumption b}. Let $\mu$ be a probability distribution on interval partitions that is supported on $g^*(\mathcal{I}_H(E))$, and define the probability measure $\mathbb{P}_\mu(\cdot)=\int_{\gamma \in g^*(\mathcal{I}_H(E))} \mathbb{P}_\gamma(\cdot)\mu(d\gamma)$. For any $y\in (0,\infty)$ let $\theta_y$ be the shift operator, and let $\eta, f: \mathcal{C}((0,\infty), \mathcal{I}_H) \rightarrow [0,\infty)$ be measurable where $\eta$ is measurable with respect to $\mathcal{F}_{\mathcal{I}}^y$. Then we have that 
    \begin{equation*}
        \mathbb{P}_\mu\big[\eta f \circ \theta_y\big] = \mathbb{P}_\mu \big[\eta \mathbb{P}_{\gamma^y}[f]\big].
    \end{equation*}
    \begin{proof}  
        This proof uses the same arguments as the analogous result \cite[Corollary 6.13]{construction}; specifically, by Lemma \ref{lemma continuity away from zero} we have that the sample paths a.s. are in $\Omega$ and the simple Markov property for the skewer process on a larger probability space is given in Proposition \ref{proposition simple markov skewer}. 
    \end{proof}
\end{corollary}
We now show the analogous result for diffusions under Assumption \ref{assumption c} and prove Theorem \ref{theorem interval partitions regular}. We let $\mathcal{I}_{\text{fin}} \subset \mathcal{I}_H$ be the interval partitions that consist of finitely many intervals and $\mathcal{I}_{\text{fin}}(g(E)) = \mathcal{I}_{\text{fin}} \cap \mathcal{I}_H(g(E))$. For any interval partition $\gamma \in \mathcal{I}_{\text{fin}}(g(E))$ we define $\mathbb{P}_\gamma$ as before and note that the scaffolding-and-spindles pair $(\mathbf{N}_\gamma, \mathbf{X}_\gamma)$ can always be constructed for such initial interval partitions $\gamma$. Let $\Omega_0\coloneqq \mathcal{C}([0,\infty), \mathcal{I}_H)$ be the space of continuous paths on $\mathcal{I}_H$ indexed by $[0,\infty)$. 
\begin{proof}[Proof of Theorem \ref{theorem interval partitions regular}]
    Let $\mathbf{Z}$ be a diffusion under Assumption \ref{assumption c}, and let $\gamma \in \mathcal{I}_{\text{fin}}(g(E))$. Let $(\mathbf{N}_\gamma,\mathbf{X}_\gamma)$ be a scaffolding-and-spindles pair started at $\gamma$ as defined in Definition \ref{definition starting interval partition}. For $y \geq 0$ let $\gamma^y=\text{skewer}(y, \mathbf{N}_\gamma, \mathbf{X}_\gamma)$. Then $(\gamma^y, y \geq 0)$ is continuous by the same reasoning as in Proposition \ref{proposition continuity assumption c}; for each $y$ we have only concatenated finitely many interval partitions.
    \par
    In this setting let $(\mathcal{F}_{\mathcal{I}}^y, y \geq 0)$ be the right-continuous filtration generated by $(\gamma^y, y \geq 0)$ on the smaller space $\Omega_0$, rather than on $\Omega$ as we had done previously. The simple Markov property for the skewer process from Proposition \ref{proposition simple markov skewer} applies in this setting as for fixed $y>0$ the interval partition $\gamma^y \in \mathcal{I}_{\text{fin}}(g(E))$ a.s. To see this, note that for each $U \in \gamma$ the SPLP $\mathbf{X}_U$ is of bounded variation and so only crosses a fixed level finitely many times a.s., \cite[Corollary VII.5]{bertoin1996levy}, and therefore $\gamma^y$ consists of finitely many intervals as it is the concatenation of finitely many interval partitions that each consist of finitely many intervals. To conclude the proof for the setting described in Theorem \ref{theorem interval partitions regular} we note that, for $n \in \mathbb{N}$, if we stop the PRM $\mathbf{N}$ of a scaffolding-and-spindles pair $(\mathbf{N},\mathbf{X})$ at the $n^{\text{th}}$ hitting time of $0$ then $\text{skewer}(0,\mathbf{N},\mathbf{X}) \in \mathcal{I}_{\text{fin}}(g(E))$ a.s.
\end{proof}
We used $\Omega_0 \subset \Omega$ under Assumption \ref{assumption c} whereas under Assumption \ref{assumption b} we were required to use $\Omega$ as we did not have a proof of continuity at time $0$. However, if we restrict the class of initial interval partitions then we can use $\Omega_0$ to state the simple Markov property under Assumption \ref{assumption b}, and we do this in the following subsection. This subsection is not necessary to prove Theorem \ref{theorem diffusion} but we include it to show the class of initial interval partitions for which we have a proof of continuity at time $0$.
\subsubsection{Continuity at time $0$ under Assumption \ref{assumption b}}\label{section continuity initial partition}
To prove the results of this section we define certain classes of interval partitons. For any $q>0$ and any class of interval partitions $\mathcal{I} \subseteq \mathcal{I}_H$ we define 
\begin{equation*}
    \mathcal{I}^{(q)} \coloneqq \{ \beta = (U_i, i \in I) \in \mathcal{I}: \sum_{i \in I} \text{Leb}(U_i)^q < \infty\}.
\end{equation*}
Note that for $q \geq 1$ we have that $\mathcal{I}^{(q)} = \mathcal{I}$. We now state a corollary of Theorem \ref{theorem continuity}.
\begin{lemma}\label{lemma restriction to I^q space}
    Consider a diffusion that satisfies Assumption \ref{assumption B'} with a stopped scaffolding-and-spindles pair $(\mathbf{N}, \mathbf{X})$ on some probability space $(\Omega, \mathcal{A}, \mathbf{P})$, and let $q>\alpha^+$. Then, for $(\text{skewer}(y, \mathbf{N}, \mathbf{X}), y \geq 0)$, there exists a modification $(\beta^y, y \geq 0)$ which is continuous on $(\mathcal{I}^{(q)}_H(E), d'_H)$.
    \begin{proof}
        We first note that for $q \in (0,1)$ we can apply Theorem \ref{theorem continuity} with the diffusion $\mathbf{Y}$ with the transformation $h(x)=x^q$, to produce the continuous IP-evolution $(\phi^y, y \geq 0)$. Also, note that the map $(h^{-1})^*$ is continuous on the space of interval partitions $\mathcal{I}_H$ (this ultimately follows because the map $h^{-1}$ is Lipschitz on $[0,1]$, see also Lemma \ref{lemma g star is continuous}). Therefore $((h^{-1})^*(\phi^y), y \geq 0)$ is $\mathbf{P}$-a.s. continuous in $\mathcal{I}_H$. 
        \par
        Now, $(\phi^y), y \geq 0)$ is a modification of $h^*(\overline{\text{skewer}}(\mathbf{N}, \mathbf{X}))$ and therefore $((h^{-1})^*(\phi^y), y \geq 0)$ is a modification of $\overline{\text{skewer}}(\mathbf{N}, \mathbf{X})$. Also, $(\beta^y, y \geq 0)$ is a modification of $\overline{\text{skewer}}(\mathbf{N}, \mathbf{X})$ and so is also a modification of $((h^{-1})^*(\phi^y), y \geq 0)$. As both $(\beta^y, y \geq 0)$ and $((h^{-1})^*(\phi^y), y \geq 0)$ are both continuous processes they are therefore indistinguishable from one another, and therefore $(\beta^y, y \geq 0)$ is a continuous process in $\mathcal{I}_H^{(q)}(E)$ a.s.     
    \end{proof}
\end{lemma}
As in the previous sections, we consider a diffusion $\mathbf{Z}$ under Assumption \ref{assumption b}, which is equal in distribution to $g(\mathbf{Y})$ where $\mathbf{Y}$ is a diffusion under Assumption \ref{assumption B'} and $g:E \rightarrow g(E)$ is a continuous bijection as in Assumption \ref{assumption b}. Let $q_0$, $\alpha^+$, $\alpha^-$ be the associated constants and let $q \in (\alpha^+, q_0)$. Finally, let $m_{\mathbf{Y}}$ and $m_{\mathbf{Z}}$ be the speed measure densities of $\mathbf{Y}$ and $\mathbf{Z}$. We now state a corollary of Lemma \ref{lemma diversity is local time} and Lemma \ref{lemma restriction to I^q space}.
\begin{corollary}\label{corollary space of continuous initial states}
    Consider a stopped scaffolding-and-spindles pair $(\mathbf{N}, \mathbf{X})$ for the diffusion $\mathbf{Z}$ which satisfies Assumption \ref{assumption b}, and let $(\gamma^y, y \geq 0)$ be a continuous modification of $\overline{\text{skewer}}(\mathbf{N},\mathbf{X})$. Then for any fixed $y \geq 0$ we have that
    \begin{equation*}
        \gamma^y \in g^*(\mathcal{I}^{(q)}_{m_{\mathbf{Y}}}(E)) = \Big\{\gamma = (V_i, i \in I) \in g^*(\mathcal{I}_{m_{\mathbf{Y}}}(E)) : \sum_{i \in I} g^{-1}(\text{Leb}(V_i))^q<\infty \Big\} \hspace{0.25cm} \mathbf{P}-\text{a.s.}
    \end{equation*}
\begin{proof}
    Let $\Tilde{\mathbf{N}}$ be the Poisson random measure of untransformed spindles. The $m_{\mathbf{Y}}$-diversity of $\overline{\text{skewer}}(\Tilde{\mathbf{N}}, \mathbf{X})$ at any fixed level has already been established in Lemma \ref{lemma diversity is local time}. This shows that $\text{skewer}(y,\mathbf{N},\mathbf{X}) \in g^*(\mathcal{I}_{m_\mathbf{Y}}(E))$ a.s. and so $\gamma^y \in g^*(\mathcal{I}_{m_\mathbf{Y}}(E))$.
    \par
    It remains to show that $\gamma^y \in g^*(\mathcal{I}^{(q)}_H(E))$. To do this we note that the process $(\gamma^y, y \geq 0)$ is a modification of $g^*(\overline{\text{skewer}}(\Tilde{\mathbf{N}},\mathbf{X}))$. Therefore $((g^{-1})^*(\gamma^y), y \geq 0)$ is a modification of $\overline{\text{skewer}}(\Tilde{\mathbf{N}},\mathbf{X})$. By Lemma \ref{lemma restriction to I^q space} $\overline{\text{skewer}}(\Tilde{\mathbf{N}},\mathbf{X})$ has a continuous modification that takes values in $\mathcal{I}^{(q)}_H(E)$, and the result follows as for fixed $y \geq 0$ we have that $(g^{-1})^*(\gamma^y)=\text{skewer}(y,\Tilde{\mathbf{N}}, \mathbf{X})$ a.s.
\end{proof}
\end{corollary}
The following lemma concerns the initial spindles of the IP-evolution, i.e. the spindles that are not equal to $0$ at time $0$. 
\begin{lemma}\label{lemma continuity initial spindles}
    Consider a diffusion that satisfies Assumption \ref{assumption B'}. Let $(a_i, i \in I)$ be a potentially infinite sequence of positive numbers such that $\sum_{i \in I}a_i^{q \wedge 1}<\infty$. Consider a family of independent diffusions $(f_i, i \in I)$ distributed as independent diffusions started from $(a_i, i \in I)$ stopped at $0$ on some common probability space $(\Omega, \mathcal{A}, \mathbf{P})$, i.e. $f_i \sim \mathbb{Q}^0_{a_i}$. Then we have that $\sum_{i \in I} A(f_i)^q < \infty$ $\mathbf{P}$-a.s.
    \begin{proof}
        By the coupling given in \eqref{equation coupling of diffusions} we have for $0<a<y<\varepsilon_0$ that $\mathbb{Q}^0_a(A>y) \leq {^{-2\alpha^-}}\mathbb{Q}^0_a(A>y)$ which in turn is equal to $(a/y)^{1+\alpha^-}$ (by \eqref{equation scale function property} and as $x^{1+\alpha^-}$ for $x \in [0,\infty)$ is a scale function for a squared Bessel process with dimension $-2\alpha^-$). Therefore we have
        \begin{equation*}
            \sum_{i \in I:a_i<\varepsilon_0} \mathbf{P}(A(f_i)>\varepsilon_0) \leq \varepsilon_0^{-1-\alpha^-} \sum_{i \in I} a_i^{1+\alpha^-} < \infty,
        \end{equation*}
        and so by the Borel-Cantelli lemma and by the fact that only finitely many $a_i$ are such that $a_i\geq \varepsilon_0$, we have that a.s. only finitely many of the incomplete spindles $f_i$ will have amplitude greater than $\varepsilon_0$.
        \par
        Furthermore, again with the upper bound for $\mathbb{Q}^0_a(A>y)$ stated above, we have
        \begin{equation*}
            \mathbf{E}(A(f_i)\mathbbm{1}_{A(f_i)<\varepsilon_0}) = \int_{y=0}^{\varepsilon_0} \mathbf{P}(y \leq A(f_i) < \varepsilon) dy.
        \end{equation*}
        This is less than or equal to 
        \begin{equation*}
            \int_{y=0}^{\varepsilon_0} \mathbb{Q}_{a_i}^0(A>y) dy = 
            \int_{y=0}^{\varepsilon_0} \Big((a_i/y)^{1+\alpha^-} \wedge 1 \Big) dy = \Big(\frac{1}{\alpha^-} + 1\Big) a_i - \frac{a_i^{1+\alpha^-}}{\alpha^-\varepsilon_0^{\alpha^-}} \leq \Big(\frac{1}{\alpha^-} + 1\Big) a_i.
        \end{equation*}
        Therefore, when $q<1$, we have that
        \begin{equation*}
            \mathbf{E}(A(f_i)^q\mathbbm{1}_{A(f_i)<\varepsilon_0}) \leq \Big(\mathbf{E}(A(f_i)\mathbbm{1}_{A(f_i)<\varepsilon_0})\Big)^q \leq C_q \Big(\frac{1}{\alpha^-} + 1\Big)^q a_i^q
        \end{equation*}
        where $C_q$ is the H\"older constant from Assumption \ref{assumption b}. Therefore we have that $\sum_{i \in I} \mathbf{E}(A(f_i)^q:A(f_i)<\varepsilon_0)$ is finite, and so $\sum_{i \in I:A(f_i)<\varepsilon_0} A(f_i)^q$ is finite $\mathbf{P}$-a.s.
    \end{proof}
\end{lemma}
In the next lemma we give a result that compares certain interval partitions with $\text{skewer}(0,\mathbf{N}, \mathbf{X})$ for a particular stopped scaffolding-and-spindles pairs $(\mathbf{N}, \mathbf{X})$, which we then use in a coupling argument in the proposition after to prove continuity of a IP-evolutions with certain initial interval partitions. 
\begin{lemma}[Interval partition domination]\label{lemma dominate by matching}
Consider a diffusion $\mathbf{Z}$ that satisfies Assumption \ref{assumption b} and let $\gamma \in \overline{\mathcal{I}}_{m_{\mathbf{Z}}}(g(E))$. Recall the notation $\mathcal{D}$ and $\overline{\mathcal{D}}$ from Definition \ref{definition m-diversity} and let $(\mathbf{N},\mathbf{X})$ be a scaffolding-and-spindles pair stopped at $\tau^0(\overline{\mathcal{D}}^{m_{\mathbf{Z}}}_\gamma(\infty)+1)$ (the inverse local time of $\mathbf{X}$ at level $0$), and let $\Tilde{\gamma}\coloneqq\text{skewer}(0,\mathbf{N}, \mathbf{X})$. Firstly we have that
\begin{equation*}
    \mathcal{D}_{\Tilde{\gamma}}^{m_\mathbf{Z}}(\infty) = \overline{\mathcal{D}}_{\gamma}^{m_\mathbf{Z}}(\infty) + 1 \text{   a.s.}
\end{equation*}
Also, with positive probability we have that there exists a matching between $\gamma$ and $\Tilde{\gamma}$ such that every interval of $\gamma$ is matched with a larger interval of $\Tilde{\gamma}$ (but that does not necessarily respect the total-order of the intervals). We say that $\Tilde{\gamma}$ dominates $\gamma$ by matching in this case. Finally we have that $g^*(\overline{\mathcal{I}}_{m_{\mathbf{Y}}}(E)) = \overline{\mathcal{I}}_{m_{\mathbf{Z}}}(g(E))$ and also $g^*(\mathcal{I}_{m_{\mathbf{Y}}}(E)) = \mathcal{I}_{m_{\mathbf{Z}}}(g(E))$.
\begin{proof}
    This first part of this lemma generalises \cite[Lemma 6.8]{construction} to interval partitions with finite $m_{\mathbf{Z}}$-diversity and the proof is the same. The result further generalises to any interval partition with a finite limsup-$m_{\mathbf{Z}}$-diversity as intervals can be added to an interval partition with finite limsup-$m_{\mathbf{Z}}$-diversity to produce an interval partition with a finite $m_{\mathbf{Z}}$-diversity and, therefore, the domination result carries over.
    \par
    For the second part first note that if a pair of interval partitions $\hat{\beta}$, $\hat{\gamma} \neq \emptyset$ satisfy $g^*(\hat{\beta})=\hat{\gamma}$ then $\hat{\beta} \in \mathcal{I}_{m_{\mathbf{Y}}}(E)$ if and only if $\hat{\gamma} \in \mathcal{I}_{m_{\mathbf{Z}}}(g(E))$ and also $\mathcal{D}_{\hat{\beta}}^{m_\mathbf{Y}}(\infty)=\mathcal{D}_{\hat{\gamma}}^{m_\mathbf{Z}}(\infty)$ (and similarly for limsup diversities). Now, for $\Tilde{\gamma}$ defined above there exists $\Tilde{\beta}$ such that $g^*(\Tilde{\beta})=\Tilde{\gamma}$ by the aggregate health summability in the transformed and untransformed cases shown in Proposition \ref{corollary levy measure assumption b}. Therefore we have for the $\gamma \in \overline{\mathcal{I}}_{m_{\mathbf{Z}}}(g(E))$ defined above that $\gamma \in g^*(\overline{\mathcal{I}}_{m_{\mathbf{Y}}}(E))$. This can be seen as $\sum_{V \in \gamma} g^{-1}(\text{Leb}(V)) \leq \norm{\Tilde{\beta}}$. The generality of $\gamma$ here means that $\overline{\mathcal{I}}_{m_{\mathbf{Z}}}(g(E)) \subseteq g^*(\overline{\mathcal{I}}_{m_{\mathbf{Y}}}(E))$. The result follows as it is immediate that $g^*(\overline{\mathcal{I}}_{m_{\mathbf{Y}}}(E)) \subseteq \overline{\mathcal{I}}_{m_{\mathbf{Z}}}(g(E))$. The other case $g^*(\mathcal{I}_{m_{\mathbf{Y}}}(E)) = \mathcal{I}_{m_{\mathbf{Z}}}(g(E))$ follows by the same reasoning.
\end{proof}
\end{lemma}
\begin{proposition}\label{proposition continuous from fixed initial state}
    Let $\gamma = (V_i, i \in I) \in g^*(\overline{\mathcal{I}}^{(q)}_{m_\mathbf{Y}}(E))$ and let $(\mathbf{N}_\gamma, \mathbf{X}_\gamma)$ be a scaffolding-and-spindles pair for $\gamma$ as in Definition \ref{definition starting interval partition}. Then there exists a modification of $\overline{\text{skewer}}(\mathbf{N}_\gamma, \mathbf{X}_\gamma)$ that is continuous on $[0,\infty)$.
    \begin{proof}
    As in the proof of Theorem \ref{theorem continuity} first assume that $g(E)$ is closed. Let $\gamma \in g^*(\overline{\mathcal{I}}_{m_\mathbf{Y}}(E)) = \overline{\mathcal{I}}_{m_\mathbf{Z}}(g(E))$ and let $\Tilde{\gamma}$ be a random interval partition such that $\mathcal{D}_{\Tilde{\gamma}}^{m_\mathbf{Z}}(\infty) = \overline{\mathcal{D}}_{\gamma}^{m_\mathbf{Z}}(\infty) + 1$ and that dominates $\gamma$ as in Lemma \ref{lemma dominate by matching}. We describe a coupling for a scaffolding-and-spindles pair $(\mathbf{N}_\gamma, \mathbf{X}_\gamma)$ and another stopped scaffolding-and-spindles pair $(\mathbf{N}, \mathbf{X})$. 
    \par
    Let $(\Tilde{V}_i, i \in \Tilde{I})$ be the intervals of $\Tilde{\gamma}$ and let $I \subset \Tilde{I}$ so that $(V_i, i \in I)$ denotes the matching intervals of $\gamma$. For each $i \in I$, we consider an independent Brownian motion $(B^i_y, y \geq 0)$ and we define $f_{\Tilde{V}_i}$ and $f_{V_i}$ to be strong solutions to the stochastic differential equations
    \begin{equation*}
        df_{\Tilde{V}_i} (y) = \text{Leb}(\Tilde{V}_i) + \mu_{\mathbf{Z}}(f_{\Tilde{V}_i} (y))dy + \sigma_{\mathbf{Z}}(f_{\Tilde{V}_i} (y))dB^i_y; \hspace{0.5cm} df_{V_i} (y) = \text{Leb}(V_i) + \mu_{\mathbf{Z}}(f_{V_i} (y))dy + \sigma_{\mathbf{Z}}(f_{V_i} (y))dB^i_y,
    \end{equation*}
    and we stop $f_{\Tilde{V}_i}$ and $f_{V_i}$ when they hit $0$. With this construction we have that $(f_{\Tilde{V}_i}, i \in I)$ are independent diffusions with the same dynamics and with initial states $(\Tilde{V}_i, i \in I)$, and similarly for $(f_{V_i}, i \in I)$. We also have for each $i \in I$ that $\zeta(f_{V_i})\leq \zeta(f_{\Tilde{V}_i})$ and for each $y>0$ that $f_{V_i}(y)\leq f_{\Tilde{V}_i}(y)$. We then define a pair of point measures $(\mathbf{N}, \mathbf{N}_\gamma)$. For $z>0$, let $\mathbb{Q}^0_z$ denote the law of the diffusion $\mathbf{Z}$ started from $z$, and let $\nu$ be the Pitman-Yor excursion measure of $\mathbf{Z}$. Consider independent Poisson random measures $(\mathbf{N}_i, i \in I)$ with intensity $\text{Leb} \otimes \nu$ and scaffolding $\xi(\mathbf{N}_i)=\mathbf{X}_i$ for $i \in I$; we use this to define
    \begin{equation*}
        \mathbf{N}_{\Tilde{V}_i} \coloneqq \delta(0,f_{\Tilde{V}_i}) + \mathbf{N}_i|_{[0,\Tilde{T}_i]} \hspace{0.25cm} \text{   where   } \hspace{0.25cm} \Tilde{T}_i\coloneqq \inf\{t>0: \mathbf{X}_i(t) = -\zeta(f_{\Tilde{V}_i})\}.
    \end{equation*}
    We further define $\mathbf{N}_\gamma \coloneqq \bigstar_{i \in I} \mathbf{N}_{V_i}$ where $\mathbf{N}_{V_i}$ is given by
    \begin{equation}\label{equation construction remaining spindles}
        \mathbf{N}_{V_i} \coloneqq \delta(0,f_{V_i}) + \mathbf{N}_i|^\leftarrow_{(T_i,\Tilde{T_i}]} \hspace{0.25cm} \text{   with   } \hspace{0.25cm} T_i \coloneqq \inf\{ t \geq 0: \mathbf{X}_i(t) = \zeta(f_{V_i})-\zeta(f_{\Tilde{V}_i}) \}
    \end{equation}
    and $\Tilde{T_i}$ is as stated above. Let $T=\sum_{i \in I}T_i$ and note $T$ is a $\mathbf{P}$-a.s. finite random variable. Finally let $\mathbf{N} \coloneqq\bigstar_{i \in I} \mathbf{N}_{\Tilde{V}_i}$ and $\mathbf{N}_{\gamma} \coloneqq\bigstar_{i \in I} \mathbf{N}_{V_i}$.
    \par
    The rest of the proof uses the same ideas as the proof of Theorem \ref{theorem continuity}. Let $q$, $\kappa>0$ such that $\alpha^++\kappa<q<q_0$ where $\alpha^+$, $q_0$ are as in Assumption \ref{assumption b} and let $\rho \coloneqq (q-\alpha^+-\kappa) \wedge q/2$. We recall $\mathbf{N}_{M,\varepsilon_0}$ from \eqref{equation n with m and epsilon} and define $\mathbf{N}_{\gamma, M,\varepsilon_0}$:
    \begin{equation*}
        \mathbf{N}_{M,\varepsilon_0}\coloneqq \mathbbm{1}_{\sup_{u \in \mathbb{R}}\ell^u(T)<M} \sum_{\substack{(t,f) \text{   atom of   } \\ \mathbf{N}: A(f)< \varepsilon_0}} \delta_{(t,f)}; \hspace{0.25cm} \mathbf{N}_{\gamma, M,\varepsilon_0}\coloneqq \mathbbm{1}_{\sup_{u \in \mathbb{R}}\ell^u(T)<M} \sum_{\substack{(t,f) \text{   atom of   } \mathbf{N}_\gamma: \\  A(f)< \varepsilon_0, f(0)=0}} \delta_{(t,f)}
    \end{equation*}
    where $\ell^u$ for $u \in \mathbb{R}$ denotes the local time of $\xi(\mathbf{N})$. We note that the set of atoms of $\mathbf{N}_{\gamma, M,\varepsilon_0}$ is a subset of the set of atoms of $\mathbf{N}_{M,\varepsilon_0}$ and so we can use the bound in \eqref{equation sum of differences x and y} as $\mathbf{N}$ is a Poisson random measure stopped at $\tau^0(\overline{\mathcal{D}}^{m_{\mathbf{Z}}}_\gamma(\infty)+1)$. The constant in front of the bound in \eqref{equation sum of differences x and y} will change from $A$ to $A/p$ where $p\coloneqq \mathbf{P}(\Tilde{\gamma} \text{   dominates   } \gamma)>0$ where we defined domination in Lemma \ref{lemma dominate by matching} (note that the fact that the domination of $\gamma$ by $\Tilde{\gamma}$ does not respect the order of the intervals in the interval partition does not impact the upper bound for the sum).
    \par
    At this point we need to incorporate the initial spindles and the spindles with amplitude greater than $\varepsilon_0$. In the proof of Theorem \ref{theorem continuity}, we only had the spindles with amplitude greater than $\varepsilon_0$ to add in, and as there were only finitely many of these we could concatenate finitely many continuous IP-evolutions to finish the proof. In this setting however as we have potentially infinitely many initial spindles we need to take care to incorporate these. 
    \par
    We adapt the proof of the Kolmogorov-\v{C}entsov continuity criterion from \cite[Theorem 1.2.1]{revuz2013continuous}. We define
    \begin{equation*}
        D_m \coloneqq \{i2^{-m}: i \in \mathbb{N} \cap [0,2^m) \} \text{   for   } m \in \mathbb{N}; \hspace{0.5cm} D \coloneqq \bigcup_{m \in \mathbb{N}} D_m
    \end{equation*}
    and
    \begin{equation*}
        M_\rho \coloneqq \sup_{x, y \in D, x \neq y} |x-y|^{-\rho} \Big\{ \sum_{\substack{(t,f) \text{   atom of   } \mathbf{N}_{M,\varepsilon_0}:, \\ \mathbf{X}_{t-}<x<\mathbf{X}_t \text{   or   } \\ \mathbf{X}_{t-}<y<\mathbf{X}_t}} 
        \Big|g\big(f(y-\mathbf{X}_{t-})\big) - g\big(f(x-\mathbf{X}_{t-})\big)\Big|\Big\}.
    \end{equation*}
    Then it follows from the proof in \cite[Theorem 1.2.1]{revuz2013continuous} that $M_\rho$ is finite $\mathbf{P}$-a.s. as we have shown the bound in \eqref{equation sum of differences x and y} with $\mathbf{N}_{\gamma, M,\varepsilon_0}$ in place of $\mathbf{N}_{M,\varepsilon_0}$. Let $\mathbf{N}_{\gamma,M}\coloneqq \mathbbm{1}_{\sup_{u \in \mathbb{R}}\ell^u(T)<M}\mathbf{N}_\gamma$. From this point we show that $\overline{\text{skewer}}(\mathbf{N}_{\gamma,M},\mathbf{X}_\gamma)$ is $d'_H$-continuous on $D$. Let $\varepsilon>0$, let $K \in \mathbb{N}$, $i_1, \dots, i_K \in I$ such that $\sum_{i \in I: i \neq i_1, \dots, i_K} A(f_{V_i}) < \varepsilon/4$, and let $K' \in \mathbb{N}$ and $(t_1, f_1), \dots, (t_{K'},f_{K'})$ be the atoms of $\mathbf{N}_\gamma$ with amplitude greater than or equal to $\varepsilon_0$. There exists $\delta_1>0$ such that for $x$, $y \in D$ with $|x-y|<\delta_1$ we have
    \begin{equation*}
        \sum_{j=1}^K \Big|g(f_{V_j}(y))-g(f_{V_j}(x))\Big| < \frac{\varepsilon}{4}; \hspace{0.5cm} \sum_{j=1}^{K'}  \Big|g\big(f_j(y-\mathbf{X}_{t_j-})\big) - g\big(f_j(x-\mathbf{X}_{t_j-})\big)\Big| < \frac{\varepsilon}{4}.
    \end{equation*}
    Furthermore, for $x$, $y \in D$ with $|x-y|<\delta_2 \coloneqq (\varepsilon/4M_\rho)^{1/\rho}$ we have
    \begin{equation*}
        \sum_{\substack{(t,f) \text{   atom of   } \mathbf{N}_{M,\varepsilon_0}:, \\ \mathbf{X}_{t-}<x<\mathbf{X}_t \text{   or   } \\ \mathbf{X}_{t-}<y<\mathbf{X}_t}} 
        \Big|g\big(f(y-\mathbf{X}_{t-})\big) - g\big(f(x-\mathbf{X}_{t-})\big)\Big| < \frac{\varepsilon}{4}
    \end{equation*}
    Let $\delta= \min \{ \delta_1, \delta_2\}$. Then for all $x$, $y \in D$ with $|x-y|<\delta$ the sum of $|g(f(y-\mathbf{X}_{t-})) - g(f(x-\mathbf{X}_{t-}))|$ over all atoms $(t,f)$ in $\mathbf{N}_{\gamma,M}$ is less than $\varepsilon$ $\mathbf{P}$-a.s. This is an upper bound for the distortion between $\text{skewer}(x,\mathbf{N}_{\gamma,M},\mathbf{X}_\gamma)$ and $\text{skewer}(y,\mathbf{N}_{\gamma,M},\mathbf{X}_\gamma)$ for the correspondence that matches intervals with each other if they both are generated by the same spindle. This is a correspondence as we have only finitely many spindles with lifetime greater than $|x-y|$; it is the same correspondence used in Theorem \ref{theorem continuity}. We can take $M$ to infinity at this stage as by Proposition \ref{corollary bi-continuous} the local time of $\mathbf{X}$ is continuous and therefore bounded. Therefore we have that $\overline{\text{skewer}}(\mathbf{N}_\gamma,\mathbf{X}_\gamma)$ is continuous on $D$ $\mathbf{P}$-a.s. To finish, we can define for $y \in [0,\infty)$
    \begin{equation*}
        \gamma^y \coloneqq \lim_{\substack{s \rightarrow y \\ s \in D}} \text{skewer}(s,\mathbf{N}_\gamma, \mathbf{X}_\gamma)
    \end{equation*}
    and $(\gamma^y, y \geq 0)$ is a continuous modification of $\overline{\text{skewer}}(\mathbf{N}_\gamma, \mathbf{X}_\gamma)$. This modification is in the appropriate subspace $(\mathcal{I}_H(g(E)),d'_H)$ as this subspace is closed. Now assume $g(E)$ is open. To complete the argument in this case we avoid taking modifications of interval partition evolutions with intervals great than $\varepsilon_0$. As in Theorem \ref{theorem continuity} there are only finitely many spindles with amplitude greater than $\varepsilon_0$ and so we can consider stopping times in \eqref{equation stopping times Sn} with $\mathbf{N}_\gamma$ in place of $\mathbf{N}$. Denoting these stopping times by $(S_n^\gamma, n \in \mathbb{N})$ we can obtain finitely many skewer processes $(\overline{\text{skewer}}(\mathbf{N}_\gamma|_{(S_n^\gamma, S_{n+1}^\gamma)}), n \leq P)$ for some $P \in \mathbb{N}$ and concatenate these with the skewer processes generated by the individual spindles with amplitude larger than $\varepsilon_0$ as in \eqref{equation concatenate}.
    \end{proof}
\end{proposition}
\begin{remark}
    The previous proposition shows path continuity at $0$ for certain initial states but not for all possibles ones. This is a partial response to the claim just after \cite[Corollary 6.19]{construction} that states \say{We believe that $d_H'$-path-continuity extends to time $0$} which is in reference to their branching interval partition evolutions generated from self-similar diffusions. As their self-similar diffusions satisfy our Assumption \ref{assumption b} we have confirmed that this belief is correct for at least some initial states. 
\end{remark}
\begin{proposition}[Simple Markov property under Assumption \ref{assumption b} on $\Omega_0$]\label{proposition simple markov property}
    Consider a diffusion under Assumption \ref{assumption b}. Let $\gamma \in g^*(\overline{\mathcal{I}}^{(q)}_{m_\mathbf{Y}}(E))$ and let $(\mathbf{N}_\gamma, \mathbf{X}_\gamma)$ be a scaffolding-and-spindles pair with initial state $\text{skewer}(0,\mathbf{N}_\gamma, \mathbf{X}_\gamma)=\gamma$ a.s. as in Definition \ref{definition starting interval partition}. Let $(\gamma^y, y \geq 0)$ be a modification of $\overline{\text{skewer}}(\mathbf{N}_\gamma, \mathbf{X}_\gamma)$ that is continuous on $[0,\infty)$ which exists by Proposition \ref{proposition continuous from fixed initial state}. Recall the filtration $(\mathcal{F}_{\mathcal{I}}^y, y \geq 0)$ from Corollary \ref{corollary simple markov property}. We define $(\Omega_0, \mathcal{F}, (\mathcal{F}_{\mathcal{I}}^y, y \geq 0), \mathbb{P}_\gamma)$, a filtered probability space, where $\Omega_0$ is the set of all paths over $\mathcal{I}_H$ that are continuous on $[0,\infty)$ and with initial states in $g^*(\overline{\mathcal{I}}^{(q)}_{m_\mathbf{Y}})$, $\mathcal{F}\coloneqq\bigcup_{y \geq 0} \mathcal{F}_{\mathcal{I}}^y$, and also $\mathbb{P}_\gamma(\cdot)\coloneqq\mathbf{P}((\gamma^y, y \geq 0) \in \cdot)$. Finally, let $\mu$ be a probability distribution on interval partitions that is supported on $g^*(\overline{\mathcal{I}}^{(q)}_{m_\mathbf{Y}}(E))$, and define the probability measure $\mathbb{P}_\mu(\cdot)=\int_{\gamma \in g^*(\mathcal{I}_H(E))} \mathbb{P}_\gamma(\cdot)\mu(d\gamma)$. Then $(\gamma^y, y \geq 0)$ is a simple Markov process on $(\Omega_0, \mathcal{F}, (\mathcal{F}_{\mathcal{I}}^y, y \geq 0), \mathbb{P}_\mu)$.
    \begin{proof}
        This proof uses the same arguments as the analogous result \cite[Corollary 6.13]{construction}. In particular, we have continuity of the sample paths from Proposition \ref{proposition continuous from fixed initial state} with initial states in $g^*(\overline{\mathcal{I}}^{(q)}_{m_\mathbf{Y}}(E))$ as well as the simple Markov property for the skewer process given in Proposition \ref{proposition simple markov skewer}. Furthermore, at each level we have, for $y \in (0,\infty)$, that the $\text{skewer}(y, \mathbf{N}, \mathbf{X}) \in g^*(\mathcal{I}^{(q)}_{m_\mathbf{Y}}(E)) \subset g^*(\overline{\mathcal{I}}^{(q)}_{m_\mathbf{Y}}(E))$ a.s. by Corollary \ref{corollary space of continuous initial states}. Therefore, for any fixed $y \in (0,\infty)$, the continuous modification takes values in $g^*(\overline{\mathcal{I}}^{(q)}_{m_\mathbf{Y}}(E))$ a.s. as required. 
    \end{proof}
\end{proposition}
\subsection{Strong Markov property and proof of Theorem \ref{theorem diffusion}}\label{subsection strong markov}
In this section we assume that the diffusion satisfies Strong Assumption \ref{assumption b}, that is Assumption \ref{assumption b} along with the assumption that either $\liminf_{x \rightarrow 0} g(x)/x >0$ or $\limsup_{x \rightarrow 0} g(x)/x <\infty$. Note that the identity transformation $g(x)=x$ satisfies these conditions, and therefore Strong Assumption \ref{assumption b} is more restrictive than Assumption \ref{assumption b} but more general than Assumption \ref{assumption B'}. Before we begin the proof we give a class of transformations $g$ that do not satisfy $B5$ of Strong Assumption \ref{assumption b} to show that this really is a restriction.
\begin{example}[Transformation $g$ that does not satisfy $B5$]\label{example g counterexample}
    Consider a function $g:[0,1/32] \rightarrow [0,1/32]$ such that for all $n \geq 3$ we have:
    \begin{equation*}
        g(2^{-2^n})=2^{-2^n}; \hspace{0.25cm} g(2^{-2^n-1})=2^{-2^n-n}; \hspace{0.25cm} g(2^{-2^{n+1}+1}) = 2^{-2^{n+1}+n},
    \end{equation*}
    and where for all $n \geq 3$ $x \in [2^{-2^{n+1}}, 2^{-2^{n+1}+1}]$, $[2^{-2^{n+1}+1},2^{-2^n-1}]$, $[2^{-2^n-1}, 2^{-2^n}]$ that $g$ is smooth and strictly increasing and $g(0)=0$. Then $g$ is a smooth continuous strictly increasing bijection and moreover
    \begin{equation*}
        \liminf_{x \downarrow 0} \frac{g(x)}{x} \leq \lim_{n \rightarrow \infty} \frac{g(2^{-2^n-1})}{2^{-2^n-1}} = \lim_{n \rightarrow \infty} 2^{-n+1} = 0; \hspace{0.25cm} \limsup_{x \downarrow 0} \frac{g(x)}{x} \geq \lim_{n \rightarrow \infty}\frac{g(2^{-2^{n+1}+1})}{2^{-2^{n+1}+1}} = \lim_{n \rightarrow \infty} 2^{n-1} = \infty.
    \end{equation*}
\end{example}
We state key results about the transformations $g^*$ and $(g^{-1})^*$ which apply under Strong Assumption \ref{assumption b}.
\begin{lemma}\label{lemma g star is continuous}
    Let $g:E \rightarrow g(E)$ be a strictly increasing continuous function with $g(0)=0$ where $E\subseteq[0,\infty)$ is of the form $E=[0,c)$ or $E=[0,c]$ for some $c>0$ and $g(E)\subseteq[0,\infty)$. If $\limsup_{x \rightarrow 0} g(x)/x < \infty$ then the map $g^*:\mathcal{I}_H(E) \rightarrow \mathcal{I}_H(g(E))$ is injective and continuous. If instead $\liminf_{x \rightarrow 0} g(x)/x >0$ then the map $(g^{-1})^*:\mathcal{I}_H(g(E)) \rightarrow \mathcal{I}_H(E)$ is injective and continuous.
    \begin{proof}
        First suppose that $\limsup_{x \rightarrow 0} g(x)/x < \infty$. Let $\varepsilon>0$, let $\beta \in \mathcal{I}_H(E)$ and let $(a_j)_{j \in \mathbb{N}}$ be the interval lengths of $\beta$ in a non-decreasing order and let $(U_j, j \in \mathbb{N})$ be the corresponding intervals of $\beta$. If $g(c)$ is finite let $D=(a_1+g(c))/2$ and if $g(c)$ is infinite let $D=a_1+1$. Let $C\coloneqq\sup_{x \in (0,D]} g(x)/x$.  We note that $\sum_{j \in \mathbb{N}}g(a_j) \leq C\sum_{j \in \mathbb{N}}a_j$ which is finite so $g^*(\beta)$ is a non-trivial interval partition. The fact that $g^*$ is injective then follows immediately from $g$ being injective.
        \par
        Let $K \in \mathbb{N}$ be large enough so that $\sum_{j>K} a_j < \varepsilon/4C$. As $g$ is continuous, then $g$ is uniformly continuous on $[0,D]$, and therefore let $\Tilde{\varepsilon}>0$ be small enough so that for $x$, $y \in [0,D]$ with $|x-y|< \Tilde{\varepsilon}$ we have $|g(x)-g(y)|<\varepsilon/2K$.
        \par
        Let $\delta \coloneqq \frac{1}{2} \min\{ \Tilde{\varepsilon}, a_K/2, \varepsilon/4C, D-a_1\}$. Let $\gamma \in \mathcal{I}_H(E)$ such that $d'_H(\beta, \gamma)<\delta$. Consider a correspondence $(\Tilde{U}_j, \Tilde{V}_j)_{1 \leq j \leq \Tilde{K}}$ from $\beta$ to $\gamma$ with distortion less than $\delta$. Since $\delta < a_K$, we must have that $\Tilde{K} \geq K$ and $\{U_j\}_{1 \leq j \leq K} \subseteq \{\Tilde{U}_j\}_{1 \leq j \leq \Tilde{K}}$. Let $(V_j)_{1 \leq j \leq K}$ denote the intervals of $\gamma$ paired in the correspondence with the respective $U_j$, and let $b_j\coloneqq\text{Leb}(V_j)$. 
        \par
        Consider the correspondence of $g^*(\beta)$ and $g^*(\gamma)$ that uses the corresponding interval pairs to $(U_j, V_j)_{1 \leq j \leq K}$. This has distortion
        \begin{equation*}
            \sum_{j=1}^K |g(a_j) - g(b_j)| + \max \Big\{ \norm{g^*(\beta)} - \sum_{j=1}^K g(a_j), \norm{g^*(\gamma)} - \sum_{j=1}^K g(b_j) \Big\}.
        \end{equation*}
        First, let $1 \leq j \leq K$. Then $|a_j - b_j| < \tilde{\varepsilon}$ and therefore we have $|g(a_j) - g(b_j)| < \varepsilon / 2K$. Also note that $\norm{\beta} - \sum_{j=1}^K a_j < \varepsilon / 4C$, and therefore we have
        \begin{equation*}
            \norm{\gamma} - \sum_{j=1}^K b_j \leq d'_H(\beta, \gamma) + \norm{\beta} - \sum_{j=1}^K a_j < \frac{\varepsilon}{2C}.
        \end{equation*}
        From this we can conclude that
        \begin{equation*}
            \max \Big\{ \norm{g^*(\beta)} - \sum_{j=1}^K g(a_j), \norm{g^*(\gamma)} - \sum_{j=1}^K g(b_j) \Big\} < \frac{\varepsilon}{2}
        \end{equation*}
        and the first result follows. To prove the second result, we show that if $\liminf_{x \rightarrow 0}g(x)/x>0$ then
        \begin{equation*}
            \limsup_{x \rightarrow 0} \frac{g^{-1}(x)}{x} \leq \Big(\liminf_{y \rightarrow 0} \frac{g(y)}{y}\Big)^{-1} < \infty
        \end{equation*}
        and the result then follows as previously. To see this let $(x_n, n \geq 1)$ be a positive sequence that converges to $0$. Then $(g^{-1}(x_n), n \geq 1)$ is also a sequence that converges to $0$. Therefore we have that 
        \begin{equation*}
            \liminf_{n \rightarrow \infty} \frac{x_n}{g^{-1}(x_n)} = \liminf_{n \rightarrow \infty} \frac{g(g^{-1}(x_n))}{g^{-1}(x_n)} \geq \liminf_{y \rightarrow 0} \frac{g(y)}{y} >0.
        \end{equation*}
        Therefore $\limsup_{n \rightarrow \infty} g^{-1}(x_n)/x_n \leq 1/(\liminf_{y \rightarrow 0} g(y)/y)<\infty$ and the second result follows in the same way as the first result.
    \end{proof}
\end{lemma}
\begin{remark}
    In fact the conditions in Lemma \ref{lemma g star is continuous} are necessary as well as sufficient in the sense that if $\limsup_{x \rightarrow 0} g(x)/x=\infty$ then the map $g^*$ is not injective and if $\liminf_{x \rightarrow 0} g(x)/x=0$ then the map $(g^{-1})^*$ is not injective (they both map certain interval partitions to $\emptyset$ as the transformed intervals are not summable). 
\end{remark}
We now give certain bounds for the diffusion under Strong Assumption \ref{assumption b} which we will use to prove continuity in the initial state in Proposition \ref{proposition continuity in the initial state}.
\begin{lemma}\label{lemma all under z}
Consider a diffusion $\mathbf{Z}$ that satisfies Strong Assumption \ref{assumption b} with untransformed diffusion $\mathbf{Y}$ and spatial transformation $g$ such that $\liminf_{x \rightarrow 0} g(x)/x>0$, and let $z>0$ and let $\varepsilon >0$. Then there exists $k_7=k_7(\varepsilon, z)>0$ such that for any interval partition $\gamma \in \mathcal{I}_H(g(E))$ such that $\norm{\gamma}<k_7$ we have that 
\begin{equation*}
    \mathbf{P}(\zeta^+(\mathbf{X}_\gamma) <z)>1-\varepsilon,
\end{equation*}
where $\mathbf{X}_\gamma$ is constructed as in Definition \ref{definition starting interval partition} and $\zeta^+(X)\coloneqq\sup_{t \in [0, \text{len}(X)]}X(t)$ for $X \in \mathcal{D}$. 
\begin{proof}
Let $\varepsilon$, $z>0$. The injectivity of $(g^{-1})^*$ from Lemma \ref{lemma g star is continuous} leads to the surjectivity of $g^*$ and so $g^*(\mathcal{I}_H(E))=\mathcal{I}_H(g(E))$. Also we note that the untransformed diffusion satisfies Corollary \ref{corollary assumption b' starting interval partition}. Therefore we can create a scaffolding-and-spindles pair $(\mathbf{N}_\gamma, \mathbf{X}_\gamma)$ with an initial interval partition $\gamma \in \mathcal{I}_H(g(E))$ as in Definition \ref{definition starting interval partition} where we construct the scaffolding by adding $\zeta(f_{V_i})$ jumps for each interval $V_i$ of the initial interval partition $\gamma$. We will compare this to an SPLP with a single initial jump $P=\sum_{U \in \gamma}\zeta(f_U)$ at the beginning of the process which we run until it hits level $0$ and we will denote this process as $\mathbf{X}_P$. By a natural coupling we have $\mathbf{X}_P(t) \geq \mathbf{X}_\gamma(t)$ for all $t \geq 0$ and therefore the event $\{\zeta^+(\mathbf{X}_\gamma)>z\}$ is a subset of by $\{\zeta^+(\mathbf{X}_P)>z\}$. We construct $\mathbf{X}_P$ by considering a L\'evy process started at $\sum_{U \in \beta}\zeta(f_U)$ such that $\mathbf{X}_P-\sum_{U \in \beta}\zeta(f_U)$ has the same distribution as an SPLP $\mathbf{X}$ with L\'evy measure $\nu(\zeta \in \cdot)$. 
\par
Let $W$ be the scale function for the L\'evy process $\mathbf{X}$. This function is continuous and increasing with $W(0)=0$, \cite[Theorem VII.8]{bertoin1996levy}, and therefore there exists $x>0$ such that $W(x)/W(z) \leq \varepsilon/2$. Also note that for the empty interval partition $\emptyset$ we have $d'_H(\gamma, \emptyset) = \norm{\gamma}$. Recall that we can apply Proposition \ref{proposition starting interval partition} to the untransformed diffusion by Corollary \ref{corollary assumption b' starting interval partition}. By Lemma \ref{lemma g star is continuous} $(g^{-1})^*$ is continuous at $\emptyset$ and therefore there exists $k_7>0$ such that for all $\gamma \in \mathcal{I}_H(g(E))$ with $\norm{\gamma}<k_7$ we have $\norm{(g^{-1})^*(\gamma)}<\varepsilon x/2(1+k_3)$ where $k_3$ is constant from Proposition \ref{proposition starting interval partition}. Recall $\varepsilon'$ from the same proposition and take $k_7<g(\varepsilon')$. Then, by \eqref{equation sum of starting diffusion lifetimes}, we therefore have that
\begin{equation*}
    \mathbf{E}\Big(\sum_{i \in I} \zeta(f_{V_i})\Big) \leq (1+k_3) \sum_{i \in I} g^{-1}(\text{Leb}(V_i)) = (1+k_3) \norm{(g^{-1})^*(\gamma)}.
\end{equation*}
Therefore, for $\gamma \in \mathcal{I}_H(g(E))$ such that $||\gamma||<k_7$, we have that\begin{equation*}
    \mathbf{P}\big(\zeta^+(\mathbf{X}_\gamma)>z \big) \leq \mathbf{P}\big(\zeta^+ (\mathbf{X}_P) > z\big) \leq \mathbf{P}\Big(\sum_{i \in I} \zeta(f_{V_i}) > x\Big) + \frac{W(x)}{W(z)} \leq \varepsilon.
\end{equation*}
\end{proof}
\end{lemma}
In the following lemma, for $x \in g(E)$, let $\mathbb{Q}^0_x(\cdot)$ denote the law of the diffusion under Assumption \ref{assumption b}.
\begin{lemma}\label{lemma uniform stopping time}
    Consider a diffusion under Assumption \ref{assumption b}, let $D \in g(E)$ and let $d_0$, $d_1>0$. Then there exists $d_2>0$ such that for all $a$, $b \in [0,D]$ with $a<b<a+d_2$ we have 
    \begin{equation*}
        \mathbb{Q}^0_b(T_a<d_1)>1-d_0.
    \end{equation*}
    \begin{proof}
        This result uses the Markov inequality. We first note that
        \begin{equation*}
            \mathbb{Q}^0_b(T_a>d_1) \leq \frac{\mathbb{E}_{\mathbb{Q}^0_b}(T_a)}{d_1}
        \end{equation*}
        and recall from Corollary \ref{corollary assumption b' starting interval partition} that the function $r(x):=\mathbb{E}_{\mathbb{Q}^0_{g^{-1}(x)}}(\zeta)$ has a derivative that is bounded above near $0$ (we have the transformation $g^{-1}$ here as the corollary we use only applies to untransformed diffusions). Therefore, setting $C:=\sup_{x \in (0,g^{-1}(D))} r'(x)$ we have for all $0 \leq a < b \leq D$  
        \begin{equation*}
            \mathbb{E}_{\mathbb{Q}^0_b}(T_a) = r(g^{-1}(b)) - r(g^{-1}(a)) \leq C|g^{-1}(b) - g^{-1}(a)|.
        \end{equation*}
        Now, $g^{-1}$ is uniformly continuous on $[0,D]$ and so there exists $d_2>0$ such that for $a$, $b \in [0,D]$ with $b>a$ and $b-a<d_2$ we have $g^{-1}(b) - g^{-1}(a)<d_0d_1/C$ and the result follows.
    \end{proof}
\end{lemma}
\begin{lemma}[Continuity in the initial state for $\overline{\text{skewer}}(\mathbf{N}_U, \mathbf{X}_U)$]\label{lemma continuity intial state single interval}
    Let $d_3$, $d_4$, $z>0$, $D \in g(E)$ and let $a \in [0,D]$. Then there exists $d_5>0$ such that for $b\in [0,D]$ with $|a-b|<d_5$ we have the following. Let $U=(0,a)$ and $V=(0,b)$. Then, on some probability space $(\Omega, \mathcal{A}, \mathbf{P})$, there exists a coupled pair of scaffolding and spindles $(\mathbf{N}_U, \mathbf{X}_U)$, $(\mathbf{N}_V, \mathbf{X}_V)$ started from $U$ and $V$ as in Definition \ref{definition starting interval partition} on which we have
    \begin{equation*}
        \mathbf{P}(d'_H(\text{skewer}(z, \mathbf{N}_U, \mathbf{X}_U), \text{skewer}(z, \mathbf{N}_V, \mathbf{X}_V))>d_3)<d_4.
    \end{equation*}
    \begin{proof}
    Let $d_3$, $d_4$, $z>0$ and let $a \leq D \in g(E)$. Let $d_6>0$ be small enough that $W(z-2d_6)/W(z)>1-d_4/5$ where $W$ is the scale function for the SPLP from Lemma \ref{lemma all under z}. Let $d_7=d_2/2$ where $d_2$ is the constant obtained from Lemma \ref{lemma uniform stopping time} when applied with $d_0=d_4/5$ and $d_1=d_6$. Let $f$ be a $0$-diffusion started at $(a+d_7) \wedge D$ and note that
\begin{equation*}
    \mathbf{P}\Big(G_0 \coloneqq \Big\{T_{(a-d_7) \vee 0}(f) < d_6\Big\}\Big)>1-\frac{d_4}{5}
\end{equation*}
where $T_{(a-d_7) \vee 0}(f)$ is the hitting time of $(a-d_7) \vee 0$ by $f$. The diffusion $f$ has a finite lifetime and so is uniformly continuous and therefore there exists $d_8>0$ such that
    \begin{equation*}
        \mathbf{P}\Big(G_1 \coloneqq \Big\{|f(x)-f(y)|<\frac{d_3}{2} \hspace{0.25cm} \forall x,y \in \mathbb{R}:|x-y|<d_8\Big\}\Big)>1-\frac{d_4}{5}.
    \end{equation*}
Now consider an independent Poisson random measure $\mathbf{N}$ with scaffolding $\mathbf{X}$ which is stopped at $T_{-\zeta(f)}$. First note that
\begin{equation*}
    \mathbf{P}\Big(G_2 \coloneqq \Big\{ \sup_{t \in [T_{-\zeta(f)+2d_6},T_{-\zeta(f)}]} \mathbf{X}_t \leq z-\zeta(f)\Big\}\Big) \leq \frac{W(z-2d_6)}{W(z)} > 1 - \frac{d_4}{5}
\end{equation*}
by the independence of $f$ and $\mathbf{N}$ and the strong Markov property of $\mathbf{X}$ and where we follow convention that the supremum over the empty set is $-\infty$. Furthermore note that $(\text{skewer}(y,\mathbf{N}, \mathbf{X}), y \in \mathbb{R})$ has a uniformly continuous modification $(\gamma^y, y \in \mathbb{R})$ by Theorem \ref{theorem continuity}. Let $d_9>0$ be small enough so that
\begin{equation*}
    \mathbf{P}\Big(G_3\coloneqq \Big\{d'_H(\gamma^x, \gamma^y)<d_3/2 \text{   for all   } x,y \in \mathbb{R}:|x-y|<d_9\Big\}\Big)>1-\frac{d_4}{5}.
\end{equation*}
Now, let $d_{10}=d_2/2$ where $d_2$ is now the constant obtained from Lemma \ref{lemma uniform stopping time} when applied with $d_1=d_8 \wedge d_9$ and $d_0=d_4/5$ and let $d_5 \coloneqq d_7 \wedge d_{10}$. Then we have
\begin{equation*}
    \mathbf{P}\Big(G_4 \coloneqq \Big\{T_{(a-d_5) \vee 0}(f) - T_{(a+d_5) \wedge D}(f)<d_8 \wedge d_9\Big\}\Big)>1-\frac{d_4}{5}.
\end{equation*}
Now let $b \in [0,D]$ with $|a-b|<d_5$. Let $f_U(\cdot) \coloneqq f(T_a + \cdot)$ and let $f_{V}(\cdot) \coloneqq f(T_{b}+\cdot)$. We define
\begin{equation*}
    \Tilde{\mathbf{N}}_{U} \coloneqq \mathbf{N}|_{[0,T_{-\zeta(f_{U})}(\mathbf{N})]}; \hspace{0.25cm} \mathbf{N}_{U} \coloneqq \delta_{(0,f_{U})} + \Tilde{\mathbf{N}}_{U}; \hspace{0.25cm} \Tilde{\mathbf{N}}_{V} \coloneqq \mathbf{N}|_{[0,T_{-\zeta(f_{V})}(\mathbf{N})]}; \hspace{0.25cm} \mathbf{N}_{V} \coloneqq \delta_{(0,f_{V})} + \Tilde{\mathbf{N}}_{V}.
\end{equation*}
Note that
\begin{equation*}
    \mathbf{P}\Big(G_5 \coloneqq \Big\{ \gamma^{z-\zeta(f_U)} = \text{skewer}(z-\zeta(f_U), \mathbf{N},\mathbf{X}), \hspace{0.25cm} \gamma^{z-\zeta(f_V)} = \text{skewer}(z-\zeta(f_V), \mathbf{N},\mathbf{X})\Big\}\Big) = 1
\end{equation*}
by the spatial homogeneity of the scaffolding-and-spindles pair and the independence of $f_U$ and $\mathbf{N}$ and of $f_V$ and $\mathbf{N}$. Therefore we have
\begin{equation*}
    \text{skewer}(z, \Tilde{\mathbf{N}}_U, \mathbf{X}_U) = \gamma^{z-\zeta(f_U)}, \hspace{0.25cm} \text{skewer}(z, \Tilde{\mathbf{N}}_V, \mathbf{X}_V) = \gamma^{z-\zeta(f_V)} \text{   on   } G_0 \cap G_2 \cap G_5. 
\end{equation*}
Furthermore we have that
\begin{equation*}
    |f_U(z)-f_V(z)|< d_3/2 \text{   on   } G_1 \cap G_4; \hspace{0.25cm} \text{   and   } \hspace{0.25cm} d'_H(\gamma^{z-\zeta(f_U)}, \gamma^{z-\zeta(f_V)})< d_3/2 \text{   on   } G_3 \cap G_4.  
\end{equation*}
Therefore, as
\begin{equation*}
    d'_H(\text{skewer}(z, \mathbf{N}_U, \mathbf{X}_U), \text{skewer}(z, \mathbf{N}_V, \mathbf{X}_V)) \leq |f_U(z)-f_V(z)| + d'_H(\text{skewer}(z, \Tilde{\mathbf{N}}_U, \mathbf{X}_U), \text{skewer}(z, \Tilde{\mathbf{N}}_V, \mathbf{X}_V))
\end{equation*}
then it follows that
\begin{equation*}
    d'_H(\text{skewer}(z, \mathbf{N}_U, \mathbf{X}_U), \text{skewer}(z, \mathbf{N}_V, \mathbf{X}_V)) < d_3 \text{   on   } G \coloneqq G_0 \cap G_1 \cap G_2 \cap G_3 \cap G_4 \cap G_5,
\end{equation*}
and as $\mathbf{P}(G)>1-d_4$ the result follows.
    \end{proof}
\end{lemma}
\begin{proposition}[Continuity in the initial state]\label{proposition continuity in the initial state}
Consider a diffusion under Strong Assumption \ref{assumption b} for which $\liminf_{x \rightarrow 0} g(x)/x >0$. Let $z>0$ and let $f:\mathcal{I}_H \rightarrow [0,\infty)$ be a bounded and continuous function. Then the map $\gamma \mapsto \mathbb{E}_\gamma(f(\gamma^z))$ is continuous on $(\mathcal{I}_H(g(E)), d'_H)$.
\begin{proof}
This was shown for self-similar diffusions in \cite[Proposition 6.15]{construction}. We adapt the proof to generalise to our setting.
\par
As in Lemma \ref{lemma all under z} we are under the assumption $\liminf_{x \rightarrow 0} g(x)/x >0$ and so $g^*(\mathcal{I}_H(E))=\mathcal{I}_H(g(E))$. Let $\gamma \in \mathcal{I}_H(g(E))$, $\varepsilon$, $z>0$, and let $(\Omega, \mathcal{F}_{\mathcal{I}}, \mathbb{P}_\gamma)$ be the probability space as defined in Section \ref{section simple markov}, where $\Omega \coloneqq \mathcal{C}((0,\infty), \mathcal{I}_H)$, $\mathbb{P}_\gamma$ is the probability measure for $(\gamma^y, y \geq 0)$ and $\mathcal{F}_{\mathcal{I}}$ is the natural filtration for $(\gamma^y, y \geq 0)$. As $f$ is continuous at $\gamma^z$ we can define a $(0,\infty)$-valued random variable $A_\varepsilon(\gamma^z)$ on $(\Omega, \mathcal{F}_{\mathcal{I}}, \mathbb{P}_\gamma)$ for which we have that for $\beta \in \mathcal{I}_H(g(E))$ we have $d'_H(\gamma^z,\beta)<A_\varepsilon(\gamma^z)$ implies that $|f(\gamma^z)-f(\beta)|<\varepsilon/3$. Let $M\coloneqq \sup_{\beta \in \mathcal{I}_H(g(E))} f(\beta)$. Then there is $\Tilde{\varepsilon}>0$ such that $\mathbb{P}_\gamma(A_\varepsilon(\gamma^z)<\Tilde{\varepsilon})<\varepsilon/3M$. We take $\Hat{\gamma} \in \mathcal{I}_H(g(E))$ and consider a probability space $(\Tilde{\Omega}, \mathcal{A}, \mathbf{P})$ for a pair of processes $(\gamma^y, y \geq 0)$ and $(\Hat{\gamma}^y, y \geq 0)$ which have marginal distributions $\mathbb{P}_\gamma$ and $\mathbb{P}_{\Hat{\gamma}}$ respectively; we will later specify an explicit coupling. By Lemma \ref{lemma g star is continuous} note that $g^*$ is surjective so as the untransformed diffusion satisfies Corollary \ref{corollary assumption b' starting interval partition} it is possible to construct IP-evolutions $(\gamma^y, y \geq 0)$ and $(\hat{\gamma}^y, y \geq 0)$ with initial states $\gamma$ and $\hat{\gamma}$. We first note that
\begin{equation*}
    \Big|\mathbb{E}_\gamma(f(\gamma^z)) - \mathbb{E}_{\Hat{\gamma}}(f(\Hat{\gamma}^z))\Big| \leq \mathbf{E}\Big(|f(\gamma^z) - f(\Hat{\gamma}^z)|\Big). 
\end{equation*}
Then we can further bound the expectation above by
\begin{equation*}
    \mathbf{E}\Big(|f(\gamma^z) - f(\Hat{\gamma}^z)|\Big(\mathbbm{1}_{d'_H(\gamma^z,\Hat{\gamma}^z)<\Tilde{\varepsilon}, A_\varepsilon(\gamma^z) \geq \Tilde{\varepsilon}} + \mathbbm{1}_{d'_H(\gamma^z,\Hat{\gamma}^z)\geq\Tilde{\varepsilon}} + \mathbbm{1}_{A_\varepsilon(\gamma^z) < \Tilde{\varepsilon}} \Big)\Big).
\end{equation*}
The first of the three terms is less than or equal to $\varepsilon/3$; the second is less than or equal to $M\mathbf{P}(d'_H(\gamma^z,\Hat{\gamma}^z)\geq\Tilde{\varepsilon})$; and the third is less than or equal to $\varepsilon/3$. Therefore it suffices to show that there is some $\delta>0$ such that for $\Hat{\gamma} \in \mathcal{I}_H(g(E))$ with $d'_H(\gamma,\Hat{\gamma})<\delta$  we have a coupled pair of evolutions $(\gamma^y, y \geq 0)$ and $(\Hat{\gamma}^y, y \geq 0)$ starting from these two initial states, with
\begin{equation}\label{equation coupled probability continuity bound}
    \mathbf{P}(d'_H(\gamma^z, \Hat{\gamma}^z) > \Tilde{\varepsilon}) < \frac{\varepsilon}{3M}.
\end{equation}
To this end, let $U_1$, $U_2$, $\dots$ denote the blocks of $\gamma$, listed in non-increasing order by mass and let $a_j\coloneqq\text{Leb}(U_j)$.
\par
Let $k_7(\varepsilon/9M, z)$ be the constant from Lemma \ref{lemma all under z} and let $K$ be large enough that $\sum_{j>K}a_j < \frac{1}{2}k_7$. If $a_1<g(c)$ then take $D \in g(E)$ such that $a_1<D<g(c)$ and if $a_1=g(c)$ then take $D=a_1=g(c)$. For $1 \leq j \leq K$ let $\delta_j=d_5$ as in Lemma \ref{lemma continuity intial state single interval} with $d_3=\Tilde{\varepsilon}/K$, $d_4 = \varepsilon/9MK$, $D$ as above and $a=a_j$. Now let
\begin{equation*}
    \delta < \left\{
        \begin{matrix}
            \min\{ \delta_1, \dots, \delta_K, a_K, \frac{k_7}{2}, D-a_1 \} \text{   if   } a_1<D<g(c), \\
            \min\{ \delta_1, \dots, \delta_K, a_K, \frac{k_7}{2}\} \text{   if instead   } a_1=D=g(c). \\
        \end{matrix} 
        \right.
\end{equation*}
Let $\Hat{\gamma} \in \mathcal{I}_H(g(E))$ with $d'_H(\gamma,\Hat{\gamma})<\delta$. By definition of $d'_H$, there exists a correspondence $(\Tilde{U}_j, \Tilde{V}_j)_{1 \leq j \leq \Tilde{K}}$ from $\gamma$ to $\Hat{\gamma}$ with distortion less than $\delta$. Since $\delta<a_K$, we get $\Tilde{K} \geq K$ and $\{U_j\}_{1 \leq j \leq K} \subset \{\Tilde{U}_j\}_{1 \leq j \leq \Tilde{K}}$. Let $(V_j)_{1 \leq j \leq K}$ denote the terms paired with the respective $U_j$ in the correspondence. For $1 \leq j \leq K$, let $b_j\coloneqq\text{Leb}(V_j)$. Let $(V_j, j >K)$ be the remaining intervals of $\Hat{\gamma}$. Note that necessarily we have $|a_j-b_j|<\delta$ for $1 \leq j \leq K$. Consider coupled scaffolding and spindles pairs $(\mathbf{N}_{U_j}, \mathbf{X}_{U_j})$, $(\mathbf{N}_{V_j}, \mathbf{X}_{V_j})$ independently for each $1 \leq j \leq K$ as in Lemma \ref{lemma continuity intial state single interval}. For $j>K$ we consider independent scaffolding and spindles $(\mathbf{N}_{U_j}, \mathbf{X}_{U_j})$ and $(\mathbf{N}_{V_j}, \mathbf{X}_{V_j})$. We construct coupled scaffolding and spindles $(\mathbf{N}_\gamma, \mathbf{X}_\gamma)$ and $(\mathbf{N}_{\Hat{\gamma}}, \mathbf{X}_{\Hat{\gamma}})$ constructed from coupled scaffolding and spindles $(\mathbf{N}_{U_j}, \mathbf{X}_{U_j})$ and $(\mathbf{N}_{V_j}, \mathbf{X}_{V_j})$ for $j \in \mathbb{N}$ as in Definition \ref{definition starting interval partition}. Let $(\gamma^y, y \geq 0)$ and $(\Hat{\gamma}^y, y \geq 0)$ be the modifications of $\overline{\text{skewer}}(\mathbf{N}_\gamma, \mathbf{X}_\gamma)$ and $\overline{\text{skewer}}(\mathbf{N}_{\Hat{\gamma}}, \mathbf{X}_{\Hat{\gamma}})$ that are continuous on $(0,\infty)$ which are shown to exist in Lemma \ref{lemma continuity away from zero}. First we make the obvious note that
\begin{equation*}
    \mathbf{P}\Big(F_0 \coloneqq \Big\{ \gamma^z = \text{skewer}(z, \mathbf{N}_\gamma, \mathbf{X}_\gamma); \hspace{0.25cm} \Hat{\gamma}^z = \text{skewer}(z, \mathbf{N}_{\Hat{\gamma}}, \mathbf{X}_{\Hat{\gamma}})\Big\}\Big) = 1
\end{equation*}
as $(\gamma^y, y \geq 0)$ and $(\Hat{\gamma}^y, y \geq 0)$ are modifcations of the skewer processes. Also, for $1 \leq j \leq K$, we have $|a_j - b_j|<\delta_j$ and so
\begin{equation*}
    \mathbf{P}\Big(F_1^j \coloneqq \Big\{ d'_H(\text{skewer}(z, \mathbf{N}_{U_j}, \mathbf{X}_{U_j}), \text{skewer}(z, \mathbf{N}_{V_j}, \mathbf{X}_{V_j}))< \frac{\Tilde{\varepsilon}}{K}\Big\}\Big)> 1 - \frac{\varepsilon}{9MK} \text{   for   } 1 \leq j \leq K.
\end{equation*}
Furthermore we have
\begin{equation*}
    \mathbf{P}\Big(F_2\coloneqq \Big\{ \forall j >K, \zeta^+(\mathbf{N}_{U_j})<z\Big\}\Big) \geq 1 - \frac{\varepsilon}{9M} \text{   by Lemma \ref{lemma all under z}.}
\end{equation*}
Also we have that
\begin{equation*}
    ||\Hat{\gamma}|| - \sum_{j=1}^K b_j \leq d'_H(\gamma, \Hat{\gamma}) + ||\gamma|| - \sum_{j=1}^K a_j < \delta + \frac{1}{2}k_7 < k_7.
\end{equation*}
Therefore 
\begin{equation*}
    \mathbf{P}\Big(F_3 \coloneqq \Big\{\zeta^+(\mathbf{N}_{V_j})<z, \forall j>K \Big\}\Big) \geq 1-\frac{\varepsilon}{9M} \text{   again by Lemma \ref{lemma all under z}.}
\end{equation*}
Pulling all these together we conclude
\begin{equation*}
    \mathbf{P}\Big(F \coloneqq F_0 \cap \Big(\bigcap_{j=1}^K F_1^j\Big) \cap F_2 \cap F_3  \Big) \geq 1 - \frac{\varepsilon}{3M},
\end{equation*}
and we now show that $d'_H(\gamma^z, \Hat{\gamma}^z)<\Tilde{\varepsilon}$ on $F$ to finish the proof. First we note that on $F_2$ $\gamma^z$ consists only of intervals from spindles from $\mathbf{N}_{U_j}$ for $1 \leq j \leq K$ and similarly for $\Hat{\gamma}^z$ on $F_3$ and therefore
\begin{equation*}
    d'_H(\gamma^z, \Hat{\gamma}^z) \leq \sum_{j=1}^K d'_H(\text{skewer}(z, \mathbf{N}_{U_j}, \mathbf{X}_{U_j}), \hspace{0.125cm} \text{skewer}(z, \mathbf{N}_{V_j}, \mathbf{X}_{V_j})) \text{   on   } F_0 \cap F_2 \cap F_3,
\end{equation*}
and so $d'_H(\gamma^z, \Hat{\gamma}^z)\leq \Tilde{\varepsilon}$ on $F$ and the result follows.
\end{proof}
\end{proposition}
\begin{theorem}[Strong Markov property]\label{theorem strong markov property}
    Consider a diffusion under Strong Assumption \ref{assumption b} and then let $\gamma \in g^*(\mathcal{I}_H(E))$. Let $(\mathbf{N}_\gamma, \mathbf{X}_\gamma)$ be a scaffolding-and-spindles pair started from $\gamma$ as in Definition \ref{definition starting interval partition} and let $(\gamma^y, y \geq 0)$ be a modification that is continuous on $(0,\infty)$ as in Lemma \ref{lemma continuity away from zero}. Recall the appropriate filtered probability space a filtered probability space $(\Omega, \mathcal{F}_{\mathcal{I}}, (\mathcal{F}_{\mathcal{I}}^y, y \geq 0), \mathbb{P}_\gamma)$ as used in Corollary \ref{corollary simple markov property}. Then we have that $(\gamma^y, y \geq 0)$ is a strong Markov process on this filtered probability space. 
    \begin{proof}
        We first consider a diffusion that satisfies Strong Assumption \ref{assumption b} and we begin with the case where $\liminf_{x \rightarrow 0} g(x)/x > 0$. As in \cite[Proposition 6.17]{construction} we extend the result in Proposition \ref{proposition continuity in the initial state} to the claim of continuity of the process
        \begin{equation*}
            \gamma \mapsto \mathbb{E}_\gamma\Big(\prod_{i=1}^m f_i(\gamma^{y_i})\Big),
        \end{equation*}
        where $f_1, \dots, f_m: \mathcal{I}_H \rightarrow [0,\infty)$ are bounded and continuous. This can be shown by induction on $m$ using the simple Markov property from Corollary \ref{corollary simple markov property}. Then we show the strong Markov property from this by the method given by Kallenberg \cite[Theorem 17.17]{KallenbergOlav2021Fomp} where the continuity in the initial state stated above then works in place of the Feller property that is used in their result, as in \cite[Proposition 6.17]{construction}.
        \par
        Finally we consider the case where $\limsup_{x \rightarrow 0} g(x)/x < \infty$. We first note that Proposition \ref{proposition continuity in the initial state} applies to the identity transformation $g(x)=x$. Also we have the map $g^*:\mathcal{I}_H(E) \rightarrow g^*(\mathcal{I}_H(E))$ is continuous and bijective by Lemma \ref{lemma g star is continuous} which means that the Strong Markov property in the untransformed case leads to the strong Markov property in the transformed case.
    \end{proof}
\end{theorem}
We can now prove Theorem \ref{theorem diffusion}.
\begin{proof}[Proof of Theorem \ref{theorem diffusion}]
    First consider the case where we have a diffusion under Assumption \ref{assumption b} and we have a scaffolding-and-spindles pair $(\mathbf{N},\mathbf{X})$ for this diffusion stopped at an a.s. finite random time $T$. Then there exists a continuous modification of $\overline{\text{skewer}}(\mathbf{N},\mathbf{X})$ by Theorem \ref{theorem continuity} and the furthermore $\mathbf{X}$ has a local time that is continuous in time and space by Proposition \ref{corollary bi-continuous}.
    \par
    Now, instead consider a diffusion under Assumption \ref{assumption B'} and let $s>0$. Let $(\Tilde{\mathbf{N}},\Tilde{\mathbf{X}})$ be a scaffolding-and-spindles pair for this untransformed diffusion that is stopped at $\tau^0(s)$. Then by Corollary \ref{corollary assumption b' starting interval partition} the diffusion satisfies \eqref{equation starting IP condition}. Therefore, as the diffusion also satisfies $B1$ we can apply Proposition \ref{proposition clade poisson random measure} as $\tau^0(s)$ satisfies the conditions of the proposition. This gives a regular conditional distribution of the clade component for level $0$ given $\beta^0\coloneqq\text{skewer}(0,\Tilde{\mathbf{N}},\Tilde{\mathbf{X}})$ and by Proposition \ref{proposition simple markov skewer} we have a simple Markov property for $\overline{\text{skewer}}(\Tilde{\mathbf{N}},\Tilde{\mathbf{X}})$ on the probability space $(\Omega, \mathcal{A},\mathbf{P})$ for $\Tilde{\mathbf{N}}$. We also have aggregate health summability at level $0$ of the transformed diffusion by Proposition \ref{corollary levy measure assumption b} and Theorem \ref{theorem aggregate mass summability}. Therefore by Theorem \ref{theorem strong markov property} we have a strong Markov property for a modification of the transformed skewer process on the space of paths continuous on $(0,\infty)$ which we denote as $(\Omega, \mathcal{F},\mathbb{P}_\mu)$ where $\mu$ is the probability distribution on $g^*(\mathcal{I}_H(E))$ of $g^*(\beta^0)$.
\end{proof}

%% file: appendix.tex
\section{Appendix}\label{appendix}
\begin{lemma}\label{lemma existence disintegration}
Let $\nu$ be a $\sigma$-finite measure on $\mathcal{E}$. Then, for any measurable map $S: \mathcal{E} \rightarrow [0,\infty]$, there exists a $S$-disintegration $(\nu_t)_{t \in [0,\infty]}$. This disintegration is uniquely determined up to an almost sure equivalence, i.e. if $(\nu_t^*)_{t \in [0,\infty]}$ is another $S$-disintegration of $\nu$ then $\nu (S^{-1}(t \in [0,\infty] : \nu_t \neq \nu_t^*)) = 0$. 
\end{lemma}
Lemma \ref{lemma green differential} restates Lemma \ref{lemma solution green differential}. The Green's function for the $\uparrow$-diffusion is given by
\begin{equation*}
        G^\uparrow_{(a,w)}(x,v) = 2 \frac{s^\uparrow(x)-s^\uparrow(a)}{s^\uparrow(w)-s^\uparrow(a)} \left(s^\uparrow(w)-s^\uparrow(v)\right)m^\uparrow(v) \mathbbm{1}_{x<v<w} + 2 \frac{s^\uparrow(w)-s^\uparrow(x)}{s^\uparrow(w)-s^\uparrow(a)} \left(s^\uparrow(v)-s^\uparrow(a)\right)m^\uparrow(v) \mathbbm{1}_{a<v<x},
\end{equation*}
Recall the functions $U_{(a,w)}$ and $V_{(a,w)}$ from \eqref{equation U differential solution}:
\begin{equation*}
    U_{(a,w)}(x) = \mathbb{E}_{\mathbb{Q}^\uparrow_x}\big(\big(T_a \wedge T_w\big)^2\big), \hspace{0.25cm} V_{(a,w)}(x) = 2 \mathbb{E}_{\mathbb{Q}^\uparrow_x}\big(T_a \wedge T_w\big).
\end{equation*}
We have from Section \ref{subsection Green's functions} that for $0<a<w<c$ that $U_{(a,w)}$ is the solution to the following differential equation, from \eqref{equation U differential}, in $x \in (a,w)$
\begin{equation*}
    \frac{1}{2} \sigma^2(x) U''_{(a,w)}(x) + \mu^\uparrow(x) U'_{(a,w)}(x) + V_{(a,w)}(x) = 0,
\end{equation*}
with boundary conditions $U_{(a,w)}(a)=U_{(a,w)}(w)=0$. 
\begin{lemma}[Proof of Lemma \ref{lemma solution green differential}]\label{lemma green differential}
We claim that the solution to the differential equation in \eqref{equation U differential}, with boundary conditions given by $U(a)=U(w)=0$, is the following:
\begin{equation*}
    U_{(a,w)}(x) \coloneqq  \mathbb{E}_{\mathbb{Q}^\uparrow_x}\big(\big(T_a \wedge T_w\big)^2\big) = k(a,w) (s^\uparrow (w) - s^\uparrow(x)) + \int_{y=x}^w (s^\uparrow)'(y) \int_{z=a}^y \frac{2V_{(a,w)}(z)}{(s^\uparrow)'(z)\sigma^2(z)}dzdy,
\end{equation*}
where $k(a,w)$ is a constant given by
\begin{equation*}
    k(a,w) = \frac{-1}{s^\uparrow(w) - s^\uparrow(a)} \int_{y=a}^w (s^\uparrow)'(y) \int_{z=a}^y \frac{2V_{(a,w)}(z)}{(s^\uparrow)'(z)\sigma^2(z)}dzdy,
\end{equation*}
and where $V_{(a,w)}(x)$ is 
\begin{equation*}
    V_{(a,w)}(x) = 2 \mathbb{E}_{\mathbb{Q}^\uparrow_x}\big(T_a \wedge T_w\big) = 2 \int_{v=0}^w G^\uparrow_{(a,w)}(x,v)dv,
\end{equation*}
where $G^\uparrow$ is the Green's function for the $\uparrow$-diffusion. 
\begin{proof}
We first note that the $\uparrow$-diffusion satisfies $A3$ and that the infinitesimal drift $\mu^\uparrow$ and infinitesimal variance $(\sigma^2)^\uparrow$ are continuous functions on $(0,c)$. Therefore the differential equation in \eqref{equation U differential}, from \cite[Equation 15.3.37]{karlin1981second}, has solution given in \eqref{equation U differential solution}.
\par
We now check that the boundary conditions are satisfied:
\begin{equation*}
    U_{(a,w)}(w) = k(a,w) (s^\uparrow (w) - s^\uparrow(w)) + \int_{y=w}^w (s^\uparrow)'(y) \int_{z=a}^y \frac{2V_{(a,w)}(z)}{(s^\uparrow)'(z)\sigma^2(z)}dzdy=0,
\end{equation*}
\begin{equation*}
    \begin{split}
        U_{(a,w)}(a) = & \frac{-1}{s^\uparrow(w) - s^\uparrow(a)} \int_{y=a}^w (s^\uparrow)'(y) \int_{z=a}^y \frac{2V_{(a,w)}(z)}{(s^\uparrow)'(z)\sigma^2(z)}dzdy (s^\uparrow (w) - s^\uparrow(a)) \\
        & + \int_{y=a}^w (s^\uparrow)'(y) \int_{z=a}^y \frac{2V_{(a,w)}(z)}{(s^\uparrow)'(z)\sigma^2(z)}dzdy = 0. \\
    \end{split}
\end{equation*}
Now we show that the differential equation itself is solved by this solution. Taking a derivative with respect to $x$ gives
\begin{equation*}
    U'_{(a,w)}(x) = - k(a,w) (s^\uparrow)'(x) - (s^\uparrow)'(x) \int_{z=0}^x \frac{2V_{(a,w)}(z)}{(s^\uparrow)'(z)\sigma^2(z)} dz.
\end{equation*}
Then multiplying by $1/(s^\uparrow)'(x) = \exp\left(\int_{t=b}^x \frac{2\mu^\uparrow(t)}{\sigma^2(t)} dt\right)$ gives
\begin{equation*}
    \exp\left(\int_{t=1}^x \frac{2\mu^\uparrow(t)}{\sigma^2(t)} dt\right) U'_{(a,w)}(x) = \frac{1}{(s^\uparrow)'(x)} U'_{(a,w)}(x) = - k(a,w) - \int_{z=0}^x \frac{2V_{(a,w)}(z)}{(s^\uparrow)'(z) \sigma^2(z)} dz.
\end{equation*}
Finally taking another a derivative with respect to $x$ gives
\begin{equation*}
    \exp\left(\int_{t=b}^x \frac{2\mu^\uparrow(t)}{\sigma^2(t)} dt\right) \left(U''_{(a,w)}(x) + \frac{2\mu^\uparrow(x)}{\sigma^2(x)}U'_{(a,w)}(x)\right) = \frac{-2V_{(a,w)}(x)}{(s^\uparrow)'(x)\sigma^2(x)},
\end{equation*}
and multiplying by $(s^\uparrow)'(x) = \exp\left(\int_{t=b}^x \frac{-2\mu^\uparrow(t)}{\sigma^2(t)} dt\right)$, adding $\frac{2V_{(a,w)}(x)}{\sigma^2(x)}$, then multiplying by $\frac{1}{2}\sigma^2(x)$ gives the differential equation in the form of \eqref{equation U differential}. Therefore the result follows by standard results for the uniqueness of solutions to differential equations.
\end{proof}
\end{lemma}
This final result generalises \cite[Lemma 8]{supplementA} to our setting. The other relevant results from this supplementary paper apply without any modification.
\begin{lemma}
    For fixed $y \in \mathbb{R}$, the map $N \mapsto \text{skewer}(y,N, \xi_N)$ is measurable on the set $\mathcal{N}^{\text{sp}}_{\text{fin}}$. Recall that we map $(y,N, \xi_N)$ to the empty interval partition $\emptyset$ in the event that the aggregate mass is infinite; i.e. when $M^y_{N,\xi(N)}(\infty)=\infty$, see Definition \ref{definition aggregate mass}. 
    \begin{proof}
        The proof \cite[Lemma 8]{supplementA} that the map restricted to the finite aggregate mass summability case is measurable generalises to our setting because the $\sigma$-algebra of $(\mathcal{I}_H, d'_H)$ is the same as $(\mathcal{I}_H, d_H)$ (\cite[Theorem 2.3(b)]{metric2020diversity}), where for $\beta$ and $\gamma \in \mathcal{I}_H$ we define $d_H(\beta, \gamma)$ to be the Hausdorff distance between $[0,\norm{\beta}]\backslash \beta$ and $[0,\norm{\gamma}]\backslash \gamma$. The extension to the map to include the case where the aggregate mass is infinite follows as the subset $\mathcal{N}^{\text{sp}}_{\text{fin}}$ for which the aggregate mass at level $y$ is finite is measurable. 
    \end{proof}
\end{lemma}
\subsection{Excursion intervals}
This section follows \cite[A Excursion intervals]{construction} very closely. All their measure theoretic results hold in our setting, and we state the results and state why they generalise to our case.
\begin{proposition}\label{proposition properties of levy process of unbounded variation}
Let $\mathbf{X}$ be an SPLP with jump measure that satisfies \eqref{equation levy measure unbounded variation}. For $y \in \mathbb{R}$ recall the concept of the sets of excursion intervals $V^y(\mathbf{X})$ and $V^y_0(\mathbf{X})$ from Definition \ref{definition set of excursion intervals}. Let $N \in \mathcal{N}^{\text{sp}}$ and $y \in \mathbb{R}$. Then we have the following a.s.:
\begin{enumerate}
    \item $V^y(\mathbf{X})=\Big\{[a,b] \subset [0,\infty)| a<b; \mathbf{X}_{a-}=y=\mathbf{X}_b; \text{   and   } \mathbf{X}_t \neq 0 \text{   for   } t \in (a,b)\Big\}$.

    \item For $I$, $J \in V^y_0(\mathbf{X})$, $I \neq J$, the set $I \cap J$ is either empty or a single shared endpoint.
    \item If two interval $[a,b]$, $[b,c] \in V^y_0(\mathbf{X})$ share an endpoint $b$ then $\mathbf{X}$ does not have a jump at $b$.
    \item For every $t \notin \bigcup_{I \in V^y_0(\mathbf{X})} I$, we have that $\mathbf{X}_{t-}=\mathbf{X}_t=y$.
    \item $\text{Leb}\Big([0,\infty) \backslash \bigcup_{I \in V^y_0(\mathbf{X})}I\Big) = 0$.
    \item There are no degenerate excursions of $\mathbf{X}$ about level $y$, that is to say that for all excursions we have $0<T^+_0<\zeta$ a.s.
    \item Local times $(\ell^y(t), t \geq 0)$ of $\mathbf{X}$ exist and are injective with respect to the excursion intervals in the sense that, for $[a,b]$, $[c,d] \in V^y_0$, $\ell^y(a) \neq \ell^y(c)$ unless $[a,b]=[c,d]$.
    \item If $y>0$, we have $T^y > T^{\geq y} \coloneqq \inf\{t \geq 0: \mathbf{X}_t \geq y \}$.
    \item For $[a,b] \in V^y(N)$, the process $N|^\leftarrow_{I^y_N(a,b)}$ is a bi-clade. Furthermore, the set $\{N|^\leftarrow_{I^y_N(a,b)}: [a,b] \in V^y_0(N)\}$ partitions the spindles of $N$.
\end{enumerate}
\begin{proof}
    These results are given in the self-similar cases in \cite[Proposition A.2-4]{construction}. The proofs in \cite[Proposition A.2]{construction} immediately generalise. As for the results of Proposition \cite[Proposition A.3]{construction}, they rely on the exit properties of L\'evy processes characterised by Millar \cite{millar1973exit} which applies to our general setting. The final part, from \cite[Proposition A.4]{construction}, is a result of the careful constructions of the intervals $I^y_N(a,b)$ and also applies in our general setting.
\end{proof}
\end{proposition}
We conclude this section with definitions, specifically concerning incomplete excursions, to formally define the cutoff processes stated in \eqref{equation cutoff one} and \eqref{equation cutoff two}. We first define $T^{\geq y}$ and $T^{\geq y}_*$ in the same way as we have defined $T^y$ and $T^y_*$ and then with these define $N^{\leq y}_{\text{first}}$, $N^{\leq y}_{\text{last}}$, $N^{\geq y}_{\text{first}}$, $N^{\geq y}_{\text{last}}$:
\begin{equation*}
    T^{\geq y} \coloneqq \inf ( \{t \in [0,\text{len}(N)]: \xi_N(t) \geq y\} \cup \{\text{len}(N)\}); \hspace{0.25cm} T^{\geq y}_* \coloneqq \sup(\{t \in [0,\text{len}(N)]: \xi_N(t-) \leq y \} \cup \{0\});
\end{equation*}
\begin{equation*}
    \begin{split}
        & N^{\leq y}_{\text{first}} = \Big(N|_{[0,T^{\geq y})} + \mathbbm{1}_{\xi_N(T^{\geq y}-)<y}\delta_{(T^{\geq y}, \check{f}^y_{T^{\geq y}})} \Big) \mathbbm{1}_{T^{\geq y} \neq T^{\geq y}_*};  N^{\geq y}_{\text{last}} = N|^\leftarrow_{(T_*^{\geq y}, \text{len}(N)]} + \mathbbm{1}_{\xi_N(T_*^{\geq y}-)>y}\delta_{(0, \hat{f}^y_{T_*^{\geq y}})};\\
        & N^{\leq y}_{\text{last}} = N|^\leftarrow_{[T^y_*,T^{\geq y}_*)} + \mathbbm{1}_{y \neq \xi_N(\text{len}(N)); \xi_N(T_*^{\geq y}-)<(y \wedge \xi_N(T^{\geq y}_*))}\delta_{(T_*^{\geq y} - T^y_*, \check{f}^y_{T_*^{\geq y}})}; \\
        & N^{\geq y}_{\text{first}} = \Big(N|^\leftarrow_{(T^{\geq y},T^y]} + \mathbbm{1}_{y \neq 0;\xi_N(T^{\geq y}-)>y \vee 0}\delta_{(0, \hat{f}^y_{T^{\geq y}})} \Big) \mathbbm{1}_{T^{\geq y} \neq T^{\geq y}_*}. \\
    \end{split}
\end{equation*}